\pdfoutput=1
\RequirePackage{ifpdf}
\ifpdf % We~are running pdfTeX in pdf mode
\documentclass[pdftex]{sigma}
\else
\documentclass{sigma}
\fi

\numberwithin{equation}{section}

\usepackage{tabularx,mdwtab}

\usepackage{bbm}
\newcommand{\One}{\mathbbmss{1}}

\newcommand{\Mod} {{\sf Mod}}

\newcommand {\fbgl}{{\mathfrak{bgl}}}

\newcommand {\fel}{{\mathfrak{el}}}

\newcommand {\fwk}{{\mathfrak{wk}}}

\newcommand {\See}{{\mathbb{S}}}

\newcommand{\un}{\underline{N}}
\newcommand {\Cee}{{\mathbb C}}
\newcommand {\Eee}{{\mathbb E}}
\newcommand {\Fee}{{\mathbb F}}
\newcommand{\Kee}{{\mathbb K}}
\newcommand{\Zee}{{\mathbb Z}}
\newcommand{\fc}{\mathfrak{c}}
\newcommand{\fe}{\mathfrak{e}}
\newcommand{\fg}{\mathfrak{g}}
\newcommand {\ev}{{\bar0}}
\newcommand {\od} {{\bar1}}
\newcommand {\bcdot}{\mathbin{\hbox{\raise.4ex\hbox{\bf.}}}}
\newcommand {\tto}{\longrightarrow}

\newtheorem{Theorem}{Theorem}[section]
\newtheorem*{Theorem*}{Theorem}

\newtheorem{Lemma}[Theorem]{Lemma}

\newtheorem{Conjecture}[Theorem]{Conjecture}
 { \theoremstyle{definition}

\newtheorem{Remark}[Theorem]{Remark}
\newtheorem{Remarks}[Theorem]{Remarks}
\newtheorem*{OpenProblem}{Open Problem}
\newtheorem*{OpenProblems}{Open Problems}
\newtheorem*{Comment}{Comment}
\newtheorem*{Hypothesis}{Hypothesis}
 }

\begin{document}

\allowdisplaybreaks

\newcommand{\arXivNumber}{0807.3054}

\renewcommand{\PaperNumber}{031}

\FirstPageHeading

\ShortArticleName{Deformations of Symmetric Simple Modular Lie (Super)Algebras}

\ArticleName{Deformations of Symmetric Simple Modular\\ Lie (Super)Algebras}

\Author{Sofiane BOUARROUDJ~$^{\rm a}$, Pavel GROZMAN~$^{\rm b}$ and Dimitry LEITES~$^{\rm ac}$}

\AuthorNameForHeading{S.~Bouarroudj, P.~Grozman and D.~Leites}

\Address{$^{\rm a)}$~New York University Abu Dhabi, Division of Science and Mathematics,\\
\hphantom{$^{\rm a)}$}~P.O.~Box 129188, United Arab Emirates}
\EmailD{\href{mailto:sofiane.bouarroudj@nyu.edu}{sofiane.bouarroudj@nyu.edu}, \href{mailto:dimleites@gmail.com}{dimleites@gmail.com}}

\Address{$^{\rm b)}$~Deceased}

\Address{$^{\rm c)}$~Department of Mathematics, University of Stockholm, SE-106 91 Stockholm, Sweden}

\ArticleDates{Received November 16, 2016, in final form February 02, 2023; Published online May 29, 2023}

\Abstract{We say that a~Lie (super)algebra is ``symmetric'' if with every root (with respect to the maximal torus) it has the opposite root of the same multiplicity. Over algebraically closed fields of positive characteristics (up to 7 or 11, enough to formulate a general conjecture), we computed the cohomology corresponding to the infinitesimal deformations of all known simple finite-dimensional symmetric Lie (super)algebras of rank $<9$, except for superizations of the Lie algebras with ADE root systems, and queerified Lie algebras, considered only partly. The moduli of deformations of any Lie superalgebra constitute a supervariety. Any infinitesimal deformation given by any odd cocycle
is integrable. All deformations corresponding to odd cocycles are new. Among new results are classifications of the cocycles describing deforms (results of deformations) of the 29-dimensional Brown algebra in characteristic~3, of Weisfeiler--Kac algebras and orthogonal Lie algebras without Cartan matrix in characteristic~2. Open problems: describe non-isomorphic deforms and equivalence classes of cohomology theories. Appendix: For several modular analogs of complex simple Lie algebras, and simple Lie algebras indigenous to characteristics 3 and 2, we describe the space of cohomology with trivial coefficients. We show that the natural multiplication in this space is very complicated.}


\Classification{17B50; 17B55; 17B56; 17B20}

\setcounter{tocdepth}{2}

{\small \tableofcontents}

\subsection*{Summary} A long introduction, split into three parts to make it more digestible, contains a review of the conjectural classification of simple finite-dimensional modular Lie algebras and Lie superalgebras over (algebraically closed) fields $\Kee$ of characteristic $p>0$ and various partial results, some of which, although declared done long ago, are still incomplete, see formula~\eqref{HamSt}.

Section~\ref{Scoh} contains a discussion of the cohomology theory in general and in particular its part needed to classify deformations, derivations and central extensions of Lie algebras and superalgebras over $\Kee$ for $p>0$; with peculiarities indigenous to $p=3$ and $2$. Here, we formulate several vital problems, regrettably still open.

Next sections contain our results: classification of cocycles describing deformations of Lie algebras and Lie superalgebras considered in characteristics sufficient to formulate a general conjecture speaking of serial cases; to prove the claim analytically is an open problem. To distinguish the isomorphism classes of the deformed algebras is also an open problem. It is clear that one should use the method applied by Chebochko (consider orbits of the algebraic (super)groups naturally acting in the spaces of cohomology); this is a long work to be done.

Several cases remain to be computed, e.g., deformations of superizations in characteristic~2 of the orthogonal series (in particular, the periplectic superalgebras, cf.\ with~\cite{FSSc}) and of the $\fe$-type superalgebras.

\section{Introduction (with ``super'' in the background)}\label{Sintro}

Hereafter, $\Kee$ is an algebraically closed field of characteristic
$p>0$ and $\fg$ is a~finite-dimensional Lie (super)algebra;
$\Zee_{+}$ is the set of non-negative integers, the elements of $\Zee/p$ are denoted with a~bar over the non-negative integer representatives of the class of an element, e.g., $\Zee/2=\{\ev,\od\}$. The parity function and characteristic are denoted by the same letter $p$; no misunderstanding can occur.

Recall that the Lie algebras \textit{over fields of
positive characteristic} are often called ``modular''. This term is not related with
modular elements in lattices and appeared much earlier than \textit{modular
elements} became a~term in programming. Although overused now, it is shorter
than the italicized attribution; and (a referee reminds us) it became widespread after the book~\cite{Sel} with the term in the title was published: ``the use of ``modular'' is quite common in the representation theory of finite groups, where it means ``representations modulo $p$'' as opposed the ``ordinary'' representations; and its usage goes back to at least 1941''.

This paper precedes \cite{BGLL1} and
both of them are sequels to \cite{BGL2}, where we classified finite-dimensional modular Lie
superalgebras of the form $\fg(A)$ with indecomposable symmetrizable
Cartan matrix $A$ over $\Kee$.

The paper \cite{BGL2} contains
also precise definitions of notions such as \textit{root}, \textit{Cartan
matrix}, and \textit{Dynkin diagram} over $\Kee$, where they sometimes differ
from their namesakes over~$\Cee$. Although this is not used in this text, we recommend also the recent papers~\cite{BLLoS, KLS} containing additional arguments in favor of the suggested definition of the root.

In this paper we continue taking care of terminology, make definitions of certain notions---for example, of Lie superalgebra in terms of the functor of points---precise, and continue to get results needed to understand deformations parametrized by supervarieties---an attempt to superize a part of Skryabin's pape~\cite{SkR}.
\begin{gather*}%\label{Hereafter}
\text{Hereafter, the unspecified $i$ in $\fg(A)^{(i)}/\fc$ is the smallest such that $\fg(A)^{(i)}=\fg(A)^{(i+1)}$},
\end{gather*}
usually stabilization occurs at $i=1$ or $2$. In particular, for $A$ non-invertible, \cite{BGL2} lists the simple Lie (super)algebras
$\fg(A)^{(i)}/\fc$, where $\fc$ is the center of $\fg(A)$ and $\fg(A)^{(i)}$. If $A$ is invertible, then~$\fg(A)$ is simple.

We overview several aspects of deformation theory.
Rudakov, Kostrikin, Dzhumadildaev, Kuznetsov and his students (notably, Chebochko
and Ladilova)
are main contributors to the \textit{classification} of \textit{deforms}, i.e., the results of deformations. (The term \textit{deform} was suggested by
M.~Gerstenhaber in analogy with \textit{transform}, the result of
transformation, or philosopher's term \textit{construct}, something constructed by the
mind.)

Here we consider the Lie
(super)algebras listed in \cite{BGL2} and several other
``\textit{symmetric}'' Lie (super)algebras, the ones that with every
root (with respect to the maximal torus) contain its opposite of the same multiplicity; otherwise the Lie (super)algebra is said to be ``non-symmetric''.\footnote{Can an algebra be symmetric w.r.t.\ one maximal torus and non-symmetric w.r.t.\ another? It is an interesting question, but, probably, better not go into this now.}

Deformations of non-symmetric Lie (super)algebras will be considered in detail elsewhere, here we only
give a~sweeping overview needed to formulate a~sketch of the super version of the Kostrikin--Shafarevich conjecture; for most unexpected results, see \cite{BGLLS, GZ, KrLe,SkT1}.

Assuming our target audience contains not only experts in superalgebras, we recall what the Lie algebra (co)homology is---if considered not as a~derived functor, but as something that can be explained to the computer. In doing so we realized, with the help of a~wonderfully helpful and inquisitive referee and A.~Lebedev, that the task is more difficult than we thought it was; as a result, we formulate desirable claims as open problems.

\subsection{Basic definitions. Subtleties}\label{ssBasDef} Recall the sign rule of linear algebra (for \textit{neighboring} elements):

\medskip

\begin{minipage}[l]{14cm}
``If something of parity $a$ is moved past something of parity $b$ the sign $(-1)^{ab}$ accrues. Formulas defined only on homogeneous elements are extended to arbitrary elements via linearity.''
\end{minipage}

\medskip

Studying deformations, derivations, and central extensions of Lie \textit{super}algebras we encounter several new phenomena; blind application of the sign rule to the non-super definitions might cause trouble.

\subsection{Super commutativities: anti-, skew-, antiskew-, and ``just so''}\label{ssAnti}

Recall, see \cite{Ls, Gr}, that whereas in pre-super era the adjectives anti-commuting and skew-commuting meant the same thing, in the super setting we MUST introduce a~new terminology and use different terms to denote different notions in the following cases where, for all elements~$a$,~$b$ of a~superalgebra $A$, we have
\begin{alignat*}{4}% \label{supersym_rule}
& ba=(-1)^{p(b)p(a)}ab \quad &&(\text{super symmetry}),\quad&& (s)&\\
& ba=-(-1)^{p(b)p(a)}ab \quad &&(\text{super anti-symmetry}),\quad &&(as)&\\
& ba=(-1)^{(p(b)+1)(p(a)+1)}ab \quad &&(\text{super skew-symmetry}),\quad && (ss)\\
& ba=-(-1)^{(p(b)+1)(p(a)+1)}ab \quad &&(\text{super antiskew-symmetry}).\quad && (ass)
 \end{alignat*}
In other words: our ``anti'' designates the overall minus sign, whereas ``skewness'' can be straightened (made symmetric) by reversing the parity of $A$.

 In the case of characteristic~2, the signs disappear, so \textit{super commutativity} and \textit{super antiskew-commutativity} turn into
\[
\text{$ab=ba$ for all $a$, $b$ and $a^2=0$ for $p(a)=\od$,}
\]
whereas \textit{super anti-commutativity} and \textit{super skew-commutativity} turn into
\[
\text{$ab=ba$ for all $a$, $b$ and $a^2=0$ for $p(a)=\ev$.}
\]

\textit{In this paper, all commutative algebras and supercommutative superalgebras are supposed to be associative with~$1$; their morphisms should send~$1$ to~$1$, and the morphisms of supercommutative superalgebras should preserve parity}.

Let $p\neq 2$. Let $C$ be a~supercommutative superalgebra and $V$ a~$C$-module. The \textit{tensor superalgebra} is $T^{\bcdot}_C(V):=\oplus_{n\geq 0}T^{n}_C(V)$, where $T^{0}_C(V)=C$ and $T^{n}_C(V)=V\otimes_C\dots\otimes_C V$ with $n>0$ factors. The
\textit{super symmetric} algebra of~the $C$-module $V$ is the $C$-algebra
\begin{equation*}%\label{1Req18}
S^{\bcdot}_C(V):=T^{\bcdot}_C(V)/I_s,
\end{equation*}
where the two-sided
ideal $I_s$ is generated by elements of~the form (we skip the $\otimes$ between factors)
\begin{equation*}% \label{xy_px_py}
 xy-(-1)^{p(x)p(y)}yx \qquad \text{for any $x, y\in V = T^1_C(V)$.}
\end{equation*}
If $p=2$, we need to add (as the careful referee noted) to generators of $I_s$ those of the form
\[
\text{$x^2$ for any $x\in V_\od = T^1_C(V)_\od$.}
\]

In the algebra $S^{\bcdot}_C(V)$, we consider the grading
\begin{equation*}% \label{grading_1}
 S^{\bcdot}_C(V) =
 \mathop{\oplus}\limits_{n\geq 0}S^n_C(V),\qquad \text{where $S^n_C(V) = T^n_C(V) \pmod{I_s}$.}
\end{equation*}
Since $I_s\subset \oplus _{n\geq 2}T^n_C(V)$, we have $S^1_C(V) =
T^1_C(V) = V$. Denote by
\begin{equation*} %\label{can_hom_1}
 i\colon \ V\tto S^1_C(V)\subset S^{\bcdot}_C(V)
\end{equation*}
the corresponding embedding. Clearly, the
elements of~$V$ generate $S^{\bcdot}_C(V)$, which is a~supercommutative $C$-algebra.

\subsection{The super anti-symmetric vs. super exterior algebra of a module}\label{ss1.9.4}
Define the \textit{super exterior
algebra} $E^{\bcdot}_C(V):=S^{\bcdot}_C(\Pi(V))$, where $\Pi$ is the change of parity functor. If $V$ is a~superspace over a~field $\Fee$ or over $\Zee$, it does not matter from which side we apply~$\Pi$ to~$V$, i.e., $\Pi(V):=\Pi(\Zee)\otimes_\Zee V$ and $V\Pi:=V\otimes_\Zee \Pi(\Zee)$
are isomorphic $\Fee$-modules. For any one-sided $C$-module $V$ over any supercommutative superalgebra $C$ with $C_\od\neq 0$, there are two inequivalent ways to make $V$ into a~two-sided $C$-module, see \cite{Ls}: the conventional one we use here, namely
\[
cv=(-1)^{p(c)p(v)}vc \qquad \text{for any $c\in C$ and $v\in V$,}
\]
or by setting $cv=(-1)^{p(c)(p(v)+1)}vc$. Over $C$ with $C_\od\neq 0$, it is important to distinguish $C$-module $\Pi(V)$ from $V\Pi$, and the left dual ${}^*V$ from the right dual $V^*$. The package \textsc{SuperLie}, see~\cite{Gr}, distinguishes all these cases.

The algebra $E^{\bcdot}_C(V)$
is a~supercommutative superalgebra. The composition of~the
$C$-module morphisms $\Pi({\rm id})\colon V\tto\Pi(V)$ and
$i\colon \Pi(V)\tto S^{1}_C(\Pi(V))$ determines the {\it
canonical odd homomorphism} $\iota\colon V\tto
E^{1}_C(V)$.

Observe that if $p=2$ and $V$ is purely even, then
\begin{equation}\label{EinS}
 E^{\bcdot}_\Kee(V)\subset S^{\bcdot}_\Kee(V),
\end{equation}
where
\begin{gather*}
E^{\bcdot}_\Kee(V):=T^{\bcdot}_\Kee(V)/(x\otimes x \ \text{for any $x\in V$}),\\
S^{\bcdot}_\Kee(V):=T^{\bcdot}_\Kee(V)/(x\otimes y+ y\otimes x \ \text{for any $x, y\in V$})
\end{gather*}
are the conventional exterior algebra and symmetric algebra, respectively. It should be noted in~\eqref{EinS} that $E^{\bcdot}_\Kee(V)$ embeds in $S^{\bcdot}_\Kee(V)$ as a subspace, but not as a subalgebra: the algebra $E^{\bcdot}_\Kee(V)$ contains nilpotent elements, while $S^{\bcdot}_\Kee(V)$ is an integral domain. Since $E^{\bcdot}_\Kee(V)$ is commutative when $p=2$, it is actually a~quotient algebra of the free commutative algebra $S^{\bcdot}_\Kee(V)$.
This fact hints at the possible existence---only if $p=2$---of a~cohomology theory with symmetric, rather than anti-symmetric, product of cochains. Such cohomology did appear, naturally, several times, see \cite{DzZ, LZ, Z} and references therein.

The super-exterior and super-symmetric powers
of a~superspace (we skip the subscript $\Kee$) are designated
\[
\Eee^k(V):=\See^k(\Pi(V))=\mathop{\oplus}\limits_{0\leq i\leq k}S^i(V_\od)\otimes E^{k-i}(V_\ev).
\]

\subsubsection{Another superization of the exterior algebra}\label{ssExtLam} Apply the sign rule to the
anti-symmetry condition, i.e., let $I_a$ be the two-sided ideal of~$T_C^{\bcdot}(V)$ generated by
\begin{gather*}% \label{xy_px_py_2}
 xy + (-1)^{p(x)p(y)} yx\text{ for
 any $x, y\in V= T^1_C(V)$},\\
\text{and $x^2$ for any $x\in V_\ev$ if $p=2$}.
\end{gather*}
Let
\begin{equation*}%\label{Lambda_C_1}
 \Lambda_C^{\bcdot}(V) := T_C^{\bcdot}(V)/I_a.
\end{equation*}

Hereafter we identify $\ev$ with 0 and $\od$ with 1 in similar formulas where we add elements of $\Zee/2$ with elements of $\Zee$ thus making expressions \eqref{ab_1_pa_pb}, \eqref{a_dega_pa} meaningful.

The superalgebra $\Lambda_C(V)$ is not
supercommutative. However, it is \textit{graded-commutative} (a.k.a. ``colored") with
respect to $(\Zee/2\, \times \Zee)$-grading: the $\Zee/2$-degree
is parity, the $\Zee$-degree is induced by the standard grading on $T_C^{\bcdot}(V)$, i.e., for any $a, b\in\Lambda_C^{\bcdot}(V)$, we have
\begin{equation}
 \label{ab_1_pa_pb}
 ab=(-1)^{\deg(a)\deg(b)+p(a)p(b)}ba.
\end{equation}
Clearly, $\Lambda^n_C(V)\simeq E^n_C(V)$ as spaces, hence,
\begin{equation}\label{Nekl}
\begin{split}
&\text{$\Lambda^{\bcdot}(V)\simeq E^{\bcdot}(V)$ as spaces; moreover,}\\
& \text{$\Lambda^{\bcdot}(V)$ and $E^{\bcdot}(V)$ are ``equivalent'' algebras
 in the following sense: }
\end{split}
\end{equation}
by passing to the total degree (considered, together with $\deg(a)$, as an element of $\Zee/2$)
\begin{equation}
 \label{a_dega_pa}
 \deg_{\rm tot}(a):=\deg(a)+p(a),
\end{equation}
and slightly modifying the multiplication in $\Lambda^{\bcdot}(V)$ by setting
\begin{equation}\label{a_prod_b_pa}
 a*b:=(-1)^{p(a)\deg(b)}ab,
\end{equation}
we see that the superalgebra obtained is the supercommutative
superalgebra $E^{\bcdot}(V)$:
\begin{equation*} %\label{supercommutative_1}
 a*b=(-1)^{ \deg_{\rm tot}(a)\deg_{\rm tot} (b)}b*a.
\end{equation*}

The fact \eqref{Nekl} was known to experts in deformation theory since at least 1950s. Deligne gave categorical (in both senses) arguments exorcising ``colored algebras'', see \cite[pp.~62--64]{Del}. Molotkov elaborated the fact \eqref{Nekl} and generalized the first theorems of Nekludova (for the commutative and associative case) and Scheunert (for the Lie case, see \cite{Sch}), see Molotkov's Appendix on ``colored algebras'', see~\cite{Ls}, containing a~proof of the following theorem.

%\sssbegin[(There are only two classes of graded algebras)]
\begin{Theorem}
For any finitely generated commutative group $G$, any $G$-graded commutative and associative $($resp.\ Lie$)$ algebra is naturally equivalent $($under analogs of changes \eqref{a_dega_pa} and~\eqref{a_prod_b_pa}, an equivalence defined in~{\rm \cite{Sch})} to either commutative \textup{($G=\{1\}$)} and associative $($resp.\ Lie$)$ algebra, or a~super commutative \textup{($G=\Zee/2$)} and associative $($resp.\ super Lie$)$ algebra.
\end{Theorem}

\subsection{Lie superalgebra, pre-Lie superalgebra, Leibniz superalgebra}\label{sring}
 Generalization of notions of this subsection to super rings and modules over them is immediate.

For any $p$, a~\textit{Lie superalgebra} is a
superspace $\fg=\fg_\ev\oplus\fg_\od$ such that the even part
$\fg_\ev$ is a~Lie algebra, the odd part $\fg_\od$ is a
$\fg_\ev$-module (made into the two-sided one by
\textit{anti}-symmetry, i.e., $[y,x]=-[x,y]$ for any $x\in
\fg_\ev$ and $y\in
\fg_\od$), and a~\textit{squaring} $x\mapsto x^2$ and the \textit{bracket} are defined on $\fg_\od$ via a~linear map $s\colon S_2(\fg_\od)\tto\fg_\ev$, where $S_2(\fg_\od)$ is the subspace of $\fg_\od\otimes\fg_\od$ spanned by elements of the forms $x\otimes x$ and $x\otimes y + y\otimes x$ for any $x,y\in\fg_\od$:
\begin{gather*}%\label{squaring0}
x^2 := s(x\otimes x), \\ %\label{squaring1}\\
{}[x,y] := s(x\otimes y + y\otimes x).%\label{squaring2}
\end{gather*}

The linearity of the $\fg_\ev$-valued function $s$ implies that
\begin{gather*}%\label{squaring3}
\text{$(ax)^2=a^2x^2$ for any $x\in
\fg_\od$ and $a\in \Kee$, and} \\ %\label{squaring4}\\
{}[\cdot,\cdot]\text{~is a~bilinear form on $\fg_\od$ with values
in $\fg_\ev$.} %\label{squaring5}
\end{gather*}

For $p\neq 2,3$, the anti-commutativity and Jacobi identities required for Lie algebras turn under superization into
\begin{gather*}
[x,y]=-(-1)^{p(x)p(y)}[y,x] \quad \text{for any $x,y\in\fg$},\\
[x, [y,z]]=[[x,y],z]+(-1)^{p(x)p(y)}[y,[x,z]]\quad \text{for any $x, y,z\in\fg$}. % and $p\neq 2,3$}.
\end{gather*}

For $p\neq 2,3$, the Jacobi identity for three equal to each other odd elements is
\[
3[x,[x,x]]=0,
\]
which, if $p\neq 3$, is equivalent to
\begin{gather}\label{p=3}
[x,[x,x]]=0 \quad \text{for any $x\in\fg_\od$}.
\end{gather}

If $p=3$, we have to add condition~\eqref{p=3}, separately, to the anti-symmetry and Jacobi identity amended by the sign rule.

If $p=2$, the \textit{anti-symmetry} condition for $p=2$ is
\[
[x,x]=0\quad \text{for any $x\in\fg_\ev$}.
\]

The \textit{Jacobi identity} involving odd elements takes the form of the following two conditions:
\begin{gather}
\big[x^2,y\big]=[x,[x,y]]\quad \text{for any~} x\in\fg_\od, y\in\fg_\ev,\nonumber\\ %\label{6b}\\
\big[x^2,x\big]=0\quad \text{ for any $x\in\fg_\od$.}\label{6}
\end{gather}

The superalgebra satisfying the Jacobi identity, and without any restriction on symmetry is called a~\textit{Leibniz superalgebra}, introduced by Loday, see \cite{BeMh, IKV}.

Over $\Zee/2$, the condition \eqref{6} must (as explained in \cite{KLLS}) be replaced with a~more general one:
\begin{gather}\label{JIgen}
\big[x^2,y\big]=[x,[x,y]] \quad \text{for any~} x,y\in\fg_\od.
\end{gather}
For any other ground field, this condition is equivalent to condition \eqref{6}.

Since, for any Lie superalgebra $\fg$, we want $\operatorname{\mathfrak{der}} \fg$ to be a~Lie superalgebra, we require the following generalization of the condition~\eqref{JIgen}
\[
D\big(x^{2}\big) = [Dx,x]\quad \text{for any $x\in\fg_\od$ and $D\in \operatorname{\mathfrak{der}} \fg$.}
\]

For $p=2$, the definition of the derived algebra of the Lie superalgebra $\fg$ changes: set $\fg^{(0)}:=\fg$ and
\begin{gather*}%\label{deralg}
\fg^{(i+1)}=\big(\fg^{(i)}\big)^{(1)}:=\big[\fg^{(i)},\fg^{(i)}\big]+\operatorname{Span}\big\{g^2\mid
g\in\big(\fg^{(i)}\big)_\od\big\} \quad \text{for any $i\geq 0$}.
\end{gather*}

A \textit{representation} of the Lie superalgebra $\fg$ in the \textit{module} (superspace) $V$ is any Lie superalgebra morphism $\rho\colon \fg\tto \mathfrak{gl}(V)$ such that $\rho([x,y])=[\rho(x), \rho(y)]$ for any $x,y\in\fg$, and $\rho\big(x^2\big)=(\rho(x))^2$ for any $x \in\fg_\od$.

Naively, the Lie superalgebra $\fg$ is \textit{simple}, if $\dim\fg>1$ and $\fg$ has no proper ideals.

\subsubsection{Pre-Lie superalgebras. An open problem}\label{Kann}

 The superalgebra satisfying the anti-symmetry and the Jacobi identity, but not the condition~\eqref{p=3} is called \textit{pre-Lie superalgebra}. We know only two examples of pre-Lie superalgebras, see \cite{BBH}.

The method used in \cite{Kan} to get several exceptional simple Lie superalgebras with Cartan matrix in characteristic~$3$ and $5$ from simple Lie algebras over $\Cee$ yields pre-Lie superalgebras when applied to $\mathfrak{br}(3)$ and its deforms listed in Section~\ref{LeBr3a}.

\begin{OpenProblem}
 Describe the quotients of these pre-Lie superalgebras modulo center.
\end{OpenProblem}

\subsection{The Kostrikin--Shafarevich conjecture}\label{ssKSh} In the 1960s, as a~step towards
classification of simple finite-dimensional Lie algebras over~$\Kee$, Shafarevich, together with his student, Kostrikin, first
considered the \textit{restricted} simple modular Lie algebras. Their
self-restriction was, perhaps, occasioned by the fact that only
restricted Lie algebras correspond to algebraic groups, see~\cite{Vi4}.

The restricted case is much easier than the general one. Being a~geometer, albeit an algebraic one, Deligne
recently gave us advice to look, if $p>0$, at the restricted Lie
(super)algebras and modules over them before venturing into the wilderness of nonrestricted
ones because (in our words) ``the problem of classifying
nonrestricted simple Lie algebras and their representations is
definitely very tough (and, perhaps, is not even reasonable) whereas
restricted algebras are related to geometry", see Deligne's appendix
to~\cite{LL} for his own words and several interesting problems.

Whatever the reasons for preferring restricted Lie algebras to nonrestricted ones, one sometimes \textit{has} to consider the nonrestricted Lie
algebras as well, e.g., to describe the restricted Lie algebras, as in \cite{BW}, or because the results of \textit{Cartan--Tanaka--Shchepohkina prolongs} (briefly: CTS-prolongs), see \cite{Shch}, are more natural than their restricted versions (although these prolongs are mostly non-restricted), or---the referees remind us---in the study of other topics, e.g., in the study of $p$-groups, see~\cite{Kos1}.

In the late 1960s, for restricted simple finite-dimensional Lie algebras over $\Kee$ for $p>7$, Kostrikin and Shafarevich
formulated a~conjectural method producing all such Lie algebras; they also gave
an explicit list of $\Zee$-graded simple algebras. The conjecture turned out to hold for $p>5$ as well, and even for $p=5$ if one adds Melikyan's
example to the list. Having incorporated the deformations, one got a conjectural list of all, not only restricted, simple finite-dimensional Lie algebras over $\Kee$ for $p>3$, now proved, see~\cite{BGP, S}.

\begin{itemize}\itemsep=0pt
\item The KSh-method of getting $\Zee$-graded simple finite-dimensional Lie algebras is as follows. Take simple Lie algebras over $\Cee$ of one of the two types
\begin{gather}\label{KSh1}
\begin{split}
& \text{(a) finite-dimensional, i.e., of the form $\fg(A)$ with Cartan matrix $A$},\\
& \text{(b) vectorial, i.e., of vector fields with polynomial coefficients.}
\end{split}
\end{gather}
\end{itemize}

In case (a), select an integer form $\fg(A)_\Zee$ of $\fg(A)$ (with smallest structure constants, i.e., take the standard Chevalley basis) and take $\fg(A)_\Zee\otimes_\Zee\Kee$. This algebra is simple (perhaps, modulo center). The same method is applicable to Lie superalgebras $\fg(A)$ over $\Cee$; for their list, see~\cite{CCLL}. Observe that the classification of finite-dimensional modular Lie superalgebras $\fg(A)$
with indecomposable $A$ and their simple subquotients is performed in \cite{BGL2} by a~different method (induction on rank).

In case (b), let us immediately superize this part of the KSh-method of producing examples, replace polynomial coefficients by divided powers in $a$ indeterminates, $m$ of which are even and $a-m$ of which are odd. The powers of the even indeterminates are bounded by the shearing vector $\un = (N_1,\dots, N_m)$, forming the supercommutative superalgebra (here $p^{\infty}:=\infty$)
\begin{equation*}%\footnotesize
\label{u;N}
\mathcal{O}(m; \un|n):=\mathbb{K}[u;
\un]:=\operatorname{Span}_{\mathbb{K}}\left(u^{(\underline{r})}\mid r_i
\begin{cases}< p^{N_{i}}&\text{for $i\leq m$}\\
=0\text{ or 1}&\text{for $m<i\leq a$}\end{cases}\right).
\end{equation*}
Introduce
\textit{distinguished} partial derivatives each of them serving as several
partial derivatives at once, for each of the generators $u,u^{(p)}, u^{(p^2)}, \dots$:
\[
\partial_i\big(u_j^{(k)}\big):=\delta_{ij}u_j^{(k-1)}\quad \text{for all $k$.}
\]
\begin{itemize}\itemsep=0pt
\item The (\textit{general}) \textit{Lie superalgebra of vector fields} is defined to be
\[
\mathfrak{vect}(m;\un|n):=\left\{\sum f_i\partial_i\mid f_i\in \mathcal{O}(m; \un|n)\right\}.
\]
The Lie algebra $\mathfrak{vect}(m;\un)$ is usually denoted by $W(m;\un)$ and called \textit{Jacobson--Witt} Lie algebra.
\end{itemize}

The other finite-dimensional $\Zee$-graded simple vectorial Lie algebras for $p>5$ are the \textit{derived} (first or second) of the three classical series of its subalgebras
\begin{itemize}\itemsep=0pt
\item $\mathfrak{svect}(m; \un|n)$ a.k.a.\ $S(m; \un|n)$ of divergence-free vector fields,

\item $\mathfrak{h}(2m; \un|n)$ a.k.a.\ $H(2m; \un|n)$ of Hamiltonian vector fields,

\item $\mathfrak{k}(2m+1; \un|n)$ a.k.a.\ $K(2m+1; \un|n)$ of contact vector fields.
\end{itemize}

\subsection{Too many deformations of Lie (super)algebras}\label{ssTooManyDefs}
The classification of deforms over $\Kee$ is aggravated by the
following phenomenon not as widely known as it deserves; until
recently only a~few examples could be found in the literature.

Sometimes, the deformed algebra is isomorphic to the
initial one even though the cocycle corresponding to the linear part
of the global deformation represents a~nontrivial cohomology class
of $H^2(\fg; \fg)$!
We called such deforms \textit{semi-trivial};
for examples of semi-trivial deforms when $\operatorname{char}\Kee=0$, see~\cite{Ca, Ri}. For a~characterization of a~wide range of such
 deforms when $\operatorname{char}\Kee>0$, see \cite{BLLS, BLW}. We say that any nontrivial and not
semi-trivial deform is a~\textit{true deform}. The Lie (super)algebra is said to be \textit{rigid} if it has no true deformations.

Rudakov proved that for $p>3$ the simple Lie algebras of the form $\fg(A)$ and $\mathfrak{psl}$ are rigid, see~\cite{Ru}. So this part of the KSh-list is clear.

\subsection{The vague part of the KSh-list}\label{ssKShVague} It consists of the deforms of the simple $\Zee$-graded vectorial Lie algebras.
Observe that the classification of simple \textit{restricted} Lie algebras obtained in \cite{BW} (the final list) \textit{does not contain all explicit expressions of $p$-structures to this day}: e.g., for the deforms of
 $\mathfrak{svect}(m; \One|n)$, the description was only recently obtained, see \cite{BKLLS}. For several families of Hamiltonian algebras their $p$-structures remain to be described:
 \begin{equation}\label{HamSt}
\begin{split}
& \text{``The problem of restrictedness is approached. \dots\ [But] the family}\\
& \text{of Hamiltonian algebras \dots\ is not yet handable'', see \cite[p.~357]{S}}.
\end{split}
\end{equation}
Moreover, even the brackets of the deformed Hamiltonian algebras are explicitly given only for $p>3$ and only in the written in Russian PhD~Thesis of S.~Kirillov \cite{Kir}, so are practically inaccessible.

\subsection[Kac's contribution to the KSh-conjecture. Quantization. Shchepochkina's example]{Kac's contribution to the KSh-conjecture. Quantization.\\ Shchepochkina's example}\label{KacFil}

Important sharpening of the vague part of the initial KSh-list of non-restricted simple Lie algebras---description of deforms of simple $\Zee$-graded Lie algebras $\fg$ and their non-simple rela\-ti\-ves---is due to Kac \cite{KfiD}, who suggested to evade computing $H^2(\fg;\fg)$ by considering much easier to compute filtered deformations of $\fg$, corresponding to $H^2(\fg_-;\fg)$, where $\fg_-=\oplus_{i<0} \fg_i$. Moreover, \textit{it seemed} that Kac reduced the latter task to description of the normal shapes of differential forms the deformed algebra preserves.

In the case of Lie algebras of Hamiltonian series, however, Kac's evasion contains a~dearth: these Lie algebras have deformations which are not filtered deformation, namely, \textit{the ones induced by quantizations of the Poisson Lie algebras}. This deviation from Kac's reduction does not contribute to new simple Lie algebras (at least, for $p>3$), and this is, probably, why it did not draw any attention, until a~new Lie algebra, not a~filtered deform of a~volume-preserving Lie algebra, was discovered for $p=2$, see Shchepochkina's example in \cite{BGLLS2}.

For $p>3$, there are two types of preserved differential forms: the volume forms and the symplectic forms, see Section~\ref{sssNF}; and the contact forms, preserved up to a factor which is a~function.

For $p=3$, in addition to the Lie algebras preserving the volume and symplectic forms, there are Skryabin algebras preserving non-integrable distributions more general than the contact distribution, see \cite{GL3}; they require extra study.

For $p=2$, the desuperization of one of Shchepochkina's exceptional simple vectorial Lie superalgebras is a~deform of the divergence-free Lie superalgebra $\mathfrak{svect}(5;\un)$, see \cite{BGLLS2}, not a~filtered deform of $\mathfrak{svect}(5;\un)$; so Kac's arguments, see~\cite{KfiD}, do not work in this case, same as in the Hamiltonian case for any $p>0$.

\subsection{Normal shapes of volume, symplectic and contact forms}\label{sssNF}

Tyurin \cite{Tyu} described normal shapes of volume forms (with an isomorphism missed) for any $p$; Wilson \cite{W} rediscovered the result in a~correct form.

For $p>2$, Skryabin \cite{Sk0, SkH} found the normal shapes of symplectic forms (a result rather difficult to prove) and contact forms. For $p=2$, there are two types of symplectic forms, see~\cite{LeP}, Skryabin found the normal shapes of alternating forms. The case of non-alternating symplectic forms is in the process of being solved, see \cite{KoKCh}. These results are obtained for symplectic forms on varieties; to superize them for symplectic and periplectic forms on supervarieties is an open problem.\looseness=-1\index{Problem, open}

\subsection[Generalization of the Kostrikin--Shafarevich conjecture to $p=3$ and 2.\\ Improvements of the KSh-method]{Generalization of the Kostrikin--Shafarevich conjecture to $\boldsymbol{p=3}$ and 2.\\ Improvements of the KSh-method}\label{ssKShImr}

Based on the results of 40 years of work by several (teams of) researchers, Premet and Strade proved the KSh-conjecture, and completed the classification of all simple finite-dimensional Lie algebras over $\Kee$ for $p>3$,
see \cite{BGP, S}.
The initial KSh-method of getting the complete list of simple finite-dimensional Lie algebras over $\Kee$ was meantime improved.
\begin{itemize}\itemsep=0pt
\item In \cite{KD} (written mainly by Dzhumadildaev
after his teacher and co-author, Kostrikin, died), we find the following
refinement of the KSh-method producing simple Lie algebras: take ``standard'' Lie algebras and their true deforms. The list of ``standard'' Lie algebras is the same list of $\Zee$-graded simple Lie algebras as in the original KSh-method plus Poisson algebras (which can be considered as divergence-free subalgebra of the Lie algebra of contact vector fields) whose deforms are (if $p=5$) Melikyan algebras.

\item The next refinement of the KSh-method (for details, see \cite{BGLLS, BGLLS2, BGL2, Lt}) is
\begin{equation}\label{KShCTS}
\begin{split}
& \text{\emph{Consider as
``standard'' only Lie algebras with indecomposable Cartan matrix}}\\
& \text{\emph{and the $($complete or partial$)$ Cartan--Tanaka--Shchepochkina prolongs}}\\
& \text{\emph{$($CTS-prolongs for short$)$ $($for its most lucid definition and examples,}} \\
& \text{\emph{see {\rm \cite{BGLLS, Shch})} of their non-positive parts in one of the $\Zee$-gradings.}}
\end{split}
\end{equation}
We thus get
the ``classical'' types, and
examples indigenous for $p=3$, such as algebras discovered by Frank, Ermolaev, and, most of all, Skryabin. (Together with one new algebra, Skryabin algebras are more
lucidly than previously described and interpreted in~\cite{GL3}.)
Ladilova described filtered deforms of several Skryabin algebras, see \cite{LaD,LaF, LaY,LaZ}.
\end{itemize}

Let us summarize: the latest simplification of the KSh-method \eqref{KShCTS} emphasizes the role of
\begin{equation*}%\label{mainIng}
\begin{minipage}[c]{14cm}
(i) Lie algebras of the
form $\fg(A)$ with indecomposable Cartan matrix $A$,\\
(ii) the CTS
prolongs, especially \textit{partial} ones,\\
(iii) deformations. \end{minipage}
\end{equation*}
The role of deformations is huge for any $p$, as is clear from \cite{SkH, W}, and
the ``standard examples" are not only simple algebras of types (a) and (b) in~\eqref{KSh1} but also central extensions of simple algebras and
semisimple algebras, as shown, e.g., in \cite{GZ, KD,SkT1} and non-simple subalgebras of complete CTS-prolongs, see \cite{SkH, W}. The refinement \eqref{KShCTS} encourages us to offer the following conjecture.

\begin{Conjecture}[for $p\geq 3$] \label{DKSh}\quad
\begin{enumerate}\itemsep=0pt
\item[$(A)$] Any simple Lie algebra is one of the following:
\begin{enumerate}\itemsep=0pt
\item[$(a)$] of the form $\fg(A)^{(i)}/\fc$, in particular $\fg(A)$ for $A$ invertible,
\item[$(b)$] the $($derived algebra of$)$ $\Zee$-graded vectorial\,\footnote{Often called
by an ill-chosen term ``algebras of Cartan type''. The term is sometimes
applied to simple derived algebra of the Lie algebra of one of the 4 the classical series only, ignoring the
other (relatives of) simple vectorial Lie algebras in
characteristics 2, 3 and 5. There is no
conventional term yet for \textit{true deforms} of vectorial
algebras.} Lie algebras,
\item[$(c)$] a~true deform of an example of type $(a)$ or $(b)$.
\end{enumerate}
\item[$(B)$] Each simple $\Zee$-graded vectorial algebra, i.e., of type $(b)$ of $(A)$, is $($the derived algebra of$)$ the generalized CTS-prolong $($complete or partial$)$ of the non-positive part of $\fg(A)$ in one of its $\Zee$-gradings, see~{\rm \cite{GL3}}.
\end{enumerate}
\end{Conjecture}

\textbf{Towards the list of ``standard objects'' for $\boldsymbol{p=2}$}.
The tempting conjecture of \cite{KD} (``all simple algebras are
deforms of the ``standard" objects'') might turn out to be true if we incorporate
\textit{prolongs}, see \cite{Shch}, found after~\cite{KD} was
written; for $p=2$, see \cite{BGLL, BGLLS,
BGLLS2,Ei, ILL,LeP, SkT1}. For all these cases, there remained the problem of
description of their deforms and isomorphism classes of simple
algebras. The
paper \cite{KD} contains a~proof\footnote{This proof was written
somewhat sloppily (for example, the Poisson Lie algebra was called
Hamiltonian, etc.) but is correct and complete. Still, in
\cite{MeZu} it is double-checked in the particular case: for the
smallest coordinates of the shearing vector.} of the fact that Melikyan algebra
is a~deform of the Poisson Lie algebra.\footnote{An aside: in
\cite{BLLS}, it is shown that various Poisson Lie algebras $\mathfrak{po}(n;\un)$
with the same $|\un|$ are deforms of each other, and the same is
true for their quotients modulo center, Lie algebras of Hamiltonian
vector fields.} A~similar claim concerns several exceptional Lie
algebras for $p=3$, namely Skryabin and Ermolaev\footnote{They were not
specified in \cite{KD}, regrettably: there are several types of
algebras discovered by Skryabin and several ones discovered by
Ermolaev. The fact that a~certain Ermolaev algebra is a~deform of
(the simple derived algebra of) a~Lie algebra of contact type was
announced (without details) in \cite{KuJa} about a~decade earlier
than the same claim appeared in \cite{KD} with even less details. Elsewhere we
started to compute complete list of all deformations of the
appropriate subalgebras of the restricted Lie algebra of contact vector fields.}
ones).

\subsection{Digression: Linguistic. Simplicity vs. ``nature"}\label{parComm}

It is time to reconsider the
terminology designed ad hoc and without much care. Clearly, the
nomen\-cla\-ture should be improved, e.g., Shen indicates
reasons for being unhappy with the current terminology, see \cite{She2}. In
particular, although simple Lie (super)algebras are the first to
classify and study, several of their ``relatives" (the ones with Cartan
matrix, CTS prolongs of simple algebras and their relatives, the algebras of all derivations, central extensions of simple Lie algebras, restricted closures, etc.)
are often no less needed.

In textbooks on geometry, representation theory, and mathematical
physics a~``classical" Lie group or its Lie algebra is, usually, the
one of type GL, SL, O, Sp over $\Cee$
\textit{or} its ``most common" real form. Therefore, certain real
forms are disqualified from the fuzzy set of ``classical'' objects
and can be only found in the complete list of real forms \cite{He, OV}, such
serial Lie (or Chevalley, see~\cite{Go}) groups over other fields are not referred
to as ``classical" objects. Observe, that simplicity is a~neither necessary nor
sufficient condition for a~group (Lie algebra) to be ``classical'';
the finite groups of type PGL, PSL, and
GL are naturally (and justly, we think) called
``classical''. The 5 simple groups (and their Lie algebras) that do
not belong to any series are unequivocally called ``exceptional" meaning ``not forming infinite series'', but definitely are also ``classical''.

Rudakov's definition of Lie algebras of ``classical type'', see~\cite{Ru}, is too restrictive (so much so that we will not recall it): for $p=2$, certain simple Lie algebras with
Cartan matrix might be not of ``classical type'' in Rudakov's terms.

In many papers and in the books \cite{S,S-II,S-III} devoted to Lie algebras
over (usually, algebraically closed) fields of positive
characteristic ``\textit{classical}'' is either
\begin{itemize}\itemsep=0pt
\item any Lie algebra $\fg(A)$ with
indecomposable Cartan matrix $A$ denoted \`a la Cartan by its root
system (be it serial $A_n$, $B_n$, $C_n$, $D_n$ or exceptional
$E_n$, $F_4$ or $G_2$) \textit{or}
\item the simple quotient of the first
or second derived algebra of $\fg=\fg(A)$ modulo center (Chebochko designated this
simple subquotient $\fg^{(i)}/\fc$ by barring its relative with
Cartan matrix: $\overline{A_n}$, etc.).
\end{itemize}

However, the simple Lie algebras with Cartan matrix discovered
by Brown for $p=3$, and Weisfeiler--Kac algebras $\fwk$ for $p=2$, as well as both versions of
$\mathfrak{o}$, i.e., $\mathfrak{o}_I$ and $\mathfrak{o}_\Pi$, for $p=2$, and all deforms of these algebras, are not
counted as ``classical" at the moment. We suggest to make the nomenclature more precise. Namely, as follows.
\begin{itemize}\itemsep=0pt
\item The distinction
between the Lie algebra with Cartan matrix and its derived algebra or simple
subquotient (such as $\mathfrak{gl}(pn)$, $\mathfrak{sl}(pn)$, $\mathfrak{p}\mathfrak{gl}(pn)$, and
$\mathfrak{psl}(pn)$) was often disregarded, but not by Chebochko, see the
third paragraph of \eqref{NC}. Still, like many, she implicitly implies that
$\mathfrak{sl}(pn)$ has a~Cartan matrix, whereas it does not;
it is $\mathfrak{gl}(pn)$ that has a~Cartan matrix, see~\cite{BGL2}.

\item If $p=2$, there are two nonisomorphic types of orthogonal Lie
algebras, see \cite{LeP}:
\begin{itemize}\itemsep=0pt
\item $\mathfrak{o}_\Pi(n)$, which preserves the non-degenerate invariant
symmetric bilinear form (NIS for short) with Gram matrix $\Pi_n=\operatorname{antidiag}(1,\dots, 1)$, and

\item $\mathfrak{o}_I(2n)$, which preserves the non-degenerate invariant symmetric bilinear form with Gram matrix $1_{2n}$.
\end{itemize}

\item It is the derived
(and/or centrally extended) version of $\mathfrak{o}_\Pi(n)$ that has a~Cartan
matrix; no
``relative'' of $\mathfrak{o}_I(2n)$ has Cartan matrix, see \cite{BGL2}.
\end{itemize}

We suggest to consider the following Lie algebras and their
simple and nonsimple relatives as \textit{classical}:
\begin{equation}\label{propusk}
\begin{split}
& \text{serial:}\\
 & \text{for $p>0$:} \ \mathfrak{gl}(pn), \ \mathfrak{pgl}(pn), \ \mathfrak{sl}(pn), \ \mathfrak{psl}(pn) \ \text{and similarly}\\
& \text{any other Lie algebra of the form $\fg(A)$ and its simple subquotient if $A$ is singular}, \\
& \text{for $p=2$:} \ \mathfrak{o}_I(2n), \ \mathfrak{o}_\Pi(2n), \ \widehat{\mathfrak{o}\fc(i;2n)},\ \text{see \cite[formula~(6.13)]{BGL2}}, \ \mathfrak{o}(2n+1),\\
& \text{exceptional:} \\
& \text{for $p=3$:} \ \text{Brown algebras $\mathfrak{br}(2;\varepsilon)$ and $\mathfrak{br}(3)$},\\
& \text{for $p=2$:} \ \text{Weisfeiler--Kac algebras $\fwk(3; a)$ and $\fwk(4; a)$}.
 \end{split}
\end{equation}
Here are a~few samples of the numerous examples which \textit{at the moment} are \textit{non-classical}:
\begin{equation*}%\label{non-class}
\begin{split}
& \text{serial:} \\
& \text{for any $p>0$:} \ \text{(filtered) deforms of vectorial algebras, see \cite{BGLLS2}},\\
& \text{Kaplansky algebras}, \ \text{see \cite{BLLS}},\\
& \text{exceptional:}\\
& \text{for $p=3$:} \ \text{the deform $L(2,2)$, see \cite{BLW}},\\
& \text{for $p=2$:} \ \text{Eick algebras \cite{Ei}, Skryabin algebras \cite{GZ, SkT1}, and Shen algebras \cite{Fei,Shch,She2}},\\
& \text{analogs of Sergeev's $\mathfrak{as}$, i.e., numerous central extensions of various types of}\\
& \text{orthogonal algebras, see \cite{BGLL1}}.
\end{split}
\end{equation*}

Observe that certain \textit{vectorial Lie algebras}, together with their deforms,
not only the direct analogs of
the four series of simple complex Lie algebras, but also those that exist only
in characteristics 2, 3 and 5, are no less ``classical'' than
those of the form~$\fg(A)$ or ``relatives'' thereof.

We suggest to simplify the terminology and consider as ``basic''
objects not only the ones with Cartan matrix, but also $\mathfrak{o}_I(2n)$,
and the CTS prolongs, whereas various relatives of such Lie
(super)algebras, e.g., simple derived algebra ($\fg^{(i)}$ or
$\fg^{(i)}/\fc$) thereof, algebras of derivations, etc. should be
considered as secondary (``derived", speaking figuratively). The term ``exceptional" should be called
to its duty, even if $p>0$, to denote deviations from infinite series as the
term was originally used by Killing and Cartan or---perhaps temporarily, in the absence of classification---``weird'' Kaplansky series (for an attempt to describe them as deforms of something more ``classical'', see~\cite{BLLS}).

\subsection{Explicit cocycles vs. ecology}\label{ssEcology}

Having written the first 15
pages of this paper we were appalled by the volume
the explicit cocycles occupied. Who needs them?! Let us give just
the dimensions and save paper! But just the dimensions, as in Theorem~\ref{NC2},
are difficult to \textit{use}. Further on in her papers, Chebochko
gave a~rather explicit description of the cocycles representing the
cohomology classes that span the spaces listed in Theorem~\ref{NC2}, see \cite{Ch1, ChKu,KuKCh2}.

We are happy to be able to use electronic publications that make it possible to
explicitly give the corresponding cocycles while saving trees: only
the possibility to look at the explicit cocycles enabled us to
interpret some of the mysterious (for almost 20 years) Shen's
``variations", and check if the local deformation is integrable and
compute the global deform. Possibility to look at
explicit cocycles helped us to interpret several ``non-symmetric'' Lie
algebras, see \cite{BLLS}.

However, certain cocycles are too long, and probably will never be
used explicitly by humans. The reader can obtain them by means of
the \textsc{SuperLie} code on one's own.

\subsection{Chebochko's computations of deformations}\label{ssCheb}

In several papers
\cite{ChG2,Ch1, Ch2, KKCh, KuCh,KChA5g}, with difficult results clearly
explained, N.~Chebochko, and her scientific advisor, M.~Kuznetsov,
gave an overview of the situation for Lie algebras of the form~$\fg(A)$ and $\fg(A)^{(i)}/\fc$. In \cite{Ch1}, Chebochko writes (footnotes are ours):
\begin{equation}\label{NC}
\text{\begin{minipage}[c]{14cm} ``According to \cite{Dz}, for $p=3$,
the Lie algebra $C_2$ is the only algebra among the series $A_n$,
$B_n$, $C_n$, $D_n$ that admits nontrivial deformations. In~\cite{Ru}, it was established that over a~field of characteristic
$p>3$ all the \textit{classical}\footnote{On terminology, see the beginning of Section~\ref{parComm} of this paper.} Lie
superalgebras are rigid. In~\cite{KuCh} and~\cite{KKCh}, a~new scheme
was proposed for studying rigidity, and it was proved that the
classical Lie algebras of all types over a~field of characteristic
$p>2$ are rigid, except for the Lie algebra of type $C_2$ for $p=3$.
\\
\phantom{brr} For $p=2$, some deformations of the Lie algebra of
type\footnote{This is, actually, $\mathfrak{psl}(4)$.} $G_2$ were constructed in \cite{She1}. \dots
\\
\phantom{brr} The Lie algebras of type $A_l$ for $l+1\equiv 0\pmod
2$, $D_l$ and $E_7$, have nontrivial centers. We shall say that the
corresponding quotient algebras are of type $\overline{A_l}$,
$\overline{D_l}$ and $\overline{E_7}$, respectively."
\end{minipage}}
\end{equation}

Chebochko not only gave an
overview of the situation, she computed many new deforms:
\begin{Theorem}[Chebochko \cite{Ch1}]\label{NC2}
%\sssbegin[(Chebochko \cite{Ch1})]{Theorem}[\cite{Ch1}]
Let $L$ be a~Lie algebra over
a\emph{[n algebraically closed]} field of characteristic~$2$.
\begin{enumerate}\itemsep=0pt\setlength{\leftskip}{0.45cm}
\item[\emph{(Ch.~1)}] If $L$ is one of the types $A_l$ for $l>1$, $D_l$ for
$l\equiv
1\pmod 2$, $E_6$, $E_7$, $E_8$, $\overline{A_l}$ for $l\neq 3, 5$, or
$\overline{D_l}$ for $l\equiv 0\pmod 2$ and $l\neq 6$, then $H^2(L;
L)=0$.
\item[\emph{(Ch.~2)}] If $L$ is of type $\overline{A_l}$ for $l= 3, 5$, then $\dim
H^2(L; L)=20$.

\item[\emph{(Ch.~3)}] If $L$ is of type $D_4$, then $\dim H^2(L; L)=24$.

\item[\emph{(Ch.~4)}] If $L$ is of type $D_l$ for $l\equiv 0\pmod 2$ and $l>4$ or
$\overline{D_l}$ for $l\equiv 1\pmod 2$,
then $\dim H^2(L; L)=2l$.

\item[\emph{(Ch.~5)}] If $L$ is of type $\overline{D_6}$, then $\dim H^2(L;
L)=64$.

\item[\emph{(Ch.~6)}] If $L$ is of type $\overline{E_7}$, then $\dim
H^2(L; L)=56$.
\end{enumerate}
\end{Theorem}

This rather difficult to get result of Chebochko is expressed, as was customary at that time, in somewhat implicit form: above, only the dimension of the degree-2 cohomology is given whereas the weights of cocycle (also known to Chebochko) and their explicit form are sometimes useful.

\subsection[Difficulties computing (co)homology when $p>0$]{Difficulties computing (co)homology when $\boldsymbol{p>0}$}\label{ssDiffic} Let us
talk cohomology; for homology the situation is the same. Calculating
cohomology we somehow describe a~basis of the quotient module
$\operatorname{Ker}(d)/\operatorname{Im}(d)$, where $d$ is the codifferential.

The theorem, see \cite[pp.~28--29]{Fu},
\begin{equation*}%\label{ThH*}
\text{``\textit{The Lie $($super$)$algebra $\fg$ trivially acts on the spaces $H^i(\fg; M)$ and $H_i(\fg; M)$ for any~$i$}''}
\end{equation*}
holds for any characteristic~$p$. For $p>0$, however,
several obstacles arise:

1) If $p=0$, and $\fg$ is simple and finite-dimensional, then it
suffices to consider only cochains of weight~$(0, \dots, 0)$.

If $p>0$, to look for highest/lowest weight vectors
of weight whose coordinates are divisible by $p$ by using raising/lowering
operators does not suffice. The idea ``to seek the highest weight vectors in $\operatorname{Ker}(d)$ and
$\operatorname{Im}(d)$, separately, and compare them'' works only when $\operatorname{Ker}(d)$ is a~direct sum of submodules,
each having a~highest weight vector. This is not always so for
$p>0$. It can happen that $\operatorname{Ker}(d)$ and $\operatorname{Im}(d)$ have a~common
highest weight vector, but still are different modules.

2) With the help of \textsc{Mathematica}-based code \textsc{SuperLie}
we were able to \textit{directly}
compute $H^2(\fg;\fg)$ for $\fg$ of rank $\leq 4$ and not very big
dimension only. Kuznetsov and Chebochko were able to compute much more using the
Hochschild--Serre spectral sequence and induction on rank.

3) \textit{A~posteriori}, we see that the arguments given in \cite{Ch1, KuCh} proving absence of cocycles
of weight~0 are not applicable to Lie (super)algebras whose Cartan
matrix depends on a~parameter. \textit{We are unable, however, to indicate the requirement such examples fail to satisfy.}

4) The reduction of deformations in the KSh-method by passing from $H^2(\fg;\fg)$ to $H^2(\fg_-;\fg)$ is unjustified: e.g., it misses quantization and a~non-filtered deform found in~\cite{BGLLS2}, see Section~\ref{KacFil}.

\section{Introduction, continued: ``super'' occupies the scene}\label{SintroSuper}

\subsection{What ``Lie superalgebra'' is. The functorial definition. Questions 1)--3)}\label{ssDefFunct}

The functorial approach is needed, for example, to enable us to describe the actions of Lie (or algebraic) supergroups on Lie superalgebras or their duals. Speaking about Lie algebras and their representations in the simplest cases (finite-dimensional, over~$\Cee$) the naive definition suffices in many situations. There are, however, very important cases where the naive definition is insufficient to interpret the phenomena naturally arising in super setting even over $\Cee$, even for the most innocent-looking Lie superalgebras. For example:
\begin{equation}\label{3Q}
\begin{minipage}[l]{13cm}
\phantom{Xx}Q1) The Lie superalgebra $\mathfrak{vect}(0|n)$ of superderivations of the Grassmann algebra $\Kee[\xi]$ on $n$ generators $\xi=(\xi_1, \dots, \xi_n)$ should be the Lie superalgebra of the Lie supergroup of automorphisms of $\Kee[\xi]$. How to account for the odd parameters constituting $(\mathfrak{vect}(0|n))_\od$?\\
\phantom{Xx}Q2) The deformations of simple Lie superalgebras with odd parameter $\tau$, considered over the ground field, have an obvious ideal, although we'd like to consider these deformations as simple Lie superalgebras. \\
\phantom{Xx}How to interpret the deformations with odd parameter $\tau$ of irreducible modules over Lie superalgebras?\\
\phantom{Xx}Q3) The queertrace is odd, although it is natural to consider it as a~representation (same as any other linear functional on any Lie superalgebra $\fg$ that vanishes on the commutant $\fg^{(1)}$) whereas all representations are even by definition. How to resolve this contradiction?
\end{minipage}
\end{equation}

\subsubsection{Prerequisites for answers to the questions Q1)--Q3) in (\ref{3Q})}\label{ssPrereq}

{\bf Linear supervariety $\longleftrightarrow$ linear superspace.}
First, recall that there is a~one-to-one correspondence between every vector superspace $V$ (over the ground field $\Kee$) and the \textit{linear supermanifold} or \textit{linear supervariety}, i.e., a~ringed space ${\mathcal V}:=\big(V_\ev, {\mathcal L}\big(E^{\bcdot}\big(V_\od^*\big)\big)\big)$, where ${\mathcal L}(W)$ denotes the sheaf of sections of the vector bundle over $V_\ev$ with fiber $W$.

In various instances, e.g., dealing with actions of supergroups, it is more convenient to consider, instead of vector superspace $V$ and even the linear supermanifold ${\mathcal V}$, the functor ${\mathsf{ScommSalgs}}_\Kee\leadsto \Mod_C$ from the category of supercommutative $\Kee$-superalgebras $C$ to the category of $C$-modules represented by $V$ and ${\mathcal V}$:
\[
C\longmapsto {\mathcal V}(C)=V(C):=V\otimes C \quad \text{for any $C$},
\]
where ``any'' is understood inside a~suitable category (e.g., finitely generated over $\Kee$).

In these terms, a~\textit{Lie superalgebra} is a~vector superspace $\fg$, or a~linear supervariety (supermanifold) ${\mathcal G}$ corresponding to it, representing the functor from the category of supercommutative $\Kee$-superalgebras $C$ to the category of naively understood Lie superalgebras. To the Lie superalgebra homomorphisms (in particular, to representations) a~morphism of the respective functors should correspond.

Clearly, if $\fg$ is a~Lie superalgebra, then $\fg(C):=\fg\otimes C$ is also a~Lie superalgebra for any $C$ \textit{functorially in $C$}. The last three words mean that
\begin{equation}\label{funktor}
\begin{minipage}[l]{14cm}
for any morphism of supercommutative superalgebras $C\tto C'$, there exists a~morphism of Lie superalgebras $\fg\otimes C\tto \fg\otimes C'$ so that a~composition of morphisms of supercommutative superalgebras
\[
C\tto C'\tto C''
\]
goes into the composition of Lie superalgebra morphisms
\[
\fg(C) \tto \fg(C') \tto \fg(C'');
\]
the identity map goes into the identity map, etc.
\end{minipage}
\end{equation}

An \textit{ideal} $\mathfrak{h}\subset\fg$ represents the collection of ideals $\mathfrak{h}(C) \subset \fg(C)$ for every $C$.

\subsubsection[Affine superscheme (after \cite{L0} based on the Russian version of \cite{MaAG} pre\-print\-ed in 1968)]{Affine superscheme\\ (after \cite{L0} based on the Russian version of \cite{MaAG} pre\-print\-ed in 1968)}\label{Spec}

An affine superscheme $\operatorname{Spec} C$, where $C$ is a~supercommutative superalgebra or a~superring
is defined literally as the affine scheme: its points are prime ideals defined literally as in the commutative case, i.e., $\mathfrak{p}\subsetneq C$ is \textit{prime}\footnote{K.~Coulembier pointed out to us that the so far conventional definition in the non-commutative case is at variance with the common sense: at the moment, if \eqref{primeID} holds, $\mathfrak{p}$ is called (say, in Wikipedia) \textit{completely prime} while it would be natural to retain the term \textit{prime}, as is done in \cite{L0} and by P.~Deligne et al.\ in \cite{Del}, since the definition is the same despite the fact that supercommutative rings are not commutative, whereas the term \textit{prime} is (so far) applied to any ideal $P\subsetneq R$ in a~non-commutative ring $R$ which
for any two ideals $A$ and $B$ in $R$ satisfies the following version of \eqref{primeID}:
\[
\text{if $AB\subset P$, then either $A\subset P$, or $B\subset P$}.
\]}
\begin{equation}\label{primeID}
\text{if $a,b\in C$ and $ab\in \mathfrak{p}$, then either $a\in \mathfrak{p}$, or $b\in\mathfrak{p}$}.
\end{equation}
The affine scheme is endowed with Zariski topology and the structure sheaf, defined as in the commutative case, see~\cite{MaAG}.

Observe that in all works we saw published later than \cite{L0}, some published decades later, the only nontrivial point---how to define \textit{localization}?---was either ignored or answered wrongly. This question that caused a lot of trouble in the more non-commutative geometry than super geometry, cf.\ \cite{KoRo, Ro} and references therein. Should we consider one-sided denominators, say $d^{-1}n$ or $nd^{-1}$? In the supercommutative case, $d^{-1}n$ is naturally equivalent to $nd^{-1}$ for any numerator $n\in C$ and not only homogeneous with respect to parity denominator $d\in C\setminus \mathfrak{p}$ for any prime ideal $\mathfrak{p}$ in the supercommutative superalgebra (superring) $C$. This is easy to prove, but is worth stating explicitly.

{\bf Answers to questions (\ref{3Q}).} 1) Recall the definition of supermanifolds, see, e.g., \cite[Chapter~1]{Del}, \cite[Chapter~4, Section~1, no.~5]{MaG}. Same as manifolds are glued from coordinate patches locally diffeomorphic to a~ball in ${\mathbb R}^n$, supermanifolds are \textit{ringed spaces},\footnote{To understand the definition of ringed spaces, read~\cite{MaAG}---a shortest and most lucid source.} i.e., pairs ${\mathcal M}:=(M, \mathcal{O}_{\mathcal M})$, where $M$ is an $m$-dimensional manifold and $\mathcal{O}_{\mathcal M}$ is the structure sheaf of ${\mathcal M}$, locally isomorphic to ${C^\infty(U)\otimes \Lambda^{\bcdot} (n)}$, where $U$ is a~domain in $M$, and $\Lambda^{\bcdot} (n)$ is the Grassmann algebra with $n$ anticommuting generators. Pairs ${\mathcal U}:=(U, C^\infty(U)\otimes \Lambda^{\bcdot} (n))$ are called \textit{superdomains}. Morphisms of supermanifolds ${\mathcal M}\tto {\mathcal N}$ are pairs $(\varphi, \varphi^*)$, where $\varphi\colon M\tto N$ is a~diffeomorphism and $\varphi^*\colon \mathcal{O}_{\mathcal N}\tto \mathcal{O}_{\mathcal M}$ is a~preserving the natural parity of Grassmann-valued sheaves of functions morphism of sheaves (compatible with $\varphi$).

For a~superdomain ${\mathcal U}$ of superdimension $0|n$ over any ground field $\Kee$, the underlying
domain is a~single point. However, since $C^\infty ({\mathcal U}) =
\Lambda (n)$, the Grassmann algebra with $n$ anticommuting generators $\xi_i$, the superdomain ${\mathcal U}$ has a~lot of nontrivial (parity preserving)\footnote{Actually, there are lots of automorphisms of $\Lambda^{\bcdot} (n)$ that do not preserve parity, but the modern science does not know yet what to do with them, cf.\ \cite[Appendix~1, PP.~369--379]{Ls} and \cite{LSe}. D.L.\ conjectures, see~\cite{Lq}, that by analogy with supersymmetries wider than symmetries carried by groups or Lie algebras, these general automorphisms further widen supersymmetries. Recently, U.~Iyer proved that \textit{Volichenko algebras}, defined as the inhomogeneous subalgebras of Lie superalgebras, play the role of Lie algebras for the groups of such general automorphisms, see~\cite{I}.}
automorphisms, namely the group $\operatorname{Aut}^{\ev} \Lambda^{\bcdot} (n)$. All such automorphisms are of
the form
\begin{equation*}%\label{7eq8}
\xi_j
\mapsto \sum_{r} \varphi_j^{r}\xi_{r}+ \sum_{r>1}
\sum_{i_1<\dots <i_{2r+1}} \varphi_j^{i_1\dots
i_{2r+1}}\xi_{i_{1}}\dots\xi_{i_{2r+1}} \qquad \text{for all $j$},
\end{equation*}
where the matrix $(\varphi_j^{r})$ with elements in $\Kee$ is invertible. We see that such
automorphisms constitute the Lie (or at least algebraic) group whose Lie algebra consists of the even
elements of the Lie superalgebra $\mathfrak{vect}(0|n):={\mathfrak{der}}\, \Lambda^{\bcdot} (n)$. What corresponds to the odd
vector fields from ${\mathfrak{der}}\, \Lambda^{\bcdot} (n)$? We consider the answer in the following more general setting, but only for~$\Kee={\mathbb R}$ or~$\Cee$.

Let $E$ be a~trivial vector bundle over a~domain $U$ of dimension~$m$ with fiber
$V$ of dimension~$n$; let $\Lambda^{\bcdot}(E)$ be the exterior algebra of $E$. To the bundle
$\Lambda^{\bcdot}(E)$, we assign the superdomain ${\mathcal U}= (U, C^\infty({\mathcal U}))$, where
$C^\infty({\mathcal U})$ is the superalgebra of smooth sections of $\Lambda
^{\bcdot}(E)$. Clearly, each automorphism of the pair $(U,
\Lambda ^{\bcdot}(E))$, i.e., of the vector bundle $\Lambda
^{\bcdot}(E)$, induces an automorphism of the superdomain~${\mathcal U}$.

However, \textit{not all} automorphisms of the superdomain ${\mathcal U}$ are obtained in this
way. By definition, every~\textit{morphism of superdomains}
$\varphi\colon {\mathcal U}\tto {\mathcal V}$ is in one-to-one correspondence with
a~homomorphism of the superalgebras of functions
$\varphi^*\colon C^\infty({\mathcal V})\tto C^\infty({\mathcal U})$ sending 1 to 1.

Every homomorphism $\varphi^*$ is defined on the (topological\footnote{A \textit{topological algebra} $A$ over a~topological field
$\Kee$ is a~topological vector space together with a~bilinear multiplication
$A\times A\tto A$ such that $A$ is an algebra over~$\Kee$, continuous in a~certain sense. Usually the continuity of the multiplication means that the multiplication is continuous as a~map between topological spaces
${A\times A\tto A}$, a~set~$S$ is a~\textit{generating set of a~topological algebra}~$A$ if the smallest closed subalgebra of~$A$ containing~$S$ is~$A$.}) generators of
the superalgebra, in other words: coordinates. Consider the corresponding formulas
\begin{equation}
\label{7eq9}
 \begin{cases}
 \displaystyle \varphi^*(u_i) = \varphi^0_i(u) +
 \fbox{$\sum\limits_{r\geq 1}
 \sum\limits_{i_1<\dots<i_{2r}}\varphi_i^{i_1\dots
 i_{2r}}(u)\xi_{i_{1}}\cdots \xi_{i_{2r}}$} \quad \text{ for all $i$}, \vspace{1mm}\\
\varphi^*(\xi_j) = \sum_{r\geq 0} \sum_{i_1<\dots <i_{2r+1}}
 \psi_j^{i_1\dots
 i_{2r+1}}(u)\xi_{i_{1}}\cdots\xi_{i_{2r+1}} \quad \text{for all $j$}.
\end{cases}
\end{equation}
The terms $\varphi^*(u_i) = \varphi^0_i(u)$ determine an endomorphism
of the underlying domain $U$.

The linear terms $\varphi^*(\xi_j) = \sum_{i} \varphi_j^{i}(u)\xi_{i}$
determine endomorphisms of the fiber $V$ (its own fiber over each point, as the
dependence on $u$ shows).

The higher terms in $\xi$ from the right-hand side of expression of
$\varphi^*(\xi_j)$ determine an endomorphism of the
larger fiber---the Grassmann superalgebra of $E$.

All the above, except for the boxed terms in \eqref{7eq9}, existed in Differential Geometry, and no need to
introduce a~flashy prefix ``super'' was felt.

The difference between the vector bundle $\Lambda^{\bcdot}(E)$ and the
superdomain ${\mathcal U}= (U, C^\infty({\mathcal U}))$, where
$C^\infty({\mathcal U})=C^\infty(U)\otimes \Lambda
^{\bcdot}(V)$, is most easily understood when the reader looks at the boxed terms
in (\ref{7eq9}). These terms, meaningless in the conventional Differential
Geometry, make sense in the new paradigm:

{\it The category of superdomains has more morphisms than
the category of vector bundles: morphisms with non-vanishing boxed
terms in \eqref{7eq9} are exactly the additional ones}.

However, even the boxed terms in formula~\eqref{7eq9} is not all we get in the new
category: we still did not describe any of odd parameters of
endomorphisms. To account for the odd parameters, we have to consider the functor
$C\longmapsto \operatorname{Aut}_C^{\ev}(C^\infty({\mathcal U})\otimes C)$, i.e., the parity preserving $C$-linear automorphisms of the form
\begin{gather*}%\label{7eq99}
 \varphi^*(u_i) = \varphi^0_i(u) +
 \fbox{$\sum\limits_{r\geq 1}
 \sum\limits_{i_1<\dots<i_{r}}\varphi_i^{i_1\dots
 i_{r}}(u)\xi_{i_{1}}\dots \xi_{i_{r}}$}\quad \text{for all $i$}, \\
 \varphi^*(\xi_j) = \sum\limits_{r\geq 0} \sum\limits_{i_1<\dots <i_{2r+1}}
 \psi_j^{i_1\dots
 i_{2r+1}}(u)\xi_{i_{1}}\dots\xi_{i_{2r+1}}\\
 \hphantom{\varphi^*(\xi_j) =}{}
+\fbox{$\psi^0_j(u)+
\sum\limits_{r\geq 1} \sum\limits_{i_1<\dots <i_{2r}}
 \psi_j^{i_1\dots
 i_{2r}}(u)\xi_{i_{1}}\dots\xi_{i_{2r}}$}\quad \text{for all $j$},
\end{gather*}
where
\begin{gather*}
\begin{split}
& \varphi^0_i(u), \varphi_i^{i_1\dots
 i_{2r}}(u), \psi_j^{i_1\dots
 i_{2r+1}}(u)\in C_\ev \quad \text{for all $r$ are even parameters,}\\
& \psi^0_j(u), \varphi_i^{i_1\dots
 i_{2r+1}}(u), \psi_j^{i_1\dots
 i_{2r}}(u)\in C_\od \quad \text{for all $r$ are odd parameters}
 \end{split}
 \end{gather*}
 of the supergroup of automorphisms of $C^\infty({\mathcal U})$, i.e., of the Lie superalgebra $\mathfrak{vect}(m|n)$, if considered infinitesimally.

The supervarieties isomorphic to the ringed spaces locally described by the section of the exterior algebra of a~vector bundle are called \textit{split}. All smooth supermanifolds are split (see \cite[Section~4.1.3]{MaG} more instructive than the first publications \cite{Bat, Ga}), whereas \textit{in the categories of superschemes, algebraic, and analytic supervarieties, there are more objects than in the category of vector bundles}: there are non-split supervarieties, see the papers \cite{Gre, Pa} by Green and Palamodov, who were the first to understand this and gave examples; for further examples, see \cite{MaG, Va} and papers by Onishchik and his students, see, e.g., \cite{OVs} and references therein with amendments in \cite{Lob, L3}.

2) \textit{Deformations and deforms. Odd parameters.}
Let ${\mathcal C}$ be a~finitely
generated supercommutative superalgebra or superring, let $\operatorname{Spec} {\mathcal C}$ be the affine super scheme, see Section~\ref{Spec}.

A~\textit{deformation} of
a Lie superalgebra $\fg$ over $\operatorname{Spec} {\mathcal C}$
is a~Lie superalgebra ${\mathfrak G}$ over ${\mathcal C}$ such that for some closed point $\mathfrak{p}_I\in \operatorname{Spec} {\mathcal C}$ corresponding to a maximal ideal in $I\subset {\mathcal C}$, we have ${\mathfrak G}\otimes_I\Kee\simeq\fg$ as Lie superalgebras over $\Kee$. Note that since ${\mathcal C}/I\simeq\Kee$ as $\Kee$-algebras,
this definition implies that ${\mathfrak G}\simeq\fg\otimes {\mathcal C}$ as ${\mathcal C}$-modules. The deformation is \textit{trivial} if ${\mathfrak G}\simeq\fg\otimes {\mathcal C}$ as Lie superalgebras over ${\mathcal C}$, not just as modules, and \textit{non-trivial} otherwise. (Here, we superized \cite{Ru}, where the non-super case was considered.)

Generally, the \textit{deforms}---the results of a~deformation---of
a~Lie superalgebra $\fg$ over $\Kee$ are Lie superalgebras ${\mathfrak G}\otimes_{I'}\Kee$, where $\mathfrak{p}_{I'}$ is any closed point in $\operatorname{Spec} {\mathcal C}$. People are usually interested in deforms given by
$I'\neq I$, since the deform given by $I'=I$ is just $\fg$.

In particular, consider a~deformation with an odd parameter $\tau$ of a~Lie superalgebra $\fg$ over field $\Kee$. This is a~Lie superalgebra ${\mathfrak G}$ over $\Kee[\tau]$ such that ${\mathfrak G}\otimes_I\Kee\simeq\fg$, where $I=
\langle \tau\rangle$ is the only maximal ideal of $\Kee[\tau]$. This implies that ${\mathfrak G}$ isomorphic to $\fg\otimes\Kee[\tau]$ as a~\textit{module over $\Kee[\tau]$}; if, moreover, ${\mathfrak G}=\fg\otimes\Kee[\tau]$ as a~\textit{Lie superalgebra over $\Kee[\tau]$}, i.e.,
\[
[a\otimes f, b\otimes g]=(-1)^{p(f)p(b)}[a,b]\otimes fg \qquad \text{for all $a,b\in \fg$ and $f,g\in\Kee[\tau]$},
\]
then the deformation is considered \textit{trivial} (and \textit{non-trivial} otherwise). Since $\langle \tau\rangle$
is the only maximal ideal of $\Kee[\tau]$, the only deform of $\fg$ produced by a deformation with an odd parameter is $\fg$ itself. Observe that $\fg\otimes \tau$ is not an ideal of ${\mathfrak G}$: any ideal should be a~free $\Kee[\tau]$-module.

\begin{Comment} In a~sense, the people who ignore odd parameters of deformations have a~point: they consider classification of simple Lie superalgebras over the ground field $\Kee$, right? On the one hand, this is correct, but it is not right: the deforms parametrized by $\big(H^2(\fg;\fg)\big)_\od$ are no less natural than the odd part of the deformed Lie superalgebra itself. To take these parameters into account, we have to consider everything not over $\Kee$, but over $\Kee[\tau]$. We do the same, actually, when $\tau$ is even and we consider formal deformations over $\Kee[[\tau]]$. If $\tau$ is even, and the formal series converges in a~domain $D$, we can evaluate the series and consider copies $\fg_\tau$ of $\fg$ for every~$\tau\in D$.
\end{Comment}

If the parameter is formal (the series diverges) or odd, this interpretation of deforms is impossible.

3) \textit{Two examples of Lie superalgebras}. The naive Lie superalgebra $\mathfrak{gl}(a|b)$ represents the functor
\[
C\longmapsto \mathfrak{gl}_C(a|b)=\mathfrak{gl}(a|b)(C):=\mathfrak{gl}(a|b)\otimes_{\Kee}C.
\]

The naive Lie superalgebra $\mathfrak{gl}(n)$ represents the functor
\[
C\longmapsto  \mathfrak{q}_C(n)=\mathfrak{q}(n)(C):=\mathfrak{q}(n)\otimes_{\mathbb{K}}C.
\]

4) \textit{Representations of Lie superalgebras}.
By the naive definition, a~\textit{representation} of a~given Lie superalgebra $\fg$ in a~superspace $V$ is a~\textit{morphism} (i.e., an \textit{even} homomorphism) $\fg\tto \mathfrak{gl}(V)$.

In terms of the functor of points, a \textit{representation} of a~given Lie superalgebra $\fg$ in a~superspace~$V$ is a~collection of \textit{$C$-linear morphisms} of Lie superalgebras $\fg(C)\tto \mathfrak{gl}_C(V)$ which is functorial in $C$.

Observe that any \textit{trace} (for emphasis called \textit{super} trace in the super setting) defines a~representation of the Lie (super)algebra on which it is defined $\operatorname{tr}\colon\fg\tto\mathfrak{gl}(1)$.

But a~super trace can be \textit{odd}. The collection of representations is not a~vector space but a~supervariety that may have odd parameters. To account for them, introduce an odd parameter, call it $\tau$.
The odd queertrace sends the bracket to the bracket (thanks to the fact that the bracket in the target Lie algebra is identically zero) \[
\operatorname{qtr}\colon \ \mathfrak{q}_C(n)\tto \mathfrak{gl}_{C[\tau]}(1)=C_\ev\oplus \tau C_\od, \qquad (A, B) \mapsto \tau \operatorname{tr} B.
\]

\subsection[The functorial definition for $p=3$ and 2]{The functorial definition for $\boldsymbol{p=3}$ and 2}\label{ssDefFunct3,2}

 Instead of the above definition~\eqref{funktor}, one mostly encounters another definition which reduces everything to Lie (or algebraic) groups and Lie algebras: namely the functor from the category of supercommutative and associative algebras with 1 to the category of Lie algebras given by
\begin{equation}\label{fun}
C\longmapsto \fg(C):=(\fg\otimes C)_\ev,
\end{equation}
see, e.g., J.~Bernstein's lectures in \cite[Sections~1.2, 2.8, 2.9]{Del}, or
``Grassmann envelopes''\footnote{\textit{Sometimes} Berezin used in~\cite{B} the term ``Grassmann envelope'' speaking about the whole $\fg\otimes C$, not its even part.} in~\cite{B}.

The definition \eqref{funktor} did not come to mind earlier for psychological reasons only, it seems: to reduce one Lie superalgebra to the collection of Lie superalgebras considered over ALL supercommutative superalgebras $C$ did not look as a~most reasonable idea. Besides, over~$\Cee$ or~${\mathbb R}$---the ground fields for most of the physical papers---to consider the functor ${C\longmapsto \fg(C):=(\fg\otimes C)_\ev}$ was both sufficient and psychologically comfortable: we arrived to the realm without any ``super''.

For $p=3$ and $2$, the functor \eqref{fun} does not, however, catch the cubic addition \eqref{6} to the Jacobi identity, and squaring, respectively. For details, see~\cite{KLLS}.

\section[Introduction, continued: super versions of the KSh-conjecture]{Introduction, continued: super versions\\ of the KSh-conjecture}\label{SupKSh}

\subsection[Lie superalgebras over $\Cee$: simple finite-dimensional and vectorial with polynomial\\ coefficient]{Lie superalgebras over $\boldsymbol{\Cee}$: simple finite-dimensional\\ and vectorial with polynomial coefficients}\label{ssExOverC}

The KSh-method and conjecture are based on the two main ingredients: simple finite-di\-men\-sion\-al and simple vectorial with polynomial coefficients Lie algebras over $\Cee$, see Conjecture~\ref{DKSh}. All these ingredients, except for the Lie algebras of Hamiltonian vector fields that can be quantized, are rigid over $\Cee$.

Let us briefly recall their super counterparts; some of them are not rigid and the supervariety of parameters of deformation is rather involved.
\begin{itemize}\itemsep=0pt
\item The Lie superalgebras of the form $\fg(A)$ with indecomposable Cartan matrix $A$ and their simple relatives. These Lie superalgebras (or their simple relatives for $A$ non-invertible) were independently discovered by several (groups of) reserchers; V.~Kac was the first to observe that each of these Lie superalgebras (or, for $A$ non-invertible, its \textit{double extension}, for their succinct definition of this important notion, see \cite{BLS}) has a~Cartan matrix. Moreover, one Lie superalgebra can have several Cartan matrices. V.~Serganova described a~method enabling one to list all inequivalent Cartan matrices of the given Lie superalgebra. (For details and history, see \cite{CCLL}.)
\end{itemize}

The list of simple finite-dimensional Lie superalgebras of the form $\fg(A)$ and their simple relatives is as follows.
\begin{enumerate}\itemsep=0pt
\item[--] The Lie superalgebras of the form $\fg(A)$ are of series $\mathfrak{gl}$ if $A$ is not invertible, or $\mathfrak{sl}$ if $A$ is invertible, or the ortho-symplectic Lie superalgebras $\mathfrak{osp}(m|2n)$ that preserves the \textit{even} non-degenerate symmetric bilinear form on the $m|2n$-dimensional superspace.
\item[--]Two exceptional simple Lie superalgebras with Cartan matrix discovered by Kaplansky (we denote them $\mathfrak{ag}(2)$ and $\mathfrak{ab}(3)$, as suggested by their root systems, see~\cite{Se4}). The meaning (interpretation: what do they preserve?) of these Lie superalgebras was unclear until Elduque found their analogs in characteristics $p=3$ and~5 using pairs of composition algebras, and organized them, and several new simple exceptional Lie superalgebras, in super magic super rectangle (superization of Freudenthal's magic square), see~\cite{BeEl, El}.
\end{enumerate}

The serial algebras of the form $\fg(A)$ have ``purely super'' or, better say, ``odd'' analogs:
\begin{itemize}\itemsep=0pt
\item The queer Lie superalgebra $\mathfrak{q}(n)$ is an analog of $\mathfrak{gl}(a|b)$; it has its own (queer) trace, which is odd in both senses, and subquotients $\mathfrak{pq}(n)$, $\mathfrak{sq}(n)$, and $\mathfrak{psq}(n)$.

\item The periplectic Lie superalgebra $\mathfrak{pe}(n)$ is an analog of $\mathfrak{osp}$ preserving the \textit{odd} non-de\-gen\-er\-ate symmetric bilinear form on the $n|n$-dimensional superspace. Its supertraceless subalgebra $\mathfrak{spe}(n)$ is simple for $n>2$.

\item \textit{Deformations}. The simple Lie superalgebras of the above-listed types are rigid, except for $\mathfrak{osp}(4|2)$ and $\mathfrak{spe}(3)\simeq\mathfrak{svect}(0|3)$, where $\mathfrak{svect}(a|b)$ is the Lie superalgebra of divergence-free vector fields on the $a|b$-dimensional supermanifold, depending on an even and an odd parameter, respectively.

\item Vectorial simple Lie superalgebras. It is natural to consider these Lie superalgebras with an invariant (\textit{Weisfeiler}) filtration; as such they constitute 34 series and 15 exceptions. Being interested in their finite-dimensional versions in characteristics $p>0$, we consider them as abstract, their W-filtration (short for Weisfeiler filtration) forgotten:
\begin{equation}\label{summaryAbs}
\begin{minipage}[l]{14cm}
\textit{As abstract algebras, $W$-gradings associated with the $W$-filtrations disregarded, simple
$W$-graded vectorial Lie superalgebras can be collected into $14$
series $($one containing an exceptional subseries$)$ and $5$ exceptional
Lie superalgebras. Deforms of some of them are parameterized by a
supervariety}. For details, see \cite{LSh1+, LSh1, LSh3,LSh5, Sh14}, where all examples were found, and where the classification of simple $\Zee$-graded vectorial Lie algebras was announced and proved, except for the case of $\Zee$-gradings compatible with parity. This gap in the proof was filled in \cite{CCK, CaKa2, CaKa4,CK1, CK2, CK1a,K7, K10}. For a~succinct review, see \cite{LSh5} to be expanded and published shortly.
\end{minipage}
\end{equation}
Over $\Cee$, the filtered deformations of the corresponding to the 15 W-filtrations of the 5 exceptional simple vectorial Lie superalgebras (all of which are discovered by Irina Shchepochkina) are computed in \cite{GLS}.
\end{itemize}

None of the papers listed in \eqref{summaryAbs} considered complete classification of deformations, except for \cite{LSh3} in certain cases; moreover, deformations with odd parameter (known to experts for years) and rediscovered in~\cite{CK1} were thrown away by the authors of \cite{CK1}.

Although this is off-topic, we'd like to point at a~very interesting paper~\cite{CaKa1} classifying vectorial Lie superalgebras over ${\mathbb R}$; certain real forms were new and unexpected.

\subsection[Finite-dimensional modular Lie (super)algebras of the form $\fg(A)$\\ with indecomposable Cartan matrix $A$ are classified for any $p$]{Finite-dimensional modular Lie (super)algebras of the form $\boldsymbol{\fg(A)}$\\ with indecomposable Cartan matrix $\boldsymbol{A}$ are classified for any $\boldsymbol{p}$}\label{ssDefModular}

Weisfeiler and Kac~\cite{WK} gave a~classification of finite-dimensional Lie algebras
$\fg(A)$ with indecomposable Cartan matrix $A$ for any $p>0$, but
although the idea of their proof is OK, the paper has a
gap\footnote{Brown discovered the simple Lie algebra $\mathfrak{br}(3)$, see \cite{Br3}, but
failed to observe that it has Cartan matrix.
Weisfeiler and Kac missed $\mathfrak{br}(3)$ in their classification of
simple finite-dimensional Lie algebras over $\Kee$ with
indecomposable Cartan matrix \cite{WK}. It was Skryabin \cite{Sk1}
who found this gap in \cite{WK} and all (two) Cartan matrices of $\mathfrak{br}(3)$. Neither
\cite{Sk1}, nor \cite{KWK} claimed this was the only gap in \cite{WK}; this was
stated in \cite{BGL2}.} and several related notions, such as the
definition of the Lie algebra with Cartan matrix, Dynkin diagram,
and roots were vague or absent. Although the definition of Cartan matrix
given in \cite{K} is also applicable to Lie \textit{super}algebras
and modular Lie (super)algebras, it was not properly developed
and formulated at the time \cite{WK} was written; the algebras
$\fg(A)^{(i)}/\fc$, that have no Cartan matrix, are sometimes
referred to as having one, even nowadays. All these notions, and
several more, are clarified in \cite{BGLL, BGL2}, see also
\cite{CCLL}, where, for Lie superalgebras, there are given arguments in favor
of descriptive notation and against ``notation \`a la Cartan''
reasonable only for root systems of simple Lie algebras over~$\Cee$; these arguments, given
in \cite{CCLL} over~$\Cee$, work over~$\Kee$ as well:
\begin{enumerate}\itemsep=0pt
\item[(1)] $\mathfrak{gl}(n)$, $\mathfrak{sl}(n)$,
$\mathfrak{pgl}(n)$ and $\mathfrak{psl}(n)$ are of the same $A_{n}$ root type;
\item[(2)] one
Lie (super)algebra (root system) may have several inequivalent
Cartan matrices.
\end{enumerate}

All finite-dimensional Lie superalgebras $\fg(A)$ with indecomposable Cartan matrix $A$ for any $p>0$ are classified in \cite{BGL2}.

\subsection{General remarks}\label{ssGenRem}

1) The simple Lie algebras constitute
a natural first object to study, but some
of their relatives are even ``better'' having some nice properties, e.g., are restricted, have Cartan
matrix, and so on, whereas their simple ``cores'' may lack one or
all of these features, see \cite{BLLS1}.

2) Several ``mysterious'' examples of simple Lie algebras over fields
$\Kee$ for $p=3$, mainly due to Skryabin \cite{Sk}, were partly
described in \cite{S,S-II,S-III}; in \cite{GL3}, these examples are demystified
as the CTS prolongs (or deformations thereof); hidden parameters of the
shearing vectors were found in the process.

It was clear since long ago that the smaller characteristic, the
less rigid the simple Lie algebras are; Rudakov gave an example of a
3-parameter family of simple Lie algebras\footnote{The discovery of this family was
sometimes misattributed (even---for some reasons---by Rudakov himself in~\cite{Ru})
to Kostrikin, who first published one of the family's
descriptions, or to Dzhumadildaev, who studied this family.}, see~\cite{GL3, Kos}. After Rudakov's example became known, Kostrikin and
Dzhumadildaev \cite{Dz3,Dz2,Dz1, DK} studied various (e.g.,
filtered and infinitesimal) deformations of simple vectorial Lie
algebras (in~\cite{KSh}, these were dubbed ``algebras of Cartan type'' because no other vectorial Lie algebras except for the 4 series discovered by Cartan were known at that time); for a
detailed summary of the part of their results with understandable detailed
proofs, and some new results (all pertaining to the infinitesimal
deformations), see \cite{Vi1,Vi2,Vi1a, Vi4,Vi3}.

Rudakov's paper \cite{Ru} clearly showed that speaking about
deformations it is unnatural to consider the modular Lie algebras
naively, as vector spaces: Lie algebras should be viewed as algebras
in the category of varieties. This approach should, actually, be
applied even over fields of characteristic~0, but the simplicity of
the situation with simple finite-dimensional Lie algebras obscured this
(a posteriori obvious) fact. Deforms with odd parameters of
Lie superalgebras can only be viewed
as Lie algebras in the category of supervarieties; for examples over~$\Cee$, see~\cite{Lq,LSh3}.

3) Dzhumadildaev and
Kostrikin \cite{DK} claimed that every infinitesimal deformation
of the Lie algebra $\mathfrak{vect}(1; \underline{m})$, a.k.a.\ $W_1(m)$, is
integrable. Most of the infinitesimal deformations of the Lie
(super)algebras considered here are integrable; moreover, the global
deformation corresponding to a~given (homogeneous with respect to
weight) cocycle~$c$ of $\fg$ is often \textit{linear} in parameter, i.e., the deformed bracket
is of the form
\begin{equation}\label{defglob}
[x,y]_{{\rm new}}=[x,y]+t c(x, y), \text{~~where $t$ is a
parameter, for any $x,y\in\fg$}.
\end{equation}
The cocycles representing the classes of $H^2(\fg; \fg)$ correspond
to infinitesimal deformations, but if the cocycle is odd, it
certainly determines a~global deformation of the bracket \eqref{defglob} with an odd
parameter.

4) Trying to append the results of Chebochko \cite{Ch1} with deformations
of the Lie algebras she did not consider but which we consider no
less ``classical'' than the ones she did consider, we obtained, as a
byproduct, an elucidation of Shen's ``variations'', see~\cite{She1}.

Shen described seven ``variations'' $V_iG(2)$ of $\fg(2)$ and three
more ``variations'' of $\mathfrak{sl}(3)$; Shen claimed that all his examples
are simple and all but two ($V_1G(2)$ which is $\mathfrak{psl}(4)$ and
$V_7G(2)$ which is $\fwk(3; a)^{(1)}/\fc$, see~\cite{BGL2}) are new.
It was later found, see~\cite{LLg}, that the ``variations'' of
$\mathfrak{sl}(3)$ are isomorphic to $\mathfrak{sl}(3)$, the variations of dimension~15 are not simple.\footnote{Several descriptions of Shen's
``variations'' have typos: as written, these algebras do not satisfy
Jacobi identity or have ideals. Unfortunately, we were unable to
guess how to amend the multiplication tables whereas our letters to
Shen and his students remain unanswered.}

One of the algebras Shen discovered, $\fg\mathfrak{s}$, and its CTS prolong, partly
investigated by Brown (see~\cite{Bro}), are remarkable, see~\cite{BGLLS}. They qualify as ``standard''
Lie algebras in the sense of Conjecture~\ref{DKSh}.

\subsubsection{Main trouble: how to interpret numerous deformations
found?}\label{ssTooMany} In \cite{Ch1, KKCh, KK, KuCh}, a~complete description of
deformations of the ``classical'' simple Lie algebras with Cartan
matrix (and several of their non-simple, but important relatives, see
open problem~2) in Section~\ref{OP}) is performed over fields $\Kee$
for $p=3$ and 2.

Absolutely correct---in terms of the conventional definition of
Lie algebra (co)homology---computations of Dzhumadildaev, elucidated in
\cite{Vi1, Vi1a, Vi3}, imply that the vectorial Lie algebra $\mathfrak{vect}(n;
\underline{N}):={\mathfrak{der}}(\mathcal{O}(n;\underline{N}))$ has lots of
\textit{infinitesimal} deformations. The seemingly overwhelming
abundance of cocycles found in \cite{Ch1,DK, Vi1, Vi2, Vi1a, Vi3}, is less big
problem as we thought before we have read \cite{Ch1}: although there
are many nontrivial cocycles, nonisomorphic deforms of $\fg$
correspond only to the distinct orbits in the space of cocycles
under the action of the Chevalley group $\operatorname{Aut}(\fg)$ and there are
very few such orbits.

Another reason for abundance of cocycles is
appearance of numerous ``semi-trivial'' cocycles, see
\cite{BLLS, BLW}.\footnote{The situation with these too numerous cocycles
is opposite, in a~sense, to that with the Kac--Moody groups (which
``did not exist'' until a~correct definition of cohomology was used,
see~\cite{Rey}, and appendices in~\cite{FuR} not translated in~\cite{Fu}), and
Krichever--Novikov's algebras that exist despite the (correct in accepted definitions)
nonexistence theorems in~\cite{LR}. (Similar lack of understanding
was at first with Dirac's $\delta$-function which is not a~function in the
used-to-be-conventional sense.)} See also~\cite{FF}.

\subsection[Towards classification of simple finite-dimensional modular Lie superalgebras\\ over algebraically closed fields of characteristic $p>0$]{Towards classification of simple finite-dimensional modular Lie\\ superalgebras over algebraically closed fields of characteristic $\boldsymbol{p>0}$}\label{ssTowards}

\subsubsection{Examples of simple ``non-symmetric'' Lie superalgebras}

For $p>5$, the list of simple finite-dimensional Lie
\textit{super}algebras over $\Kee$ is, \textit{conjecturally}, the list of
straightforward modular analogs of the list of finite-dimensional and $\Zee$-graded
vectorial Lie superalgebras over $\Cee$, see \cite{ BGLLS2,LSh1+, LSh1} for details and \eqref{summaryAbs} for a~summary, \textit{together with filtered deformations} in the following cases
\begin{itemize}\itemsep=0pt
\item the divergence-free case,
\item the two types of super analogs of the Hamiltonian algebras,
\item the two types of super analogs of contact algebras,
\item the filtered deformations of the 15 Shchepochkina's exceptions (for these deformations over~$\Cee$, see~\cite{GLS}).
\end{itemize}

For $p=5$ and $3$, there are some indigenous examples of ``non-symmetric'' algebras~\cite{BL1, BGL1}; more might appear. We do not dare to conjecture.

For $p=2$, manifestly the most difficult case, a~\textit{miracle happens}: every finite-dimensional simple Lie
superalgebra is obtained from a~simple Lie algebra by one of the two
methods described in~\cite{BLLS1}. So,~\cite{BLLS1} contains a
classification, although modulo the classification of simple Lie
algebras. (The miracle is marred by a~threat hinted at by the recent result~\cite{KrLe}: the ``simpler'' task of classifying simple Lie algebras might be wild.)

For new simple CTS-prolongs, see \cite{BGLLS} and other examples, e.g., \cite{Ei, GZ, SkT1}.

\subsubsection{Examples of simple ``symmetric" Lie superalgebras}
\begin{itemize}\itemsep=0pt
\item The Lie superalgebras with Cartan matrix. The Lie superalgebras $\fg(A)$ with indecomposable Cartan matrix $A$ are classified, together with their simple relatives, \cite{BGL2} for any $p>0$.
\item The queer Lie
superalgebras $\mathfrak{psq}(n)$ for $p\neq 2$,
\item If $p=2$, there are lots of simple ``\textit{symmetric}'' algebras. In addition to the above, there are
\begin{enumerate}\itemsep=0pt
\item[(1)] various ortho-orthogonal and periplectic algebras;
\item[(2)] if $\fg$ is a~symmetric simple Lie algebra, then its
(\textit{partial} if $\fg$ is not restricted) \textit{queerification} $\widetilde{\mathfrak{q}}(\fg)$,
described in~\cite{BLLS1}, is a~symmetric simple Lie superalgebra.
\end{enumerate}
\end{itemize}

Various simple ``relatives'' of the Lie (super)algebras with Cartan
matrix do not have any Cartan matrix, see \cite[Warning~4.1]{BGL2}.
In \cite{ChG2, Ch1}, Chebochko described infinitesimal deforms
of the Lie algebras with indecomposable Cartan matrix (or their simple relatives), except for
those listed in~\eqref{propusk}.

\subsection{Our results}\label{ssResults}
\begin{enumerate}\itemsep=0pt
\item[1)] We append the earlier results by other researchers and consider
the cases omitted not only among the simple Lie algebras but also
among their closest relatives, which have Cartan matrix, like
$\mathfrak{gl}(pn)$, whereas neither $\mathfrak{sl}(pn)$, nor the simple Lie algebra
$\mathfrak{psl}(pn)$ have Cartan matrix. The deforms of $\mathfrak{br}(3)$ we describe here
is a~new result. (Deforms of~$\mathfrak{br}(2)$ and~$\mathfrak{o}(5)$
were known, for a~precise and correct description of isomorphism classes of these deforms, see \cite{BLW}.)

\item[2)] For $p=2$, we describe deformations
of the following ``symmetric'' Lie (super)algebras:
all known\footnote{Thanks to results of Skryabin and Eick (see~\cite{Ei, SkT1})
further explained and expounded in \cite{GGRZ, GZ}, see also \cite[Comment~2.1.1] {BGLL1}.}
``symmetric'' Lie algebras and all (known) Lie superalgebras of rank
$\leq 4$, except queerifications, see~\cite{BLLS1}, of which
we considered only the simplest ones: for illustration. In particular, we describe all
deformations of the Lie algebras named in list~\eqref{propusk}, except for the superization of the simple Lie algebras of the $ADE$ root system.
This is an \textit{open problem}.

For $p=2$, computer experiments show a~one-to-one correspondence between deformations of the restricted simple Lie algebra $\fg$ and its superizations by means of ``method~2'' for whose definition, see~\cite{BLLS1}. \textit{Open problem}: explain this correspondence.

\item[3)] In Section~\ref{App}, we give some results concerning $H^{\bcdot}(\fg)$. We stopped investigating the multiplicative structure in the space
$H^{\bcdot}(\fg)$ because it is too cumbrous and the result does not seem
to be worth the effort, as we illustrate in one case.
\end{enumerate}

In this paper we did not need, actually,
that much information about cohomology with trivial coefficients, only $H^{i}(\fg)$ for $i<3$; we
computed the whole space $H^{\bcdot}(\fg)$ out of curiosity and in the (vain, as it turned out)
hope to get a~nice-looking result.

\subsection{Open problems}\label{OP} %\textbf{Open problem}.
\index{Problem, open}

The most profound ones are listed in the next section. On a smaller scale, we consider:
\begin{enumerate}\itemsep=0pt
\item[1)] For $p=2$, the simple Lie
superalgebras which are not queerifications, the result of
\cite{BLLS1} gives us a~choice: either describe deformations of these Lie
superalgebras or describe $\Zee/2$-grading of the simple Lie
algebras they superize. Which method is easier and more appropriate
and in which cases? (Both methods lead to the same result which
helps to verify both, and (in some cases) solve the next open
problem, cf.\ Lemma~\ref{L5.2} and~\cite{KrLe}.)

\item[2)] Which of the deformations found (e.g., in this paper), and to be found, are true ones?

\item[3)] For the Lie algebras $\fg=\fg(A)$ with non-invertible Cartan
matrix $A$, Chebochko only considered $\fg^{(1)}$ and
$\fg^{(1)}/\fc$. She did not consider, e.g., the algebras we denote~$\fe(7)$ and~$\fe(7)/\fc$; her $E_7$ is our $\fe(7)^{(1)}$ and her
$\overline{E_7}$ is our $\fe(7)^{(1)}/\fc$; similar are
correspondences between Chebochko's $A_{2n-1}$ and $\overline{A_{2n-1}}$ which are,
speaking prose, $\mathfrak{sl}(2n)$, and $\mathfrak{psl}(2n)$. In the cases listed in
\eqref{NC2}, \textit{and all other cases with non-invertible Cartan
matrix}; it remains to consider the cases of $\fg$ (and, perhaps, of $\fg/\fc$ such as $\mathfrak{gl}(pn)$ and $\mathfrak{pgl}(pn)$).

The same concerns relatives of simple Lie algebras without Cartan matrix,
such as~$\mathfrak{o}_{\rm gen}(n)$ (see \cite{BGLL1}), $\mathfrak{o}_I(2n)$, and (for any $p$) Cartan prolongs---more natural
objects than their simple relatives.

\item[4)] Which of the multi-parameter deformations are integrable? Cf.\ Lemma~\ref{integra} and~\cite{BLW}.

\item[5)] Investigate classes of isomorphic deforms \`a la Kuznetsov and Chebochko \cite{KuCh}.
(Breaking news: In \cite{KuChA5} Kuznetsov and Chebochko classified non-isomorphic deforms of $\mathfrak{psl}(6)$ for $p=2$, see also some partial results in \cite{KuKCh2} which is an English version of~\cite{KuKCh1}.)

\item[6)] Describe deformations of numerous queerifications not considered in this paper for $p=2$.

\item[7)] Describe deformations of superization of the simple Lie algebras of the $ADE$ root system for $p=2$, currently performed, partly, only for the $A$-type systems.
\end{enumerate}

\section{(Co)homology of Lie superalgebras}\label{Scoh}

First, recall the classical definitions of homology and cohomology for Lie algebras:
\begin{itemize}\itemsep=0pt
\item {\it Cohomology}: Scientific definition, due to Chevalley and Eilenberg, says: given a~\textit{left} module $M$ over the Lie algebra~$\fg$, both being modules over a~commutative ring $R$ (in particular, a~ground field $\Fee$) the \textit{$n$th cohomology of $\fg$ with coefficients in $M$} is
\begin{equation}\label{coho}
H^n(\fg; M):=\operatorname{Ext}^n_{U(\fg)}(R, M),
\end{equation}
where $R$ is considered as the trivial $\fg$-module.
\end{itemize}

This definitions boils down to the following simple-minded one coming from the left-invariant de~Rham cohomology on the Lie group with Lie algebra $\fg$; this definition can be explained to a~computer: on the space $E^{\bcdot}(\fg^*)\otimes_R M$ of \textit{cochains}, the codifferential, whose square vanishes, is defined and the corresponding cohomology are called \textit{cohomology of $\fg$ with coefficients in $M$}.

\begin{itemize}\itemsep=0pt
\item {\it Homology}: Analogously, let $M$ be a~\textit{right} module over the Lie algebra $\fg$; set
\[
H_n(\fg; M):=\operatorname{Tor}_n^{U(\fg)}(R, M),
\]
which boils down to the following: on the space $M\otimes E^{\bcdot}(\fg)$ of \textit{chains}, the differential, whose square vanishes, is defined and the corresponding homology are called \textit{homology of~$\fg$ with coefficients in~$M$}.
\end{itemize}

In this paper and in \cite{BGLL1} we are dealing only with cohomology; for applications of homology, see \cite{Fu}.

\begin{OpenProblem}
\index{Problem, open} Do definitions of cochains in what follows correspond to the definition~\eqref{coho}?
\end{OpenProblem}

\begin{itemize}\itemsep=0pt
\item {\it Cohomology of Lie \textit{super}algebras}: Consider this case in detail; changes needed to treat homology are usually considered obvious, but to be on the safe side we will not go into this now and skip discussing homology. %\AL{}
According to the above, we can consider two types of the spaces of cochains of a~given Lie superalgebra $\fg$ with coefficients in the \textit{left} $\fg$-module~$M$:
\begin{equation}\label{EandLambda}
E^{\bcdot}(\fg^*)\otimes M \qquad \text{and} \qquad \Lambda^{\bcdot}(\fg^*)\otimes M.
\end{equation}
Observe that, for some reasons now unknown forever, Grozman made \textsc{SuperLie} to express the elements of the above spaces as $M\otimes E^{\bcdot}(\fg^*)$, although $M$ is a left module.
The translation of $M$ on the left of $E^{\bcdot}(\fg^*)$ is performed tacitly, via the sign rule.
\end{itemize}

Interestingly, on each of these spaces (see formula~\eqref{EandLambda}), a~codifferential, whose square vanishes, is defined: the expression for the former are used in \textsc{SuperLie}, see the description below; we are not certain how the codifferential for the latter is defined in general, but if $p=2$, the formula can be obtained from the formula for the codifferential on $E^{\bcdot}(\fg^*)\otimes M$ by replacing elements of $\Pi \fg^*$ with the corresponding elements of~$\fg^*$.
(Observe that the codifferential for~$\Lambda^{\bcdot}(\fg^*)$ given in \cite{Lc} is for some other objects, and only if $p=0$. It is not clear how these objects correspond to elements of~$\Lambda^{\bcdot}(\fg^*)\otimes M$ described here.) In particular, it means that the two cohomology spaces are in one-to-one correspondence when $p=2$.

\begin{OpenProblems}\index{Problem, open} How is the codifferential defined on elements of $\Lambda^{\bcdot}(\fg^*)\otimes M$, as it is described here, if $p\neq 2$? Are the cohomology spaces corresponding to $\Lambda^{\bcdot}(\fg^*)\otimes M$ and ${E^{\bcdot}(\fg^*)\otimes M}$ in one-to-one correspondence, when $p\neq 2$?
\end{OpenProblems}

In the package \textsc{SuperLie}, the codifferential on $M\otimes E^{\bcdot}(\fg^*)$ is defined as follows: let $e_1, \dots, e_K$ be a basis of $\fg$, and let $\alpha_1, \dots, \alpha_K$ be the dual basis of $\fg^*$. For any element
\[
{m\otimes \Pi\phi_1 \otimes \dots \otimes \Pi\phi_n\in M\otimes T^n(\Pi\fg^*)}, \qquad \text{where $m\in M$ and $\phi_1, \dots, \phi_n\in \fg^*$},
\]
 let $\langle m\otimes \Pi\phi_1 \otimes \dots \otimes \Pi\phi_n\rangle$
be the element of $M\otimes E^{\bcdot}(\fg^*)$ corresponding to its equivalence class. Then,
\begin{gather*}
d\langle m\otimes \Pi\phi_1 \otimes \dots \otimes \Pi\phi_n\rangle = \langle dm \otimes \Pi\phi_1 \otimes \dots \otimes \Pi\phi_n \nonumber \\
\hphantom{d\langle m\otimes \Pi\phi_1 \otimes \dots \otimes \Pi\phi_n\rangle =}{}
+\sum\limits_{1\leq i\leq n}(-1)^{p(m)+p(\Pi\phi_1)+\dots +p(\Pi\phi_{i-1})} m\otimes \Pi\phi_1 \otimes \cdots \nonumber\\
\hphantom{d\langle m\otimes \Pi\phi_1 \otimes \dots \otimes \Pi\phi_n\rangle =}{}
\otimes \Pi\phi_{i-1} \otimes d(\Pi\phi_i) \otimes \Pi\phi_{i+1} \otimes \dots \otimes \Pi\phi_{n}\rangle,%\label{2codiffE}
\end{gather*}
where
\[
dm = \sum\limits_{1\leq j\leq K} (-1)^{p(\Pi\alpha_j)(1+p(m))} (e_j m) \otimes \Pi\alpha_j,
\]
and
\begin{gather*}
d(\Pi\phi_i) = \sum\limits_{1\leq j<k\leq K} (-1)^{p(e_j)(1-p(e_k))} \phi_i([e_j,e_k]) \Pi\alpha_j \otimes \Pi\alpha_k\\
\hphantom{d(\Pi\phi_i) =}{} + \sum\limits_{1\leq j\leq K, \ e_j \ \text{is odd}} \phi_i\big(e_j^2\big) \Pi\alpha_j \otimes \Pi\alpha_j.
\end{gather*}

As a~very thorough referee observed, there exists linear bijections between the spaces of
super skew-symmetric and super anti-symmetric multilinear functions, e.g., one given by the assignment $f\mapsto f^a$, where
\[
f^a(X_1,\dots, X_n):=\begin{cases} (-1)^{p(X_1)+p(X_3)+\dots+ p(X_{n-1})}f(X_1,\dots, X_n)&\text{for $n$ even},\\
(-1)^{p(X_2) + p(X_4) + \dots + p(X_{n-1})}f(X_1,\dots, X_n)&\text{for $n$ odd}.
\end{cases}
\]
Therefore, one can make a complex of super skew-symmetric cochains isomorphic to that of super anti-symmetric cochains. Let us give the formulas for the corresponding codifferentials.

There are two additional definitions of cohomologies of $\fg$ with values in $M$, where $n$-cochains are considered $n$-linear functions on $\fg$ with values in $M$ with some condition on symmetry. In one of them, the cochains are \textit{super skew-symmetric}
functions, and the codifferential is given by the following expression (here, $\overline X=p(X)+1$, and the arguments with tilde should be ignored):
\begin{gather*}
d(m)(X)=(-1)^{p(m)p(X)}X(m) \qquad \text{for any $m\in M$ and $X\in\fg$}, \nonumber\\
d(f)(X_1,\dots, X_n) \nonumber\\
\qquad{}=\sum\limits_{i<j}(-1)^{\overline X_{1}+\dots+\overline X_{i}+\overline X_j(\overline X_{i+1}+\dots+\overline X_{j-1})}\!f\big(X_1, \dots, X_{i-1}, [X_i,X_j], X_{i+1}, \dots, \widetilde X_j, \dots, X_n\big) \nonumber\\
\qquad\quad{} +\sum(-1)^{p(X_j)(p(f)+n+1)+\overline X_j(\overline X_{1}+\dots+\overline X_{j-1})}X_j\big(f\big(X_1, \dots, \widetilde X_j, \dots, X_n\big)\big)\nonumber\\
\qquad\quad \text{for any $f\in E^{n-1}(\fg^*)\otimes M$ with $n>0$, and $X_1,\dots, X_n\in\fg$},%\label{2codiffA}
\end{gather*}
while in the other, the cochains are \textit{super anti-symmetric}
functions, and the codifferential is (the same notation conventions apply; these formulas are obtained from the standard formulas for the codifferential for Lie algebras by applying the sign rule):
\begin{gather*}
 d(m)(X)=(-1)^{p(m)p(X)}X(m)\text{~~for any $m\in M$ and $X\in\fg$},\nonumber\\
d(f)(X_1,\dots, X_n) \nonumber\\
\qquad{} = \sum\limits_{i<j}(-1)^{(p(X_j)(p(X_{i+1})+\dots+p(X_{j-1}))+j-1}\!f\big(X_1, {\dots}, X_{i-1}, [X_i,X_j], X_{i+1}, {\dots}, \widetilde X_j, {\dots}, X_n\big)\nonumber\\
\qquad\quad{}+\sum(-1)^{(p(X_j))(p(f)+p(X_{1})+\dots+p(X_{j-1}))+j-1}X_j\big(f\big(X_1, \dots, \widetilde X_j, \dots, X_n\big)\big)\nonumber\\
\qquad\quad{}\text{for any $f\in \Lambda^{n-1}(\fg^*)\otimes M$ with $n>0$, and $X_1,\dots, X_n\in\fg$}.%\label{2codiffS}
\end{gather*}

\begin{Hypothesis}
Supposedly, when $p=0$ or $n<p$, there is a natural isomorphism between $E^{n}(\fg^*)\otimes M$ and the space of super anti-symmetric $n$-linear functions on $\fg$ with values in $M$ such that the two corresponding definitions of codifferential agree; similarly, there is supposedly a natural isomorphism between $\Lambda^{n}(\fg^*)\otimes M$ and the space of super skew-symmetric $n$-linear functions on~$\fg$ with values in $M$ such that the two corresponding definitions of codifferential agree.
\end{Hypothesis}

\begin{OpenProblem}
How are the two isomorphisms defined?\index{Problem, open}
\end{OpenProblem}

Each of the four cohomology theories, call them $EH^{\bcdot}(\fg; M)$, $\Lambda H^{\bcdot}(\fg; M)$, and two additional theories: anti-symmetric $AH^{\bcdot}(\fg; M)$, and skew-symmetric $SH^{\bcdot}(\fg; M)$,\footnote{We will also denote the corresponding cochain spaces by $EC^{\bcdot}(\fg; M)$, $\Lambda C^{\bcdot}(\fg; M)$, $AC^{\bcdot}(\fg; M)$, and $SC^{\bcdot}(\fg; M)$.} respectively, has a~right to exist; the question is what do the corresponding invariants describe. %\ALNew{New part}
It is well-known (e.g., see~\cite{Fu}) that when $\fg$ is a Lie algebra, $AH^{2}(\fg; \fg)$ and $AH^{2}(\fg; \Kee)$ describe equivalence classes of non-trivial infinitesimal deformations and non-trivial central extensions of $\fg$, respectively. When $\fg$ is a Lie superalgebra and $p\neq 2,3$, these cohomology spaces can be used to describe those classes as well.

Let $M$ be $\fg$ or $\Kee$, then, assuming the isomorphisms between cohomology theories, $EH^{2}(\fg; M)$, $\Lambda H^{2}(\fg; M)$, and $SH^{2}(\fg; M)$ can also be used to describe those equivalence classes, but the description in terms of $AH^{2}(\fg; M)$ is the simplest one, since the bracket is a super anti-symmetric bilinear function on $\fg$.

When $p=2$ or $3$ and $\fg$ is a Lie superalgebra, $AH^{2}(\fg; M)$ and $SH^{2}(\fg; M)$ cannot be used to describe the deformations (with $M=\fg$) or central extensions (with $M$ trivial) of $\fg$, as these spaces of $2$-cohomologies are not isomorphic to the spaces of equivalence classes of non-trivial infinitesimal deformations and central extensions, respectively.
For example, when $p=3$, the ``deformations'' and ``central extensions'' described by $AH^{2}(\fg; M)$ do not necessarily satisfy condition \eqref{p=3}. When $p=2$, elements of $AC^{2}(\fg; M)$ cannot be used to describe changes to the algebraic structure of $\fg$ at all. Moreover, when $p=2$ and $\fg$ is not purely even, there is no natural isomorphism between $AC^{2}(\fg; M)$ or $SC^{2}(\fg; M)$ and the space of pairs $(b,s)$, where $b$ is an anti-symmetric bilinear map on $\fg$ with values in $M$ (i.e., a~``bracket''), and $s$ is a quadratic map on $\fg_\od$ with values in $M$ such that $b(x,y)=s(x+y)-s(x)-s(y)$ for any $x,y\in\fg_\od$ (i.e., a~``squaring'').

On the other hand, $EH^{2}(\fg; M)$ and $\Lambda H^{2}(\fg; M)$ can be used for these purposes, and this is what we do. We give results as produced by the package \textsc{SuperLie} in terms of elements of $M\otimes E^{n}(\fg^*)$. We also introduce yet another definition of cohomologies for $p=2$ in terms of functions. This definition gives $n$-cohomologies equivalent to $EH^{n}(\fg; M)$ and $\Lambda H^{n}(\fg; M)$ when $n\leq 2$, and we give an explicit isomorphism. This makes it easier to interpret our results as changes to the algebraic structure of $\fg$; the way to do it is described below.

We also present our results in terms of $EH^{2}(\fg; M)$ in characteristics other than $2$. Unfortunately, as we don't know the explicit form of isomorphism between $EH^{2}(\fg; M)$ and $AH^{2}(\fg; M)$ for Lie superalgebras $\fg$ when $p\neq 2$, we cannot say what exactly the infinitesimal deformations (or, in \cite{BGLL1}, central extensions) corresponding to the presented cochains of those Lie superalgebras are. For Lie algebras, the explicit form of the isomorphism is unclear only up to a coefficient, which is enough to find the spaces of non-trivial infinitesimal deformations and central extensions of $\fg$, knowing the corresponding spaces $EH^{2}(\fg; M)$.

\begin{OpenProblems}\index{Problem, open} When $p\neq 2$ and $\fg$ is a Lie superalgebra, is there an isomorphism between $EH^{2}(\fg; \fg)$ and the space of non-trivial infinitesimal deformations of $\fg$? Similarly, is there an isomorphism between $EH^{2}(\fg; \Kee)$ and the space of equivalence classes of non-trivial central extensions of $\fg$? If yes, what are their explicit forms? For the answer, see Section~\ref{AiP} ``added in proof''.
\end{OpenProblems}

For another definition of cohomology if $p=2$, see ``symmetric cohomology'', spoken about in~\cite{DzZ, Z} and studied in~\cite{LZ}.

\subsection[Cohomology of Lie superalgebras in degrees 1 and 2 for $p=2$\\ (written with the help of A.~Lebedev)]{Cohomology of Lie superalgebras in degrees 1 and 2 for $\boldsymbol{p=2}$\\ (written with the help of A.~Lebedev)}\label{sssZus}

For any Lie superalgebra $\fg$ in characteristic 2, its cochains needed to describe its infinitesimal deformations and central extensions cannot be described in terms of super anti-symmetric or super skew-symmetric functions on $\fg$. Here, we explain how they can be described in terms of other sorts of functions. In particular, we will need anti-symmetric bilinear and trilinear functions on $\fg$; by ``anti-symmetric'' we mean that the value of these bilinear and trilinear functions is $0$ whenever two of their arguments are equal, independently of the parity of these arguments.

This definition of cochains and codifferential produces $n$-cohomologies equivalent to those defined in terms of $E^{\bcdot}(\fg^*)\otimes M$ and $\Lambda^{\bcdot}(\fg^*)\otimes M$, when $n\leq 2$. We describe explicit correspondence between different types of cochains below.

 Recall that for every quadratic form $\mathfrak{q}$ with values in a space $M$, its polar form is the bilinear form with values in $M$ given by
\[
B_\mathfrak{q}(x,y):=\mathfrak{q}(x+y)-\mathfrak{q}(x)-\mathfrak{q}(y).
\]
In case of $p=2$, this bilinear form is anti-symmetric.
In order to define cochains of $\fg$ with values in $M$, we need a map
\begin{equation}\label{mapp}
\mathfrak p\colon \ \fg_\od \times \fg \tto M,
\end{equation}
with the following properties:
\begin{enumerate}\itemsep=0pt
\item[(i)] for any fixed $z\in \fg$, the map $x\mapsto \mathfrak{p}(x,z)$ is a quadratic form on $\fg_\od$;
\item[(ii)] for any fixed $x\in \fg_\od$, the map $z\mapsto \mathfrak{p}(x,z)$ is a linear map on $\fg$.
\end{enumerate}

Now, set
\begin{gather*}
XC^{-1}(\fg;M):= \{0\},\nonumber\\
XC^0(\fg;M):= M,\nonumber\\
XC^1(\fg;M):=\{c \mid \text{where $c\colon \fg \tto M$ is linear} \} \nonumber\\
XC^2(\fg;M):=\{(c,\mathfrak{g}) \mid \text{where $c\colon \fg\times\fg \tto M$ is anti-symmetric bilinear,}\nonumber\\
\hphantom{XC^2(\fg;M):=\{}{} \text{and $\mathfrak{q}\colon \fg_\od \tto M$ is a quadratic form such that $B_\mathfrak{q}=c|_{\fg_\od \wedge \fg_\od}$} \}, \nonumber\\
XC^3(\fg;M):=\{(c,\mathfrak{p}) \mid \text{where $c\colon \fg\times \fg \times \fg \tto M$ is anti-symmetric trilinear,}\nonumber\\
\hphantom{XC^3(\fg;M):=\{}{}
\text{$\mathfrak p$ is a map as in (\ref{mapp}) such that for all $x,y\in \fg_\od$ and $z\in \fg$ we have}\nonumber\\
\hphantom{XC^3(\fg;M):=\{}{}\mathfrak{p}(x+y,z)+\mathfrak{p}(x,z)+\mathfrak{p}(y,z)=c(x,y,z)\}. %\label{xcochains}
\end{gather*}

In later sections, we represent elements of $XC^2(\fg;M)$ as linear combinations of expressions of the form $v\otimes a\wedge b$, where $v\in M$ and $a,b\in \fg^*$. Such an expression represents the pair $(c,\mathfrak{g})$ such that
\begin{gather*}
c(x,y) = (a(x)b(y) + a(y)b(x))v \qquad \text{for all $x,y\in\fg$ and}\\
\mathfrak{q}(x)=a(x)b(x)v \qquad \text{for all $x\in\fg_\od$.}
\end{gather*}
More specifically, in the case $M=\fg$, we use expressions of the form~$u\otimes \widehat {v}\wedge \widehat {w}$, where $u$, $v$, $w$ are elements of a given basis of $\fg$, and $\widehat {u}$, $\widehat {v}$, $\widehat {w}$ are the corresponding elements of the dual basis of $\fg^*$. If $\{e_1, \dots, e_n\}$ is a basis of $\fg$ and $\{\widehat {e}_1, \dots, \widehat {e}_n\}$ is the dual basis of $\fg^*$, then the elements of the form $e_i\otimes \widehat {e}_j\wedge \widehat {e}_k$, where $1\leq i,j,k\leq n$ and either ``$j<k$'' or ``$j=k$ and $p(e_j)=\od$'', form a~basis of $XC^2(\fg;\fg)$.

Similarly, the elements of $XC^3(\fg;M)$ can be represented as linear combinations of expressions of the form
\[
v\otimes a\wedge b\wedge f, \qquad \text{where $v\in M$ and $a,b,f\in \fg^*$,}
\]
 representing the pair $(c,\mathfrak{p})$ such that
 \begin{gather*}
 c(x,y,z) = (a(x)b(y)f(z) + a(y)b(x) f(z)+ a(z)b(x) f(y)\\
 \hphantom{c(x,y,z) =}{}
 + a(z) b(y) f(x)+b(z)a(x)f(y)+b(z)a(y)f(x)) \qquad \text{for all $x,y,z\in\fg$}
 \end{gather*}
 and
 \[
 {\mathfrak p}(x,z)=(a(x) b(x) f(z)+a(z) b(x) f(x)+ a(x)b(z)f(x))v \qquad \text{for all $x\in \fg_\od$ and $z\in \fg$. }
 \]
In terms of $E^{\bcdot}(\fg)\otimes M$, the expressions $v\otimes a\wedge b$ and $v\otimes a\wedge b\wedge f$ correspond to the equivalence classes represented by $\Pi a\otimes \Pi b\otimes v$ and $\Pi a\otimes \Pi b\otimes \Pi f\otimes v$, respectively. In terms of $\Lambda^{\bcdot}(\fg)\otimes M$, $v\otimes a\wedge b$ and $v\otimes a\wedge b\wedge f$ correspond to the equivalence classes represented by $a\otimes b\otimes v$ and $a\otimes b\otimes f\otimes v$. The correspondence for $0$- and $1$-cochains is obvious. A direct check shows that under this correspondence, the operators ${\mathfrak d}^{-1}$, ${\mathfrak d}^{0}$, ${\mathfrak d}^{1}$, ${\mathfrak d}^{2}$, defined below, become the respective codifferentials on $E^{n}(\fg)\otimes M$ and $\Lambda^{n}(\fg)\otimes M$.

The elements of $XC^2(\fg; M)$ can be used to describe changes to the algebraic structure of $\fg$: a pair $(c,\mathfrak{g})$ describes a change to the bracket and the squaring given by
 \begin{gather*}
{}[x,y]' = [x,y] + c(x,y) k \qquad \text{for all $x,y\in\fg$}, \\
{} (x^2)' = x^2 + \mathfrak{q}(x) k \qquad \text{for all $x\in\fg_\od$},
 \end{gather*}
where $k$ is a multiplier whose role depends on the context. When $M=\fg$, this approach can be used to describe infinitesimal deformations of $\fg$, and in this case, $k$ is the parameter of deformation. A direct check shows that a pair $(c,\mathfrak{q})\in XC^2(\fg;\fg)$ describes an infinitesimal deformation if and only if ${\mathfrak d}^{2}(c,\mathfrak{q})=0$; this deformation is trivial if and only if there is ${c_1\in XC^1(\fg;\fg)}$ such that ${\mathfrak d}^{2}c_1 = (c,\mathfrak{q})$. Similarly, when $M=\Kee$, this approach can be used to describe central extensions of $\fg$, and in this case, $k$ is a~central element. A pair $(c,\mathfrak{q})\in XC^2(\fg;\Kee)$ describes a~central extension if and only if ${\mathfrak d}^{2}(c,\mathfrak{q})=0$; this central extension is trivial if and only if there is $c_1\in XC^1(\fg;\Kee)$ such that ${\mathfrak d}^{2}c_1 = (c,\mathfrak{q})$.

We define the \textit{codifferential}
\[
{\mathfrak d}^{-1}\colon \ XC^{-1}(\fg,M)=\{0\}\tto 0\in XC^0(\fg,M)=M.
\]

In our context, the \textit{codifferential} is the map
\[
{\mathfrak d}^0\colon \ XC^0(\fg,M)\tto XC^1(\fg,M), \qquad m \mapsto {\mathfrak d}^0(m),
\]
where ${\mathfrak d}^0(m)(x)=x\cdot m$.

Then, 1-cochains are just linear functions on $\fg$ with values in a $\fg$-module $M$. In our context, the \textit{codifferential} is the map
\[
{\mathfrak d}^1\colon \ XC^1(\fg,M)\tto XC^2(\fg,M), \qquad c \mapsto (dc, \mathfrak{q}),
\] where
\begin{gather*}
dc(x,y) = c([x,y]) + x \cdot c(y) + y \cdot c(x)\qquad \text{for all $x,y\in \fg$},\nonumber\\
\mathfrak{q}(x) = c(s(x))+x\cdot c(x) \qquad \text{for all $x\in \fg_\od$}.%\label{diffsuper}
\end{gather*}

\begin{Lemma}\label{d1}
The map ${\mathfrak d}^1$ is well-defined, i.e., $\operatorname{Im} {\mathfrak d}^1 \in XC^2(\fg,M)$.
\end{Lemma}

\begin{proof} We have
\begin{gather*}
\mathfrak{q}(x+y) + \mathfrak{q}(x) + \mathfrak{q}(y) = c( s(x+y)) + (x+y)\cdot c(x+y) +c(s(x)) \\
\hphantom{\mathfrak{q}(x+y) + \mathfrak{q}(x) + \mathfrak{q}(y) =}{} + x \cdot c(x) + c(s(y))+ y \cdot c(y) \\
\hphantom{\mathfrak{q}(x+y) + \mathfrak{q}(x) + \mathfrak{q}(y)}{} = c([x,y]) + y \cdot c(x) + x \cdot c(y)=dc(x,y).\tag*{\qed}
\end{gather*}\renewcommand{\qed}{}
\end{proof}

A 1-cocycle $c$ on $\fg$ with values in a $\fg$-module $M$ must satisfy the following conditions:
\begin{gather*}
%\label{Cond1}
x\cdot c(y)+y\cdot c(x)+c([x,y])= 0 \qquad \text{for all $x,y\in \fg$},\\
%\label{Cond2}
x\cdot c(x)+c(s(x)) = 0\qquad \text{for all $x\in \fg_\od$}.
\end{gather*}
Let the space of all 1-cocycles be denoted by $Z^1(\fg;M)$.

Let us introduce the \textit{codifferential}
\[
{\mathfrak d}^2\colon \ XC^2(\fg,M)\tto XC^3(\fg,M), \qquad (c,\mathfrak{q}) \mapsto (dc, \mathfrak{p}),
\]
where
\begin{gather}
dc(x,y,z) = x\cdot c(y,z)+c([x,y],z)+ \circlearrowleft (x,y,z) \qquad \text{for all $x,y,z\in \fg$},\nonumber\\
\mathfrak p(x,z) = x\cdot c(x,z)+z \cdot \mathfrak{q}(x)+ c(s(x),z)+ c([x,z],x]
\qquad \text{for all $x\in \fg_\od, z\in \fg$}.\label{diffsuper2}
\end{gather}
Note that this definition is consistent as the following lemma shows.

\begin{Lemma}\label{d2}
The map $\mathfrak{d}^2$ is well-defined, i.e., $\operatorname{Im} {\mathfrak d}^2 \in XC^3(\fg,M)$.
\end{Lemma}

\begin{proof}Indeed, using the fact that the polar form of $\mathfrak{q}$ is $c|_{\fg_\od \wedge \fg_\od}$ we get (for all $x,y \in \fg_\od$ and for all $z\in \fg$):
\begin{gather*}
{\mathfrak p}(x+y,z)+ {\mathfrak p}(x,z)+ {\mathfrak p}(y,z)\\
=(x+y)\cdot c(x+y,z)+ z \cdot \mathfrak{q}(x+y)+ c(s(x+y), z) + c([x+y,z], x+y]+ x\cdot c(x,z)\\
\quad{} + z \cdot \mathfrak{q}(x) + c(s(x),z) + c([x,z], x] + y \cdot c(y,z)+z \cdot \mathfrak{q}(y) + c(s(y),z) + c([y,z], y]\\
 = x\cdot c(y,z)+ y \cdot c(x,z)+ z \cdot c(x,y)+ c([x,y],z)+c([z,x],y) +c([y,z],x) =dc(x,y,z).\!\!\!\!\!\! \tag*{\qed}
\end{gather*}\renewcommand{\qed}{}
\end{proof}

Any 2-cocycle $(c,\mathfrak{q})$ of $\fg$ with values in $M$ must satisfy the following conditions:
\begin{gather}
\label{2-cocCond1}
0 = x\cdot c(y,z)+c([x,y],z)+ \circlearrowleft (x,y,z) \qquad \text{for all $x,y,z\in \fg$},\\
\label{2-cocCond2} 0 = x\cdot c(x,z)+z \cdot \mathfrak{q}(x)+ c(s(x),z)+ c([x,z],x] \qquad \text{for all $x\in \fg_\od$ and for all $z\in \fg$},\!\!\!\!
\end{gather}
Let the space of all 2-cocycles be denoted by $Z^2(\fg;M)$.

\begin{Theorem}\label{d2=0}
We have $\mathfrak{d}^2\circ \mathfrak{d}^1=0$.
\end{Theorem}

\begin{proof}To show that $d\circ d(c)=0$, see formula~\eqref{diffsuper2}, is a routine, see \cite{Fu,Wi}. Let us show that the map $\mathfrak p$ is identically zero. Indeed, we have
\begin{gather*}
 {\mathfrak p}(x,z) = x\cdot dc(x,z)+ z \cdot \mathfrak{q}(x)+ dc(s(x),z)+ dc([x,z],x] \\
\hphantom{{\mathfrak p}(x,z)}{}
 = x \cdot (x\cdot c(z)+z \cdot c(x)+c[x,z])+ z\cdot (x\cdot c(x)+c(s(x)))+s(x)\cdot c(z)\\
\hphantom{{\mathfrak p}(x,z)=}{} + z\cdot c(s(x))+ c([s(x),z])+ [x,z]\cdot c(x) +x\cdot c([x,z]) + c([[x,z], x])\\
\hphantom{{\mathfrak p}(x,z)}{} = s(x)\cdot c( z)+ [x,z]\cdot c(x)+x\cdot c([x,z]) +z\cdot c(s(x))) + s(x)\cdot c(z) \\
\hphantom{{\mathfrak p}(x,z)=}{} + z\cdot c(s(x))+ c([x,[z,x]])+ [x,z]\cdot c(x) +x\cdot c([x,z]) + c([[x,z], x])\\
\hphantom{{\mathfrak p}(x,z)}{} =0. \tag*{\qed}
\end{gather*}\renewcommand{\qed}{}
\end{proof}

We can now define deformation as follows:
\begin{gather*}
[x, y]_{{\rm new}} = [x, y]+ t c_1(x, y)+ t^2 c_2(x, y)+\cdots,\\
s_{{\rm new}}(x) = s(x) + t \mathfrak{q}_1(x)+t^2 \mathfrak{q}_2(x)+\cdots,
\end{gather*}
where $t$ is the parameter of deformation, the maps $c_i \colon\fg\times \fg\tto M$ are anti-symmetric bilinear maps, and the maps ${\mathfrak q}_i\colon \fg_\od \tto M$ are quadratic forms such that $B_{\mathfrak{q}_i}=c_i|_{\fg_\od \wedge \fg_\od}$.

The Jacobi identity for the new bracket, as well as the one involving the squaring, implies that the 2-tuple $(c_1, \mathfrak{q}_1)$ verifies the 1-cocycle conditions (\ref{2-cocCond1}) and (\ref{2-cocCond2}). Namely, $(c_1, \mathfrak{q}_1)\in Z^2(\fg,M)$. The equivalence class of deformations is captured by the cohomolgy space $H^2(\fg; M)$ as defined above.

\subsection{The Hochschild--Serre spectral sequence}\label{ssHSspecSeq} We mention the method of Hochschild--Serre spectral
sequences here because this powerful method enabled Kuznetsov and Chebochko compute (with bare hands) more, in some cases, than we computed using computer. We appealed to the HS spectral sequence just once.

Recall the first term of the Hochschild--Serre spectral
sequence corresponding to a~Lie subalgebra or subsuperalgebra
$\mathfrak{h}\subset\fg$ and abutting to the cohomology $H^{p+q}(\fg; M)$ we
are interested in, see \cite{Fu}:
\begin{equation*}%\label{SS:HSspec}
E_1^{a,b}=H^b\big(\mathfrak{h};\Eee^a(\fg/\mathfrak{h})^*\otimes M\big),
\end{equation*}
 where $\Eee^k$ is the operator of raising to the $k$th exterior power.

This gives an estimate of the number of cocycles to be tested from above, and also explicit forms of these cocycles.

In lemmas, all cocycles of weight $0$ that could represent a~nontrivial cohomology are
elements of $\big(\Eee^2(\fg_\od^*)\otimes \fg_\ev\big)^{\fg_\ev}$ or $(\fg_\od^* \otimes \fg_\od)^{\fg_\ev}$.
\begin{itemize}\itemsep=0pt
\item Kuznetsov and Chebochko \cite{Ch1, KuCh} proved that, in the cases they
considered, there are no cocycles of weight~0. There are, however, cases where such cocycles are manifest:
if the Cartan matrix has a~parameter.
\end{itemize}

\textit{We were unable to pinpoint the place in the arguments due to Kuznetsov and Chebochko where these arguments $($proving that there should be no cohomology$)$ do not work; we just give examples where such cohomology is non-zero}.

\subsection{Restricted cohomology of restricted Lie algebras and Lie superalgebras}\label{ssResH}
This subsection contains problems tackled in \cite{Ev, YChC} but left mostly open.

In \cite{Hos}, the notion of restricted cohomology $H_{\rm res}^i(\fg; M)$ of the restricted Lie algebra $\fg$ with coefficients in the restricted $\fg$-module $M$ was introduced. If $p=2$, there are, however, \textit{several} contenders for the role of restricted Lie algebra and a restricted module over it, see \cite{BLLS1} whereas only the ``traditional'' definition is considered in \cite{Hos} and \cite{Ev}. The only examples considered in the literature are filiform Lie algebras. The following general statements on comparison of restricted and ``usual'' cohomology are proved in \cite{Hos}:
\begin{enumerate}\itemsep=0pt
\item[0)] $H_{\rm res}^0(\fg; M)=H^0(\fg; M)$;
\item[1)] $H_{\rm res}^1(\fg; M)\subset H^1(\fg; M)$;
\item[2)] there is a map $H_{\rm res}^2(\fg; M)\tto H^2(\fg; M)$ with a possibly non-zero kernel.
\end{enumerate}

Interpretations of the \textit{restricted} cohomology in degrees 1 and 2 with coefficients in trivial and adjoint modules are the same as those of the ``usual'' cohomology, but are applicable to only \textit{restricted} extensions, and derivations.

In \cite{YChC}, the definitions and results of \cite{Hos} are superised; observe that the paper \cite{YChC} ignores various new, unconventional, restrictednesses of Lie algebras and Lie superalgebras in characteristic 2 introduced in \cite{BLLS1}. The above remarks lead to the following

\begin{OpenProblem}\index{Problem, open}
In characteristic 2, define restricted universal enveloping algebra and corresponding restricted cohomology for various restrictednesses of Lie algebras, Lie superalgebras, and modules over them. In all these cases, compute the kernel and the image of the map $H_{\rm res}^2(\fg; M)\tto H^2(\fg; M)$.
\end{OpenProblem}

\section{Notation in ``Results" sections}\label{sNotResult}
Recall the definition of the Lie (super)algebras with indecomposable Cartan matrices, see~\cite{BGL2}.

Observe that, studying deformations, it does not matter which of
several Cartan matrices that a~given algebra $\fg(A)$ has we
take: the simplest incarnation will do. Computations were performed,
with long gaps, during several years, and the final results are
expressed differently in some cases; all is clear and we did not
waste time for unification of final expressions.

Recall that it is meaningless to ask ``what is the parity of a~given 1-cochain?'': it depends on the problem considered. Having selected a basis of~$\fg$, for any basis element $x\in\fg$, the cochain~$\widehat x$ can be considered as the element of the dual space~$\fg^*$ if there is no need to multiply it by any other cochain (or if we consider \textit{super symmetric} product, \textit{not exterior} one, as when dealing with analogs of metrics), whereas if we consider $\widehat x$ in the algebra of cochains using the $\wedge$ product, then $\widehat x\in\Pi(\fg^*)$. The notation is chosen to replace~$dx$ since the symbol~$d$ is overused.

We denote the exterior product of cochains by wedge, the powers are denoted $(\widehat x)^{\wedge n}$.

The cocycles below, except for those of $\mathfrak{o}_I$ which is a~non-split Lie algebra, are indexed by
their weight; the superscript
enumerates linearly independent cocycles of the same weight, if there
are several of them.

We let the positive Chevalley \textit{generators} be of degree 1, and denote the
elements of Chevalley \textit{basis} they generate, by $x_i$, we denote the corresponding negative
basis elements by $y_i$; we set $h_i:=[x_i, y_i]$ for the \textit{generators} $x_i$ and $y_i$ of degree $\pm 1$ only---this is the \textit{principal} grading.

\subsection{Convention}\label{ssConven} By abuse of language, we say ``$\fg$ is rigid'' meaning ``$H^2(\fg;\fg)=0$''.

Since our algebras are symmetric with respect to the change of the
sign of the roots, it suffices to consider cocycles of only
nonpositive degrees, so we do not list cocycles of positive degrees.
By abuse of notation and to save trees, the phrase ``$H^2(\fg;\fg)$ is spanned by the cocycles~$c_i$'' means ``$H^2(\fg;\fg)$ is spanned by the classes $[c_i]$ of cocycles $c_i$; moreover, cocycles of positive degrees, symmetric to those of negative degrees, are assumed, but not given''.

The even and odd basis elements of the Lie superalgebras considered
are grouped so that the even are on the left of a~bar, the odd
ones on the right, see, e.g., formula~\eqref{nolab}.

The cocycles in ``Results'' sections is found with help of \textsc{SuperLie}
package, in which notation are at variance with formula~\eqref{EandLambda}. Namely, when (for $p\neq 2$) interpreting $f_1\wedge f_2\otimes m\in H^2(\fg; M)$ as an $M$-valued anti-symmetric bilinear function on $\fg$, one has to move $m$ to the left using the sign rule and then apply $f_1\wedge f_2$ to an element of $\fg\times \fg$.

Hereafter, we underline the \underline{odd} cocycle of definite weight; each of them is
integrable. Deforms corresponding to every odd cocycle should be
considered with an odd parameter.

\section[Results: Lie superalgebras for $p=5$]{Results: Lie superalgebras for $\boldsymbol{p=5}$}\label{SResSuperP=5}

\subsection[Simple (relatives of) Lie superalgebras with Cartan matrix for $p\geq 5$ are rigid,\\ except for $\mathfrak{osp}(4\vert 2)$ and $\mathfrak{brj}(2;5)$]{Simple (relatives of) Lie superalgebras with Cartan matrix for $\boldsymbol{p\geq 5}$ \\ are rigid, except for $\boldsymbol{\mathfrak{osp}(4\vert 2)}$ and $\boldsymbol{\mathfrak{brj}(2;5)}$}\label{ssPgeq5}

We checked this for $p=5,7,11$ for $\mathfrak{ag}(2)$, $\mathfrak{ab(}3)$, and for $p=5,7$ for the series: $\mathfrak{sl}(m|n)$ for $1\leq m< n\leq 4$ with $m+n>2$, and $\mathfrak{psl}(n|n)$ for $2\leq n\leq 4$. The orthosymplectic algebras are considered for $p\geq 3$ in Section~\ref{Lrigid}.

\begin{Lemma}\label{fbrj(2;5)} For $\fg=\mathfrak{brj}(2;5)$, we consider the Cartan matrix
\[
\begin{pmatrix}
\hphantom{-}0&-1\\
 -2&\hphantom{-}1
\end{pmatrix}
\]
and the basis even $|$ odd
\begin{gather}
\mid x_{1}, \ x_{2},\nonumber\\
x_3= [x_{1}, x_{2} ], \ x_4= [x_{2},\,x_{2} ] \nonumber\mid \\
\mid x_5= [x_{2}, [x_{1},x_{2} ] ],\nonumber \\
 x_6= [ [x_{1},x_{2} ], [x_{2},\,x_{2} ] ]\mid \nonumber\\
\mid x_7= [ [x_{1},x_{2} ], [x_{2}, [x_{1},\,x_{2} ] ] ], \ x_8= [ [x_{2},x_{2} ], [x_{2}, [x_{1},x_{2} ] ] ],\nonumber\\
 x_9= [ [x_{1},x_{2} ], [ [x_{1},x_{2} ], [x_{2},x_{2} ]  ] ]\mid \nonumber\\
\mid x_{10}= [ [x_{2}, [x_{1},x_{2} ] ], [ [x_{1},x_{2} ], [x_{2},x_{2} ] ] ].\label{nolab}
\end{gather}
Then, $H^2(\fg;\fg)$ is spanned by the cocycles
\begin{gather*}
\underline{c_{-15}}= 2y_1\otimes \widehat x_{10}\wedge \widehat
x_{10}+3y_3\otimes \widehat x_9 \wedge \widehat x_{10}+y_5\otimes \widehat x_7 \wedge
\widehat x_{10} +3y_6\otimes \widehat x_7 \wedge \widehat x_9+3y_7\otimes \widehat x_5 \wedge
\widehat x_{10}\nonumber\\
\hphantom{\underline{c_{-15}}=}{}
+4y_7\otimes \widehat x_6 \wedge \widehat x_9+y_7\otimes \widehat x _7 \wedge \widehat x_8
+2y_8\otimes \widehat x_7 \wedge \widehat x_7+4y_9\otimes \widehat x_3 \wedge
\widehat x_{10}+y_9\otimes \widehat x_6 \wedge \widehat x_7 \nonumber\\
\hphantom{\underline{c_{-15}}=}{}+4y_{10}\otimes \widehat x_1 \wedge \widehat x_{10}+4y_{10}\otimes \widehat x_3 \wedge
\widehat x_9+2y_{10}\otimes \widehat x_5 \wedge
\widehat x_7,\nonumber\\
\underline{c_{-5}}=3x_{1}\otimes \widehat x_{2}{\wedge} \widehat x_{8} +
4x_{1}\otimes \widehat x_{4}{\wedge} \widehat x_{6} + 3x_{1}\otimes \widehat x_{10} {\wedge}
\widehat y_{1} + 4x_{3}\otimes \widehat x_{4}{\wedge} \widehat x_{8} + 3x_{7}\otimes \widehat x_{8}
{\wedge} \widehat x_{8}
\nonumber\\
\hphantom{\underline{c_{-5}}=}{}
+y_{2}\otimes \widehat x_{8}{\wedge} \widehat y_{1} + 4y_{4}\otimes \widehat x_{6} {\wedge}
\widehat y_{1} + 2y_{4}\otimes \widehat x_{8}{\wedge} \widehat y_{3} + 2y_{6}\otimes
\widehat x_{4}{\wedge} \widehat y_{1} + 2y_{8}\otimes \widehat x_{2} {\wedge} \widehat y_{1}
\nonumber\\
\hphantom{\underline{c_{-5}}=}{}
+ 3y_{8}\otimes \widehat x_{4}{\wedge} \widehat y_{3} + y_{8}\otimes \widehat x_{8} {\wedge}
\widehat y_{7} + y_{10}\otimes \widehat y_{1}{\wedge} \widehat y_{1}.%\label{eqbrj2,5}
\end{gather*}
\end{Lemma}

\begin{Lemma} \label{el55} The Lie superalgebra $\fel(5; 5)$ is rigid.
\end{Lemma}

\section[Results: Lie algebras for $p=3$]{Results: Lie algebras for $\boldsymbol{p=3}$}\label{SResP>2}

\subsection[The simple (relatives of) Lie algebras with Cartan matrix for $p>3$ are rigid]{The simple (relatives of) Lie algebras with Cartan matrix\\ for $\boldsymbol{p>3}$ are rigid}\label{ssPgeq3}

This is proved in \cite{Ru}.
In what follows we list all deforms of simple (relatives of) Lie algebras with Cartan matrix for $p=3$, not considered earlier, cf.\ \eqref{NC}, but give the next result for completeness.

\subsection[The deforms of $\mathfrak{o}(5)$]{The deforms of $\boldsymbol{\mathfrak{o}(5)}$}\label{LeBr2a} For a~detailed study of deforms of $\mathfrak{o}(5)$ in characteristic~$3$, the first of which was discovered by Rudakov, see \cite{BLW}, where earlier claims were corrected. There are two classes of non-isomorphic deforms of $\mathfrak{o}(5)\simeq \mathfrak{br}(2; -1)$:
\begin{enumerate}\itemsep=0pt
\item[1)] A~parametric family $\mathfrak{br}(2; \varepsilon)$,
where $\varepsilon\neq 0$, with the Cartan matrix
\begin{equation*}%\label{br2a}
\begin{pmatrix}
\hphantom{-}2&-1\\
-2&1-\varepsilon\end{pmatrix}\qquad\text{and the basis} \ \ x_1, \ x_2, \
x_3=[x_1, x_2],\ x_4=-\operatorname{ad}_{x_2}^2(x_1),
\end{equation*}
and where
\begin{equation*}%\label{Lepsiso}
\mathfrak{br}(2; \varepsilon)\simeq
\mathfrak{br}(2; \varepsilon') \qquad \text{if and only if} \quad \varepsilon\varepsilon'=1\quad\text{(for
$\varepsilon\neq \varepsilon'$)}.
\end{equation*}
\item[2)] An exceptional simple Lie algebra $\mathfrak{L}(2,2)$. Recall the description of $\mathfrak{L}(2,2)$, see \cite[Proposition~3.2]{BLW}.
The contact bracket of two divided power polynomials $f,g\in \mathcal{O}(p,q,t;\un)$
is defined to be
\begin{equation*}%\label{cb}
\{f,g\}_{\rm k.b.}=\triangle f\cdot\partial_t g - \partial_t f\cdot\triangle g
 +\partial_p f\cdot\partial_q g - \partial_q f\cdot\partial_p g,
\end{equation*}
where $\triangle f=2f - p\partial_p f - q\partial_q f$.
\end{enumerate}

Let $\alpha$ and $\beta$ be the
simple roots of $\mathfrak{o}(5)$ and $E_\gamma$ be the root vector corresponding
to root $\gamma$. Then, a~basis of $\mathfrak{br}(2; \varepsilon)$ is expressed in
terms of generating functions of $\mathfrak{k}(3; \One)$ as
follows:
$$%\label{tab}
\renewcommand{\arraystretch}{1.4}
\begin{tabular}{|c|l|} \hline
$\deg$&the generator with its weight $\sim$ its generating function (= the Chevalley basis vector) \\
\hline \hline
$-2$&$E_{-2\alpha -\beta}=\{E_{-\alpha}, E_{-\alpha-\beta}\}_{\rm k.b.}\sim 1(=y_4)$\\
\hline
$-1$& $E_{-\alpha}\sim p(=y_2)$; $E_{-\alpha-\beta}=\{E_{-\beta}, E_{-\alpha}\}_{\rm k.b.}\sim q(=y_3)$ \\
\hline $\hphantom{-}0$& $H_{\alpha}\sim 2 \varepsilon t+ pq(=h_2)$; $H_{\beta}\sim -pq(=h_1)$; $E_{\beta}\sim p^2(=y_1)$; $E_{-\beta}\sim -q^2(=x_1)$ \\
\hline
$\hphantom{-}1$& $E_{\alpha}\sim -(1+\varepsilon)pq^2+ \varepsilon
qt(=x_2)$; $E_{\alpha+\beta}=\{E_\beta, E_\alpha\}_{\rm k.b.}\sim (1+\varepsilon)p^2q+\varepsilon p t(=x_3)$ \\
\hline $\hphantom{-}2$& $E_{2\alpha+\beta}=\{E_\alpha,
E_{\alpha+\beta}\}_{\rm k.b.}\sim \varepsilon(1+\varepsilon)p^2q^2+ \varepsilon^2t^2(=x_4)$ \\
\hline
\end{tabular}
$$

The Lie algebra, called $\mathfrak{L}(2,2)$ for (inessential for us) reasons explained in \cite{BLW}, is the following
deform of the bracket $[\cdot,\cdot]$ of $\mathfrak{o}(5)$
\begin{gather*}%\label{L-1bracket}
[\cdot,\cdot]_\lambda = [\cdot,\cdot]-c, \qquad \text{where $c= x_1\otimes \widehat y_3 \wedge \widehat y_4-
 x_3\otimes \widehat y_1 \wedge \widehat y_4 +
 x_4\otimes \widehat y_1 \wedge \widehat y_3$}.
\end{gather*}

\begin{Lemma}\label{LeBr3a} For $\fg=\mathfrak{br}(3)$, we consider the
Cartan matrix
\[
\begin{pmatrix}
\hphantom{-}2&-1&\hphantom{-}0\\
-1&\hphantom{-}2&-1\\
\hphantom{-}0&-1&\hphantom{-}\ev\end{pmatrix}
\]
and the basis
\begin{gather*}
x_1,\ x_2,\ x_3,\ x_4=[x_1,x_2],\ x_5=[x_2,x_3], \
x_6=[x_3,[x_1, x_2]],\ x_7=[x_3,[x_2,x_3]],\nonumber\\
x_8=[x_3,[x_3,[x_1,x_2]]],\
x_9=[[x_2,x_3], [x_3,[x_1,x_2]]], \
x_{10}=[[x_3,[x_1,x_2]], [x_3,[x_2,x_3]]],;\nonumber\\
x_{11}=[[x_3,[x_2,x_3]], [x_3,[x_3,[x_1,x_2]]]],\
x_{12}=[[x_3,[x_2,x_3]], [[x_2,x_3],
[x_3,[x_1,x_2]]]],\nonumber\\
x_{13}=[[x_3,[x_3,[x_1,x_2]]], [[x_2,x_3],
[x_3,[x_1,x_2]]]].%\label{br3}
\end{gather*}
Then, $H^2(\fg; \fg)$ is spanned by the cocycles \eqref{eqbr3}:
\begin{gather}
c_{-18}= 2y_{3}\otimes \widehat x_{12}
{\wedge} \widehat x_{13} + 2y_{5}\otimes \widehat x_{11} {\wedge} \widehat x_{13} +
2y_{6}\otimes \widehat x_{11} {\wedge} \widehat x_{12} + y_{7}\otimes \widehat x_{10}
{\wedge} \widehat x_{13} + y_{8}\otimes
 \widehat x_{10}
{\wedge} \widehat x_{12}
\nonumber\\
\hphantom{c_{-18}=}{}
+ y_{9}\otimes \widehat x_{10} {\wedge} \widehat x_{11} + 2y_{10}\otimes \widehat x_{7}
{\wedge} \widehat x_{13} + y_{10}\otimes \widehat x_{8} {\wedge} \widehat x_{12} +
2y_{10}\otimes \widehat x_{9} {\wedge} \widehat x_{11} \nonumber\\
\hphantom{c_{-18}=}{}
+ 2y_{11}\otimes \widehat x_{5}
{\wedge} \widehat x_{13}+
 y_{11}\otimes \widehat x_{6} {\wedge} \widehat x_{12} + 2y_{11}\otimes \widehat x_{9}
{\wedge} \widehat x_{10} + 2y_{12}\otimes \widehat x_{3}{\wedge} \widehat x_{13} +
y_{12}\otimes \widehat x_{6} {\wedge} \widehat x_{11} \nonumber\\
\hphantom{c_{-18}=}{}
+ 2y_{12}\otimes \widehat x_{8}
{\wedge} \widehat x_{10}+
 2y_{13}\otimes \widehat x_{3}{\wedge} \widehat x_{12} + y_{13}\otimes \widehat x_{5}
{\wedge} \widehat x_{11} + 2y_{13}\otimes \widehat x_{7}
{\wedge} \widehat x_{10},\nonumber\\
c_{-9}= 2x_{3}\otimes \widehat x_{1}{\wedge} \widehat x_{13} + x_{5}\otimes \widehat x_{4}
{\wedge} \widehat x_{13} + x_{7}\otimes \widehat x_{6} {\wedge} \widehat x_{13} +
y_{1}\otimes \widehat x_{4}{\wedge} \widehat x_{10} + 2y_{1}\otimes \widehat x_{6} {\wedge}
\widehat x_{9}\nonumber
\\
\hphantom{c_{-9}=}{}
+ 2y_{1}\otimes \widehat x_{13} {\wedge} \widehat y_{3} + y_{4}\otimes \widehat x_{1}
{\wedge} \widehat x_{10} + y_{4}\otimes \widehat x_{6} {\wedge} \widehat x_{8} +
y_{4}\otimes \widehat x_{13} {\wedge} \widehat y_{5} + y_{6}\otimes \widehat x_{1} {\wedge}
\widehat x_{9}\nonumber\\
\hphantom{c_{-9}=}{}+ y_{6}\otimes \widehat x_{4} {\wedge} \widehat x_{8} + 2y_{6}\otimes \widehat x_{13}
{\wedge} \widehat y_{7} + 2y_{8}\otimes \widehat x_{4}{\wedge} \widehat x_{6} + y_{9}\otimes
\widehat x_{1} {\wedge} \widehat x_{6} + 2y_{10}\otimes \widehat x_{1}
{\wedge} \widehat x_{4}\nonumber\\
\hphantom{c_{-9}=}{}+ y_{13}\otimes \widehat x_{1}{\wedge} \widehat y_{3} + y_{13}\otimes \widehat x_{4} {\wedge}
\widehat y_{5} + y_{13}\otimes \widehat x_{6}
{\wedge} \widehat y_{7},\nonumber\\
c_{-6}= 2x_{1}\otimes \widehat x_{2}{\wedge} \widehat x_{10} + 2x_{1}\otimes \widehat x_{5}
{\wedge} \widehat x_{9} + 2x_{1}\otimes \widehat x_{12} {\wedge} \widehat y_{3} +
x_{3}\otimes \widehat x_{12}{\wedge} \widehat y_{1} + 2x_{6}\otimes \widehat x_{2} {\wedge}
\widehat x_{12}
\nonumber\\
\hphantom{c_{-6}=}{}
+ x_{8}\otimes \widehat x_{5} {\wedge} \widehat x_{12} + 2y_{2}\otimes \widehat x_{5}
{\wedge} \widehat x_{7} + 2y_{2}\otimes \widehat x_{10} {\wedge} \widehat y_{1} +
2y_{2}\otimes \widehat x_{12} {\wedge} \widehat y_{6} + 2y_{5}\otimes \widehat x_{2}
{\wedge} \widehat x_{7}
\nonumber\\
\hphantom{c_{-6}=}{} + y_{5}\otimes \widehat x_{9} {\wedge} \widehat y_{1} + 2y_{5}\otimes \widehat x_{12}
{\wedge} \widehat y_{8} + y_{7}\otimes \widehat x_{2}{\wedge} \widehat x_{5} + 2y_{9}\otimes
\widehat x_{5} {\wedge} \widehat y_{1} + 2y_{10}\otimes \widehat x_{2}
{\wedge} \widehat y_{1}\nonumber\\
\hphantom{c_{-6}=}{}+ 2y_{12}\otimes \widehat x_{2}{\wedge} \widehat y_{6} + y_{12}\otimes \widehat x_{5}
{\wedge} \widehat y_{8} + 2y_{12}\otimes \widehat y_{1} {\wedge} \widehat y_{3},\nonumber\\
c_{-3}= 2x_{2}\otimes \widehat x_{3}{\wedge} \widehat x_{7} + x_{2}\otimes \widehat x_{10}
{\wedge} \widehat y_{4} + 2x_{2}\otimes \widehat x_{11} {\wedge} \widehat y_{6}
 + x_{4}\otimes \widehat x_{3}{\wedge} \widehat x_{8}
+ x_{4}\otimes \widehat x_{10} {\wedge} \widehat y_{2}
\nonumber\\
\hphantom{c_{-3}=}{}
+ x_{4}\otimes \widehat x_{11} {\wedge} \widehat y_{5} + 2x_{5}\otimes \widehat x_{11}
{\wedge} \widehat y_{4} + x_{6}\otimes \widehat x_{11} {\wedge} \widehat y_{2} +
x_{9}\otimes \widehat x_{3} {\wedge} \widehat x_{11} + 2y_{3}\otimes \widehat x_{7} {\wedge}
\widehat y_{2}\nonumber\\
\hphantom{c_{-3}=}{}+ 2y_{3}\otimes \widehat x_{8} {\wedge} \widehat y_{4} + 2y_{3}\otimes \widehat x_{11}
{\wedge} \widehat y_{9} + 2y_{7}\otimes \widehat x_{3}{\wedge} \widehat y_{2} + y_{8}\otimes
\widehat x_{3} {\wedge} \widehat y_{4} + y_{10}\otimes \widehat y_{2} {\wedge}
\widehat y_{4}\nonumber\\
\hphantom{c_{-3}=}{}+ y_{11}\otimes \widehat x_{3}{\wedge} \widehat y_{9} + y_{11}\otimes \widehat y_{2} {\wedge}
\widehat y_{6} + y_{11}\otimes \widehat y_{4} {\wedge} \widehat y_{5}.\label{eqbr3}
\end{gather}
\end{Lemma}

We find this result a~bit
too plentiful, as in \cite{Ch1, Vi1,Vi2, Vi1a,Vi3}. It has to be dealt with as in~\cite{Ch1, KuCh}, see Open problem~\ref{OP}\,(5).

\section[Results: Lie superalgebras for $p=3$]{Results: Lie superalgebras for $\boldsymbol{p=3}$}\label{SResSuperP=3}

\begin{Lemma}\label{brj(2;3)} For $\fg=\mathfrak{brj}(2;3)$, we consider the Cartan matrix
\[
\begin{pmatrix}
\hphantom{-}0&-1\\
 -2&\hphantom{-}1
\end{pmatrix}
\]
and the basis
even $|$ odd
\begin{gather*}
\mid x_{1},\ x_{2},\\
 x_3= [x_{1}, x_{2} ]\mid\\
 x_4= [x_{2},\,x_{2} ] \mid\\
\mid  x_5= [x_{2}, [x_{1},x_{2} ] ] , \\
   x_6= [ [x_{1},\,x_{2} ], [x_{2},\,x_{2} ] ] \mid \nonumber\\
\mid x_7= [ [x_{2}, x_{2} ], [x_{2}, [x_{1},x_{2} ] ] ], \\
 x_8= [ [x_{1},x_{2} ], [ [x_{1},x_{2} ], [x_{2},x_{2} ] ] ]\mid%\label{eqbrj(2;3)}
\end{gather*}
Then, $H^2(\fg;\fg)$ is spanned by the cocycles
\begin{gather*} c_{-12}=
y_{2}\otimes \widehat x_{7}{\wedge}\widehat x_{8} + y_{3}\otimes \widehat x_{7}
{\wedge}\widehat x_{7} + 2y_{4}\otimes \widehat x_{6}{\wedge}\widehat x_{8} + y_{5}\otimes
\widehat x_{6}{\wedge}\widehat x_{7} + y_{6}\otimes \widehat x_{4}{\wedge}\widehat x_{8}
\nonumber\\
\hphantom{c_{-12}=}{}
+ 2y_{6}\otimes \widehat x_{5}{\wedge}\widehat x_{7} + 2y_{7}\otimes \widehat x_{2}{\wedge}
\widehat x_{8} + 2y_{7}\otimes \widehat x_{3}{\wedge} \widehat x_{7} + 2y_{7}\otimes
\widehat x_{5}{\wedge} \widehat x_{6} + y_{8}\otimes \widehat x_{2} {\wedge} \widehat x_{7}\nonumber\\
\hphantom{c_{-12}=}{}
 +
y_{8}\otimes
 \widehat x_{4}{\wedge}
\widehat x_{6},\nonumber\\
c_{-6}=2x_{2}\otimes \widehat x_{1}{\wedge} \widehat x_{8} + x_{4}\otimes
 \widehat x_{3}{\wedge}
\widehat x_{8} + 2y_{1}\otimes
 \widehat x_{1}{\wedge}
\widehat x_{6} + 2y_{1}\otimes
 \widehat x_{3}{\wedge}
\widehat x_{5} + 2y_{1}\otimes \widehat x_{8}{\wedge} \widehat y_{2}
\nonumber\\
\hphantom{c_{-6}=}{}
+ 2y_{3}\otimes \widehat x_{1}{\wedge} \widehat x_{5} + y_{3}\otimes
 \widehat x_{8}{\wedge}
\widehat y_{4} + y_{5}\otimes \widehat x_{1}{\wedge} \widehat x_{3} + y_{6}\otimes \widehat x_{1}
{\wedge} \widehat x_{1} + y_{8}\otimes
 \widehat x_{1}{\wedge}
\widehat y_{2} \nonumber\\
\hphantom{c_{-6}=}{}
+y_{8}\otimes
 \widehat x_{3}{\wedge}
\widehat y_{4}.%\label{eqbrj2,3}
\end{gather*}
\end{Lemma}

\begin{Lemma}\label{g(1,6)} For $\fg=\fg(1,6)$, we consider the Cartan matrix
\[
 \begin{pmatrix}
 \hphantom{-}2&-1&\hphantom{-}0 \\
 -1&\hphantom{-}1&-1 \\
 \hphantom{-}0&-1&\hphantom{-}0
\end{pmatrix}
\]
and the basis even $|$ odd
\begin{gather*}
x_{1},\mid x_{2},\ x_{3}, \\
 x_5= [x_{2}, x_{2} ], \mid x_4=
 [x_{1}, x_{2} ], \\
  x_6= [x_{2}, x_{3} ],\ x_7= [x_{2}, [x_{1}, x_{2} ] ]\mid \nonumber\\
 x_8= [x_{3}, [x_{1}, x_{2} ] ],\mid x_9= [x_{3}, [x_{2}, x_{2} ] ], \\
x_{10}= [ [x_{1}, x_{2} ], [x_{1}, x_{2} ] ], \mid
x_{11}= [ [x_{1}, x_{2} ], [x_{2}, x_{3} ] ]\nonumber\\
\mid x_{12}= [ [x_{1}, x_{2} ], [x_{3},  [x_{1}, x_{2} ] ] ], \\
 x_{13}= [
 [x_{2}, x_{3} ], [x_{2},
 [x_{1}, x_{2} ] ]
 ]\mid \nonumber\\
 x_{14}= [ [x_{2}, [x_{1}, x_{2} ] ] , [x_{3}, [x_{1}, x_{2} ] ] ]\mid\\
\mid  x_{15}= [ [x_{3}, [x_{2},x_{2} ] ] , [ [x_{1},\,x_{2} ], [x_{1},x_{2} ] ] ],\nonumber\\
 x_{16}= [ [x_{3}, [x_{2},x_{2} ] ] , [ [x_{1},x_{2} ], [x_{3}, [x_{1},x_{2} ] ] ] ]\mid %\label{1,6}
\end{gather*}
Then, $H^2(\fg; \fg)$ is spanned by the cocycles:
\begin{gather*}
\underline{c_{-12}}= 2x_{2}\otimes \widehat x_{13}
{\wedge} \widehat x_{16} + 2x_{4}\otimes \widehat x_{14} {\wedge} \widehat x_{16} +
y_{2}\otimes \widehat x_{8} {\wedge} \widehat x_{16} + 2y_{3}\otimes
\widehat x_{13}{\wedge} \widehat x_{14} + y_{4}\otimes \widehat x_{6} {\wedge} \widehat x_{16} \nonumber\\
\hphantom{\underline{c_{-12}}=}{}
+ y_{6}\otimes \widehat x_{4} {\wedge} \widehat x_{16} + y_{6}\otimes \widehat x_{8} {\wedge}
\widehat x_{15} + 2y_{6}\otimes \widehat x_{11} {\wedge} \widehat x_{14} + 2y_{6}\otimes
\widehat x_{12} {\wedge} \widehat x_{13} + y_{8}\otimes \widehat x_{2}
{\wedge} \widehat x_{16} \nonumber\\
\hphantom{\underline{c_{-12}}=}{}
+ y_{8}\otimes \widehat x_{6} {\wedge} \widehat x_{15} + y_{8}\otimes \widehat x_{9} {\wedge}
\widehat x_{14} + y_{8}\otimes \widehat x_{11}{\wedge} \widehat x_{13} + y_{9}\otimes \widehat x_{8}
{\wedge} \widehat x_{14} + y_{11}\otimes \widehat x_{6}
{\wedge} \widehat x_{14} \nonumber\\
\hphantom{\underline{c_{-12}}=}{}
+ y_{11}\otimes \widehat x_{8}{\wedge} \widehat x_{13} + y_{12}\otimes \widehat x_{6}
{\wedge} \widehat x_{13} + 2y_{13}\otimes \widehat x_{3} {\wedge} \widehat x_{14} +
2y_{13}\otimes \widehat x_{6}{\wedge} \widehat x_{12} + 2y_{13}\otimes \widehat x_{8}
{\wedge} \widehat x_{11}\nonumber \\
\hphantom{\underline{c_{-12}}=}{}
+ y_{13}\otimes \widehat x_{16} {\wedge} \widehat y_{2} + 2y_{14}\otimes \widehat x_{3}
{\wedge} \widehat x_{13} + y_{14}\otimes \widehat x_{6} {\wedge} \widehat x_{11} +
y_{14}\otimes \widehat x_{8} {\wedge}\widehat x_{9} + y_{14}\otimes \widehat x_{16}
{\wedge} \widehat y_{4}\nonumber \\
\hphantom{\underline{c_{-12}}=}{}
+ 2y_{15}\otimes \widehat x_{6} {\wedge} \widehat x_{8} + y_{16}\otimes \widehat x_{2}
{\wedge} \widehat x_{8} + 2y_{16}\otimes \widehat x_{4}{\wedge} \widehat x_{6} +
y_{16}\otimes \widehat x_{13} {\wedge} \widehat y_{2} + 2y_{16}\otimes \widehat x_{14}
{\wedge} \widehat y_{4},\nonumber\\
 \underline{c_{-6}}
= x_2\otimes \widehat x_7\wedge \widehat x_{10} +x_3\otimes \widehat x_1 \wedge \widehat x_{14}
+x_3\otimes \widehat x_4 \wedge \widehat x_{12} +x_3\otimes \widehat x_8 \wedge \widehat x_{10}
 +2\,x_6\otimes \widehat x_1 \wedge \widehat x_{15}
\nonumber\\
\hphantom{\underline{c_{-6}}=}{}+ 2x_6\otimes \widehat x_7 \wedge \widehat x_{12} +x_8\otimes \widehat x_{10} \wedge \widehat
x_{12} +x_9\otimes \widehat x_4 \wedge \widehat x_{15} +2x_9\otimes \widehat x_7 \wedge
\widehat x_{14}
\nonumber\\
\hphantom{\underline{c_{-6}}=}{}
+2x_9\otimes \widehat x_{10} \wedge \widehat x_{13} + x_{13}\otimes \widehat x_{10} \wedge \widehat x_{15} +2y_1\otimes \widehat x_2 \wedge \widehat
x_{10} +2y_1\otimes \widehat x_4 \wedge \widehat x_7 \nonumber\\
\hphantom{\underline{c_{-6}}=}{}
+y_1\otimes \widehat x_{14}
\wedge \widehat y_3 +y_1\otimes \widehat x_{15} \wedge \widehat y_6+ y_2\otimes \widehat x_1 \wedge \widehat x_{10} +2y_4\otimes \widehat x_1 \wedge \widehat x_7
+y_4\otimes \widehat x_{12}\wedge \widehat y_3 \nonumber\\
\hphantom{\underline{c_{-6}}=}{}
 +y_4\otimes \widehat x_{15} \wedge \widehat
y_9+y_7\otimes \widehat x_1 \wedge \widehat x_4+ y_7\otimes \widehat x_{10} \wedge \widehat y_2 +y_7\otimes \widehat x_{12}\wedge \widehat y_6
 +2y_7\otimes \widehat x_{14} \wedge \widehat y_9 \nonumber\\
\hphantom{\underline{c_{-6}}=}{}
+2y_8\otimes \widehat x_{10} \wedge \widehat y_3
 +y_{10}\otimes \widehat x_1 \wedge \widehat x_2 + 2y_{10}\otimes \widehat x_7 \wedge \widehat y_2 +2y_{10}\otimes \widehat x_8 \wedge \widehat
y_3 +y_{10}\otimes \widehat x_{12} \wedge \widehat y_8 \nonumber\\
\hphantom{\underline{c_{-6}}=}{}
+2y_{10}\otimes \widehat
x_{13} \wedge \widehat y_9 +y_{10}\otimes \widehat x_{15} \wedge \widehat y_{13}+ 2y_{12}\otimes \widehat x_4 \wedge \widehat y_3 +2y_{12}\otimes \widehat x_7 \wedge \widehat
y_6 \nonumber\\
\hphantom{\underline{c_{-6}}=}{}
+2y_{12}\otimes \widehat x_{10} \wedge \widehat y_8 +2y_{13}\otimes \widehat
x_{10} \wedge \widehat y_9+2y_{14}\otimes \widehat x_1 \wedge \widehat y_3 + y_{14}\otimes \widehat x_7 \wedge \widehat y_9 \nonumber\\
\hphantom{\underline{c_{-6}}=}{}
+2y_{15}\otimes \widehat x_1 \wedge \widehat
y_6 +2y_{15}\otimes \widehat x_4 \wedge \widehat y_9 +2y_{15}\otimes \widehat
x_{10}\wedge
\widehat y_{13} ,\nonumber\\
\underline{c_{-3}}= 2x_1\otimes \widehat x_2\wedge \widehat x_7 +2x_1\otimes
\widehat x_4\wedge \widehat x_5 +x_1\otimes \widehat x_{13}\wedge \widehat y_3 +x_1\otimes
\widehat x_{15}\wedge \widehat y_8 +2x_3\otimes \widehat x_{13}\wedge \widehat y_1\nonumber\\
\hphantom{\underline{c_{-3}}=}{}
+ 2x_4\otimes \widehat x_5\wedge \widehat x_7 +2x_8\otimes \widehat x_2\wedge \widehat x_{13}
+2x_8\otimes \widehat x_5\wedge \widehat x_{11} +x_8\otimes \widehat x_7\wedge \widehat x_9
+x_8\otimes \widehat x_{15}\wedge \widehat y_1\nonumber\\
\hphantom{\underline{c_{-3}}=}{}
+ x_{11}\otimes \widehat x_5\wedge \widehat x_{13} +2x_{12}\otimes \widehat x_7\wedge \widehat x_{13}
+2y_2\otimes \widehat x_7\wedge \widehat y_1 +y_2\otimes \widehat x_{13}\wedge
\widehat y_8 +y_4\otimes \widehat x_5\wedge \widehat y_1\nonumber\\
\hphantom{\underline{c_{-3}}=}{}
+ y_5\otimes \widehat x_4\wedge \widehat y_1 +y_5\otimes \widehat x_7\wedge \widehat y_4 +2y_5\otimes
\widehat x_{11}\wedge \widehat y_8 +y_5\otimes \widehat x_{13}\wedge \widehat y_{11}
+y_7\otimes \widehat x_2\wedge \widehat y_1\\ %\label{eq1,3}\\
\hphantom{\underline{c_{-3}}=}{}
+ 2y_7\otimes \widehat x_5\wedge \widehat y_4 +y_7\otimes \widehat x_9\wedge \widehat y_8
+2y_7\otimes \widehat x_{13}\wedge \widehat y_{12} +2y_9\otimes \widehat x_7\wedge \widehat y_8
+y_{11}\otimes \widehat x_5\wedge \widehat y_8\nonumber\\
\hphantom{\underline{c_{-3}}=}{}
+ y_{13}\otimes \widehat x_2\wedge \widehat y_8 +y_{13}\otimes \widehat x_5\wedge \widehat y_{11}
+2y_{13}\otimes \widehat x_7\wedge \widehat y_{12} +y_{13}\otimes \widehat y_1\wedge \widehat y_3
+2y_{15}\otimes \widehat y_1\wedge \widehat y_8.\nonumber
\end{gather*}
\end{Lemma}

\begin{Lemma} \label{g(2,3)} For $\fg=\fg(2,3)$, we consider the Cartan matrix
\[
\begin{pmatrix}
 \hphantom{-}2&-1&-1 \\
 -1&\hphantom{-}2&-1 \\
 -1&-1&\hphantom{-}0
\end{pmatrix}
\]
and the basis even $|$ odd
\begin{gather*}
 x_{1},\ x_{2},\mid x_{3}, \\
 x_4= [x_{1}, x_{2} ],\mid x_5=
 [x_{1}, x_{3} ]\\
\mid x_6= [x_{2}, x_{3} ],\ x_7=
 [x_{3}, [x_{1}, x_{2} ]
 ],\nonumber\\
\mid x_8= [ [x_{1}, x_{2} ]
, [x_{1}, x_{3} ] ],\ x_9=
 [ [x_{1}, x_{2} ], [x_{2}, x_{3} ] ]\\
  \mid x_{10}=
 [ [x_{1}, x_{2} ], [x_{3}, [x_{1}, x_{2} ]
 ] ], \nonumber\\
 x_{11}= [ [x_{2}, x_{3} ], [ [x_{1}, x_{2} ],
 [x_{1}, x_{3} ] ] ]\mid %\label{2,3}
\end{gather*}
Then, $H^2(\fg; \fg)$ is spanned by the cocycles:
\begin{gather*}
\underline{c_{-9}}= 2x_{2}\otimes \widehat x_{9}
{\wedge} \widehat x_{11} + 2x_{4}\otimes \widehat x_{10} {\wedge} \widehat x_{11} +
y_{2}\otimes \widehat x_{5}{\wedge} \widehat x_{11}+ y_{3}\otimes \widehat x_{4}
{\wedge} \widehat x_{11} + 2y_{3}\otimes \widehat x_{7}{\wedge} \widehat x_{10} \nonumber\\
\hphantom{\underline{c_{-9}}=}{}
 +y_{3}\otimes \widehat x_{8}{\wedge} \widehat x_{9}+ y_{4}\otimes \widehat x_{3} {\wedge}
\widehat x_{11} + y_{5}\otimes \widehat x_{2}{\wedge} \widehat x_{11} + 2y_{5}\otimes
\widehat x_{6}{\wedge} \widehat x_{10}+ y_{5}\otimes \widehat x_{7}
{\wedge} \widehat x_{9} \nonumber\\
\hphantom{\underline{c_{-9}}=}{}
+ 2y_{6}\otimes \widehat x_{5}{\wedge} \widehat x_{10} + y_{7}\otimes \widehat x_{3}
{\wedge} \widehat x_{10}+ y_{7}\otimes \widehat x_{5}{\wedge} \widehat x_{9} + 2y_{8}\otimes
\widehat x_{3}{\wedge} \widehat x_{9} + 2y_{9}\otimes \widehat x_{3}
{\wedge}\widehat x_{8}\nonumber\\
 \hphantom{\underline{c_{-9}}=}{}
 + y_{9}\otimes \widehat x_{5}{\wedge} \widehat x_{7} + y_{9}\otimes \widehat x_{11} {\wedge}
\widehat y_{2} + 2y_{10}\otimes \widehat x_{3} {\wedge} \widehat x_{7} + y_{10}\otimes
\widehat x_{5}{\wedge} \widehat x_{6} + y_{10}\otimes \widehat x_{11}
{\wedge} \widehat y_{4} \nonumber\\
\hphantom{\underline{c_{-9}}=}{}
+ 2y_{11}\otimes \widehat x_{2} {\wedge} \widehat x_{5}+ y_{11}\otimes \widehat x_{3}
{\wedge} \widehat x_{4} + y_{11}\otimes \widehat x_{9}{\wedge} \widehat y_{2} +
2y_{11}\otimes \widehat x_{10} {\wedge} \widehat y_{4}.%\label{eqg2,3}
\end{gather*}
\end{Lemma}

\begin{Lemma}\label{Lg4,3} For $\fg=\fg(4,3)$, we select the Cartan matrix
\[
\begin{pmatrix}
 \hphantom{-}2 & -1 & \hphantom{-}0 & \hphantom{-}0 \\
 -1 & \hphantom{-}2 & -2 & -1 \\
 \hphantom{-}0 & -1 & \hphantom{-}2 & \hphantom{-}0 \\
 \hphantom{-}0 & \hphantom{-}1 & \hphantom{-}0 & \hphantom{-}0
\end{pmatrix}
\]
and the even elements of a basis
\begin{gather*}
x_1, \
x_2, \
x_3, \
x_5 = [x_1, x_2], \ x_6 = [x_2, x_3], \
x_8 = [x_2, [x_2, x_3]], \
x_9 = [x_3, [x_1, x_2]], \\
 x_{13} = [[x_1, x_2], [x_2, x_3]], \ x_{15} = [[x_1, x_2], [x_3, [x_1, x_2]]], \\
x_{23} = [[[x_2,x_4],[x_3,[x_1,x_2]]], [[x_2,x_4],[x_3,[x_1,x_2]]]]
\end{gather*}
and the odd elements of the same basis
\begin{gather*}
x_4, \
x_7 = [x_2, x_4], \ x_{10} = [x_4, [x_1, x_2]], \
x_{11} = [x_4, [x_2, x_3]], \
x_{12} = [x_4,[ x_3, [x_1, x_2]]], \\
x_{14} = [[x_2, x_3], [x_2, x_4]], \
x_{16} = [[x_2, x_4], [x_3, [x_1, x_2]]], \
x_{17} = [[x_2, [x_2, x_3]], [x_4, [x_1, x_2]]], \\
x_{18} =[[x_3,[x_1,x_2]],[x_4,[x_1,x_2]]], \
x_{19} =[[x_4,[x_1,x_2]],
 [[x_1,x_2],[x_2,x_3]]],\\
 x_{20} =[[x_4,[x_2,x_3]],
 [[x_1,x_2],[x_2,x_3]]], \
x_{21} = [[x_4,[x_3,[x_1,x_2]]],
 [[x_1,x_2],[x_2,x_3]]],\\
 x_{22} =[[[x_2,x_3],[x_2,x_4]],
 [[x_1,x_2],[x_3,[x_1,x_2]]]].
\end{gather*}
Then, $H^2(\fg; \fg)$ is spanned by the cocycles:
 \begin{gather*}
\underline{c_{-15}}= 2 x_2\otimes \widehat x_{17}\wedge \widehat x_{23}
 +2 x_5\otimes \widehat x_{19}\wedge \widehat x_{23}
 +x_6\otimes \widehat x_{20}\wedge \widehat x_{23}
 +2 x_9\otimes \widehat x_{21}\wedge \widehat x_{23}
\nonumber\\
 \hphantom{\underline{c_{-15}}=}{}
 +y_2\otimes \widehat x_{12}\wedge \widehat x_{23} +y_4\otimes \widehat x_{17}\wedge \widehat x_{21}
 +2 y_4\otimes \widehat x_{19}\wedge \widehat x_{20}
 +2 y_5\otimes \widehat x_{11}\wedge \widehat x_{23}\nonumber\\
 \hphantom{\underline{c_{-15}}=}{}
 +2 y_6\otimes \widehat x_{10}\wedge \widehat x_{23}
 +y_7\otimes \widehat x_9\wedge \widehat x_{23}
 +y_7\otimes \widehat x_{12}\wedge \widehat x_{22}
 +2 y_7\otimes \widehat x_{16}\wedge \widehat x_{21}\nonumber\\
 \hphantom{\underline{c_{-15}}=}{}
 +2 y_7\otimes \widehat x_{18}\wedge \widehat x_{20}
 +2 y_9\otimes \widehat x_7\wedge \widehat x_{23}
 +y_{10}\otimes \widehat x_6\wedge \widehat x_{23}
 +y_{10}\otimes \widehat x_{11}\wedge \widehat x_{22}\nonumber\\
 \hphantom{\underline{c_{-15}}=}{}
 +2 y_{10}\otimes \widehat x_{14}\wedge
 \widehat x_{21} +2
 y_{10}\otimes \widehat x_{16}\wedge \widehat x_{20}
 +y_{11}\otimes \widehat x_5\wedge \widehat x_{23} +y_{11}\otimes \widehat x_{10}\wedge \widehat x_{22} \nonumber\\
 \hphantom{\underline{c_{-15}}=}{}
 +2 y_{11}\otimes \widehat x_{16}\wedge \widehat x_{19}
 +2 y_{11}\otimes \widehat x_{17}\wedge
 \widehat x_{18} +2 y_{12}\otimes \widehat x_2\wedge
 \widehat x_{23} +y_{12}\otimes \widehat x_7\wedge
 \widehat x_{22} \nonumber\\
 \hphantom{\underline{c_{-15}}=}{}
 +2
 y_{12}\otimes \widehat x_{14}\wedge \widehat x_{19} +2 y_{12}\otimes \widehat x_{16}\wedge \widehat x_{17}
 +y_{14}\otimes \widehat x_{10}\wedge \widehat x_{21}
 +2 y_{14}\otimes \widehat x_{12}\wedge
 \widehat x_{19}\nonumber\\
 \hphantom{\underline{c_{-15}}=}{}
 +2 y_{16}\otimes \widehat x_7\wedge
 \widehat x_{21} +y_{16}\otimes \widehat x_{10}\wedge
 \widehat x_{20}
 +y_{16}\otimes \widehat x_{11}\wedge
 \widehat x_{19} +2
 y_{16}\otimes \widehat x_{12}\wedge \widehat x_{17} \nonumber\\
 \hphantom{\underline{c_{-15}}=}{}
 +2 y_{17}\otimes \widehat x_4\wedge \widehat x_{21} +y_{17}\otimes \widehat x_{11}\wedge \widehat x_{18} +2 y_{17}\otimes \widehat x_{12}\wedge \widehat x_{16}
 +2 y_{17}\otimes \widehat x_{23}\wedge \widehat y_2\nonumber\\
 \hphantom{\underline{c_{-15}}=}{}
 +2 y_{18}\otimes \widehat x_7\wedge \widehat x_{20}
 +y_{18}\otimes \widehat x_{11}\wedge \widehat x_{17}
 +2 y_{19}\otimes \widehat x_4\wedge \widehat x_{20}
 +2 y_{19}\otimes \widehat x_{11}\wedge
 \widehat x_{16} \nonumber\\
 \hphantom{\underline{c_{-15}}=}{}
 +y_{19}\otimes \widehat x_{12}\wedge
 \widehat x_{14} +2
 y_{19}\otimes \widehat x_{23}\wedge \widehat y_5
 +2
 y_{20}\otimes \widehat x_4\wedge \widehat x_{19} +y_{20}\otimes \widehat x_7\wedge \widehat x_{18}
 \nonumber\\
 \hphantom{\underline{c_{-15}}=}{}
 +2
 y_{20}\otimes \widehat x_{10}\wedge \widehat x_{16} +y_{20}\otimes \widehat x_{23}\wedge \widehat y_6 +2
 y_{21}\otimes \widehat x_4\wedge \widehat x_{17} +2
 y_{21}\otimes \widehat x_7\wedge \widehat x_{16} \nonumber\\
\hphantom{\underline{c_{-15}}=}{}
 +y_{21}\otimes \widehat x_{10}\wedge \widehat x_{14} +2 y_{21}\otimes \widehat x_{23}\wedge \widehat y_9 +2 y_{22}\otimes \widehat x_7\wedge \widehat x_{12} +y_{22}\otimes \widehat x_{10}\wedge \widehat x_{11}\nonumber\\
\hphantom{\underline{c_{-15}}=}{}
 +2 y_{23}\otimes \widehat x_2\wedge \widehat x_{12} +2 y_{23}\otimes \widehat x_5\wedge \widehat x_{11}
 +2 y_{23}\otimes \widehat x_6\wedge \widehat x_{10} +y_{23}\otimes \widehat x_7\wedge \widehat x_9 \nonumber\\
\hphantom{\underline{c_{-15}}=}{}
 +y_{23}\otimes \widehat x_{17}\wedge \widehat y_2 +2
 y_{23}\otimes \widehat x_{19}\wedge \widehat y_5 +y_{23}\otimes \widehat x_{20}\wedge \widehat y_6 +y_{23}\otimes \widehat x_{21}\wedge \widehat y_9.%\label{g3,3w=0,1}
\end{gather*}
\end{Lemma}

\begin{Lemma}\label{Lg2,6} The Lie superalgebras $\fel(5; 3)$, $\fg(3,3)$, $\fg(8,3)$, $\fg(2,6)$, $\fg(3,6)$, $\fg(4,6)$,
$\fg(6,6)$, $\fg(8,6)$ are rigid.
\end{Lemma}

We give the proof only for $\fg(2,6)$; the other cases were taken care of directly by means of \textsc{SuperLie}, see~\cite{Gr}.

\begin{proof} We have $\fg_\ev=\mathfrak{gl}(6)$ and $\fg_\od=R(\pi_3)$. The space $H^2(\fg;\fg)$ is ``bounded from above'' by
the direct sum of $H^0\big(\fg_\ev; \big(E^2(\fg_\od^*)\otimes
\fg\big)^{\fg_\ev}\big)$ and $H^1\big(\fg_\ev; (\fg_\od^*\otimes
\fg)^{\fg_\ev}\big)$, more exactly,
\begin{gather}
H^0\big(\fg_\ev; \big(E^2(\fg_\od^*)\otimes \fg\big)^{\fg_\ev}\big)=
\big(E^2(\fg_\od^*)\otimes \fg_\ev\big)^{\fg_\ev}\oplus
\big(E^2(\fg_\od^*)\otimes
\fg_\od\big)^{\fg_\ev},\nonumber\\
H^1\big(\fg_\ev; (\fg_\od^*\otimes \fg)^{\fg_\ev}\big)=H^1\big(\fg_\ev;
(\fg_\od^*\otimes \fg_\ev)^{\fg_\ev}\big)\oplus H^1\big(\fg_\ev;
(\fg_\od^*\otimes \fg_\od)^{\fg_\ev}\big).\label{3terms}
\end{gather} Now, observe that the weight of any element of $\fg_\od^*$ has
three nonzero coordinates equal to~$-1$, whereas the weight of any
element of $\fg_\ev$ has nonzero coordinates equal to~1 and~$-1$. No
sum of two weights of $\fg_\od^*$-type with one weight of
$\fg_\od$-type or $\fg_\ev$-type has all coordinates divisible by~3.
This takes care of the 3 summands in \eqref{3terms}, except for $H^1\big(\fg_\ev; (\fg_\od^*\otimes \fg_\od)^{\fg_\ev}\big)$.

In $H^1\big(\fg_\ev; (\fg_\od^*\otimes \fg_\od)^{\fg_\ev}\big)$, there
definitely is an invariant; indeed, there is exactly one:
\begin{gather*}%\label{inv}
\sum f_i \otimes q_i,\qquad \text{where the $f_i$ form a~basis of $\fg_\od$ and the $q_i$
form the dual basis of $\fg_\od^*$}.
\end{gather*}
Further investigation shows
that the corresponding 2-cocycle is a~coboundary.

In these cases, the space $\big(E^2(\fg_\od^*)\otimes \fg_\ev\big)^{\fg_\ev}$ is spanned by~2 elements of weight~0, one of which
is not closed, the other one is a~coboundary.
\end{proof}

\subsection[Deformations of $\mathfrak{osp}(4|2)$ for $p=0$ and $p\geq 3$]{Deformations of $\boldsymbol{\mathfrak{osp}(4|2)}$ for $\boldsymbol{p=0}$ and $\boldsymbol{p\geq 3}$}\label{Lrigid}

The ``classical'' simple Lie superalgebras
of rank $=2$ and $3$ and with indecomposable Cartan matrix are rigid if $p\geq 3$,
except for $\mathfrak{osp}(4|2)$.  For a basis, see formula~\eqref{wk3a}.

\begin{Lemma}\label{osp4_2} The corresponding cocycles look different in the $3$ cases: $p=0$, $p=3$, and $p\geq 5$. Their expressions are, respectively:

For $p=0$:
\begin{gather*} c_0=h_1\otimes \widehat x_2\wedge \widehat y_2
 + h_1\otimes \widehat x_4\wedge
 \widehat y_4-
 h_1\otimes \widehat x_5\wedge \widehat y_5 -h_1\otimes \widehat x_6\wedge \widehat y_6
 - h_3\otimes \widehat x_2\wedge
 \widehat y_2\nonumber\\
\hphantom{c_0=}{}
 +h_3\otimes \widehat x_4\wedge \widehat y_4- h_3\otimes \widehat x_5\wedge
 \widehat y_5 +
 h_3\otimes \widehat x_6\wedge \widehat y_6-2 x_1\otimes \widehat x_4\wedge
 \widehat y_2 -2x_1\otimes \widehat x_6\wedge
 \widehat y_5\nonumber\\
\hphantom{c_0=}{}-2x_3\otimes \widehat x_5\wedge
 \widehat y_2 -2x_3\otimes \widehat x_6\wedge
 \widehat y_4+2y_1\otimes \widehat x_2\wedge
 \widehat y_4+2y_1\otimes \widehat x_5\wedge
 \widehat y_6+2y_3\otimes \widehat x_2\wedge
 \widehat y_5
\nonumber\\
\hphantom{c_0=}{}
 +2 y_3\otimes \widehat x_4\wedge
 \widehat y_6. %\label{eqosp0}
\end{gather*}

For $p=3$:
\begin{gather*} c_0= h_1\otimes \widehat x_2\wedge \widehat y_2+h_1\otimes
\widehat x_3\wedge \widehat y_3+h_1\otimes \widehat x_4\wedge \widehat y_4+h_1\otimes \widehat x_5\wedge
\widehat y_5+ h_1\otimes \widehat x_7\wedge \widehat y_7\nonumber\\
\hphantom{c_0=}{}
+2 h_2\otimes \widehat x_3\wedge
\widehat y_3+ h_2\otimes \widehat x_7\wedge \widehat y_7+x_1\otimes \widehat h_3\wedge \widehat x_1+ x_1\otimes
\widehat x_4\wedge \widehat y_2+x_1\otimes \widehat x_6\wedge \widehat y_5\nonumber\\
\hphantom{c_0=}{}
+x_2\otimes \widehat h_2\wedge
\widehat x_2+2 x_2\otimes \widehat h_3\wedge \widehat x_2
+ 2 x_2\otimes \widehat x_5\wedge \widehat y_3+x_2\otimes \widehat x_7\wedge \widehat y_6+2
x_3\otimes \widehat h_2\wedge \widehat x_3\nonumber\\
\hphantom{c_0=}{}
+x_4\otimes \widehat h_2\wedge \widehat x_4+x_4\otimes
\widehat x_6\wedge \widehat y_3+2 x_4\otimes \widehat x_7\wedge \widehat y_5
+ 2 x_5\otimes \widehat h_3\wedge \widehat x_5+2 x_5\otimes \widehat x_2\wedge
\widehat x_3\nonumber\\
\hphantom{c_0=}{}
+x_6\otimes \widehat x_3\wedge \widehat x_4+x_7\otimes \widehat h_2\wedge \widehat x_7+2
x_7\otimes
\widehat h_3\wedge \widehat x_7+2 x_7\otimes \widehat x_2\wedge \widehat x_6+ 2 x_7\otimes \widehat x_4\wedge \widehat x_5\nonumber\\
\hphantom{c_0=}{}
+2 y_1\otimes \widehat h_3\wedge \widehat y_1+2
y_1\otimes \widehat x_2\wedge \widehat y_4+2 y_1\otimes \widehat x_5\wedge \widehat y_6+2 y_2\otimes
\widehat h_2\wedge \widehat y_2+y_2\otimes \widehat h_3\wedge \widehat y_2\nonumber\\
\hphantom{c_0=}{}
+ y_2\otimes \widehat x_3\wedge \widehat y_5+y_3\otimes \widehat h_2\wedge \widehat y_3+y_3\otimes
\widehat x_2\wedge \widehat y_5+y_3\otimes \widehat x_4\wedge \widehat y_6+2 y_4\otimes
\widehat h_2\wedge \widehat y_4\nonumber\\
\hphantom{c_0=}{}
+2 y_4\otimes \widehat x_3\wedge \widehat y_6+ y_5\otimes \widehat h_3\wedge \widehat y_5+2 y_5\otimes \widehat x_4\wedge \widehat y_7+y_6\otimes
\widehat x_2\wedge \widehat y_7+2 y_7\otimes \widehat h_2\wedge \widehat y_7
\nonumber\\
\hphantom{c_0=}{}
+y_7\otimes \widehat h_3\wedge
\widehat y_7.%\label{eqosp3}
\end{gather*}

For $p\geq 5$ the pattern is as follows $($verified for $p=5, 7, 11)$:
\begin{gather*} c_0=
h_1\otimes \widehat x_
2\wedge
\widehat y_2+
h_1\otimes \widehat x_4\wedge
\widehat y_4+
(p-1)
h_1\otimes \widehat x_5\wedge
\widehat y_5
+ (p-1)
h_1\otimes \widehat x_6\wedge
\widehat y_6\nonumber\\
\hphantom{c_0=}{}
 + h_2\otimes \widehat x_5\wedge
\widehat y_5+ (p-1)
h_2\otimes \widehat x_6\wedge
\widehat y_6
+(p-2)
h_2\otimes \widehat x_7\wedge
\widehat y_7
+ (p-2)
x_1\otimes \widehat x_4\wedge
\widehat y_2
\nonumber\\
\hphantom{c_0=}{}
 + (p-2)
x_1\otimes \widehat x_6\wedge
\widehat y_5 +
x_2\otimes \widehat h_2\wedge
\widehat x_2
 + x_2\otimes \widehat x_5\wedge
\widehat y_3
+ (p-1)
x_2\otimes \widehat x_7\wedge
\widehat y_6
\nonumber\\
\hphantom{c_0=}{}
+ (p-1)
x_3\otimes \widehat h_2\wedge
\widehat x_3
+ x_4\otimes \widehat h_2\wedge
\widehat x_4
+(p-1)
x_4\otimes \widehat x_6\wedge
\widehat y_3+ x_4\otimes \widehat x_7\wedge
\widehat y_5
\nonumber\\
\hphantom{c_0=}{}
+ (p-1)
x_5\otimes \widehat x_2\wedge
\widehat x_3
+
x_6\otimes \widehat x_3\wedge
\widehat x_4
+ x_7\otimes \widehat h_2\wedge
\widehat x_7
+ (p-1)
x_7\otimes \widehat x_2\wedge
\widehat x_6
\nonumber\\
\hphantom{c_0=}{}
+ (p-1)
x_7\otimes \widehat x_4\wedge
\widehat x_5
+ 2
y_1\otimes \widehat x_2\wedge
\widehat y_4
+ 2
y_1\otimes \widehat x_5\wedge
\widehat y_6
+ (p-1)
y_2\otimes \widehat h_2\wedge
\widehat y_2%\label{eqosp}
\\
\hphantom{c_0=}{}
+ 2
y_2\otimes \widehat x_6\wedge
\widehat y_7
+ y_3\otimes \widehat h_2\wedge
\widehat y_3
+
y_3\otimes \widehat x_2\wedge
\widehat y_5
+
y_3\otimes \widehat x_4\wedge
\widehat y_6
+ (p-1)
y_4\otimes \widehat h_2\wedge
\widehat y_4
\nonumber\\
\hphantom{c_0=}{}
+ (p-2)
y_4\otimes \widehat x_5\wedge
\widehat y_7 + (p-1)
y_5\otimes \widehat x_4\wedge
\widehat y_7
+ y_6\otimes \widehat x_2\wedge
\widehat y_7+ (p-1)
y_7\otimes \widehat h_2\wedge
\widehat y_7.\nonumber
\end{gather*}
\end{Lemma}

\section[Results for $p=2$, except for superizations of $ADE$ root types]{Results for $\boldsymbol{p=2}$, except for superizations of $\boldsymbol{ADE}$ root types}\label{SResP=2}

Proofs are obtained by means of
\textsc{SuperLie} package~\cite{Gr}.

Being interested in the classification of simple Lie (super)algebras we
consider the following results as preparatory: we have to
separate semi-trivial cocycles from those that determine true
deforms.  This \textit{open problem} is a separate task to be performed elsewhere.

For example, it is a~well-known folklore (for the proof, see
\cite{SkT1}) that there is only one simple 3-dimensional Lie algebra
for $p=2$, \textit{and hence} the 2 cocycles of Lemma~\ref{L5.1} (they are
integrable thanks to Lemma~\ref{integra3}) are semi-trivial, whereas the very first result of \cite{KrLe} proves that
the two cocycles of Lemma~\ref{L5.2} (of weight $\pm 2$) determine
true deforms (isomorphic to each other by symmetry) sending $\mathfrak{o}\mathfrak{o}_{I\Pi}^{(1)}(1|2)$ to $\mathfrak{o}\mathfrak{o}_{II}^{(1)}(1|2)$.

\begin{Lemma}\label{L5.1} For $\fg=\mathfrak{o}^{(1)}(3)$, the space $H^2(\fg; \fg)$ is
spanned by the cocycles:
\begin{equation*}%\label{eqo3}
c_{-2}=y\otimes \widehat h\wedge \widehat x.
\end{equation*}
\end{Lemma}

\begin{Lemma}\label{L5.2} For $\fg=\mathfrak{o}\mathfrak{o}_{I\Pi}^{(1)}(1|2)$, the space
$H^2(\fg; \fg)$ is spanned by the cocycles:
\begin{equation}\label{eqoo12}
c_{-2}=h_1\otimes (\widehat x_1){}^{\wedge 2} +x_1\otimes \widehat x_1 \wedge \widehat x_2
+y_2\otimes \widehat x_1 \wedge \widehat y_1 +y_2\otimes \widehat x_2 \wedge \widehat y_2.
\end{equation}
\end{Lemma}

\begin{Remark}\label{ssRemar} As shown in \cite{KrLe}, the cocycle \eqref{eqoo12} deforms $\mathfrak{o}\mathfrak{o}_{I\Pi}^{(1)}(1|2)$ to $\mathfrak{o}\mathfrak{o}_{II}^{(1)}(1|2)$.
\end{Remark}

\begin{Lemma}\label{L5.2a} For $\fg=\mathfrak{o}\mathfrak{o}_{II}^{(1)}(1|2)$, the space
$H^2(\fg; \fg)$ is spanned by the cocycles:
\begin{gather*}
c_1= E^{2,2}\otimes \widehat{\big(E^{1,3}+ E^{3,1}\big)}^{\wedge 2}+
 E^{1,2}\otimes \widehat{\big(E^{2,3}+ E^{3,2}\big)}\wedge
\widehat{\big(E^{1,3}+ E^{3,1}\big)}\nonumber\\
\hphantom{c_1=}{} +
 E^{2,1}\otimes \widehat{\big(E^{2,3}+ E^{3,2}\big)}\wedge
\widehat{\big(E^{1,3}+ E^{3,1}\big)},\nonumber\\
 c_2= E^{2,2}\otimes \widehat{\big(E^{1,2}+ E^{2,1}\big)}^{\wedge 2}+ \big(E^{1,2}+E^{2,1}\big)\otimes \widehat E^{3,3}\wedge
 \widehat{\big(E^{1,2}+ E^{2,1}\big)}\nonumber\\
\hphantom{c_2=}{} +
 \big(E^{1,3}+E^{3,1}\big)\otimes \widehat E^{3,3}\wedge
 \widehat{\big(E^{1,3}+ E^{3,1}\big)}\nonumber\\
\hphantom{c_2=}{}+
 \big(E^{1,3}+E^{3,1}\big)\otimes \widehat{\big(E^{2,3}+ E^{3,2}\big)}\wedge
 \widehat{\big(E^{1,2}+ E^{2,1}\big)}.%\label{eqooII12}
\end{gather*}
\end{Lemma}

\begin{Lemma}\label{L5.3} For $\fg=\mathfrak{sl}(3)$, we have $H^2(\fg;
\fg)=0$ for any $p$ $($even for $p=3)$.
\end{Lemma}

\begin{Lemma}\label{L5.4} For $\fg=\mathfrak{o}^{(1)}_I(5)\simeq \mathfrak{o}^{(1)}_\Pi(5)$, we take the Cartan matrix
\[
\begin{pmatrix}
\overline{1}&1\\
1&\ev\\
\end{pmatrix}
\]
and the basis
\begin{equation*} %\label{o5}
x_1,\ x_2,\ x_3=[x_1,x_2],\ x_4=[x_1,[x_1,x_2]].
\end{equation*}
Then, $H^2(\fg; \fg)$ is spanned by the cocycles:
\begin{gather*}
c_{-4}= h_1\otimes \widehat x_2\wedge \widehat x_4+x_1\otimes \widehat x_3\wedge
\widehat x_4+y_1\otimes \widehat x_2\wedge
\widehat x_3\nonumber\\
\hphantom{c_{-4}=}{}
+y_2\otimes \widehat h_2\wedge \widehat x_4+y_3\otimes \widehat h_2\wedge
\widehat x_3+y_4\otimes \widehat h_2\wedge \widehat x_2,\nonumber\\
c_{-2}=h_1\otimes \widehat x_4\wedge \widehat y_2 +x_2\otimes \widehat h_1\wedge \widehat x_4
+x_2\otimes \widehat h_2\wedge \widehat x_4 +x_3\otimes \widehat x_1\wedge \widehat x_4 +y_1\otimes
\widehat h_1\wedge \widehat x_1\nonumber\\
\hphantom{c_{-2}=}{}
+y_1\otimes \widehat h_2\wedge \widehat x_1 +y_1\otimes \widehat x_3\wedge \widehat y_2 +y_4\otimes
\widehat h_1\wedge \widehat y_2 +y_4\otimes \widehat h_2\wedge \widehat y_2 +y_4\otimes \widehat x_1\wedge
\widehat y_3.%\label{eqo5}
\end{gather*}
\end{Lemma}

\begin{Conjecture}[proved for $n=1,2,3,4$]\label{ssConj} For $\fg=\mathfrak{o}^{(1)}_I(2n+1)\simeq\mathfrak{o}^{(1)}_\Pi(2n+1)$, the space
$H^2(\fg; \fg)$ is similarly spanned by the cocycles of weight
$\pm 2, \pm 4,\ldots,\pm2n$.
\end{Conjecture}

\begin{Lemma}\label{L5.4pp} For $\fg=\mathfrak{o}^{(1)}_I(4)$, the space $H^2(\fg; \fg)$ is spanned by the following cocycles whose indices are just numbers of these cocycle, not degree, since the algebra is non-split:
\begin{gather*}
c_{1}= \big(E^ {1,3}+E^{3,1}\big) \otimes \big(\widehat{E^ {1,2}+E^ {2,1}}\big)\wedge
 \big(\widehat{E^ {2,3}+E^ {3,2}}\big) \nonumber\\
 \hphantom{c_{1}=}{}+
 \big(E^ {1,4}+E^ {4,1}\big) \otimes \big(\widehat{E^ {1,2}+E^ {2,1}}\big)\wedge
 \big(\widehat{E^ {2,4}+E^ {4,2}}\big) \nonumber\\
\hphantom{c_{1}=}{}
 +\big( E^ {3,4}+ E^ {4,3}\big)\otimes \big(\widehat{E^ {2,3}+E^ {3,2}}\big)\wedge
 \big(\widehat{E^ {2,4}+E^ {4,2}}\big),\nonumber\\
 c_2= \big( E^ {1,3} +E^ {3,1}\big) \otimes \big(\widehat{E^ {2,3}+E^ {3,2}}\big)\wedge
 \big(\widehat{E^ {3,4}+E^ {4,3}}\big)\nonumber\\
 \hphantom{c_2=}{}+
 (E^ {2,3}+E^ {3,2}) \otimes \big(\widehat{E^ {1,3}+E^ {3,1}}\big)\wedge
 \big(\widehat{E^ {3,4}+E^ {4,3}}\big )\nonumber \\
\hphantom{c_2=}{} +
 \big(E^ {3,4}+E^ {4,3}\big) \otimes \big(\widehat{E^ {1,3}+E^ {3,1}}\big)\wedge
 \big(\widehat{E^ {2,3}+E^ {3,2}}\big),\nonumber \\
 c_3= \big(E^ {1,2}+E^ {2,1}\big)\otimes \big(\widehat{E^ {1,3}+E^ {3,1}}\big)\wedge
 \big(\widehat{E^ {2,3}+E^ {3,2}}\big)\nonumber\\
\hphantom{c_3=}{}
 + \big(E^ {1,4}+E^ {4,1}\big)\otimes \big(\widehat{E^ {1,3}+E^ {3,1}}\big)\wedge
 \big(\widehat{E^ {3,4}+E^ {4,3}}\big)
 \nonumber \\
\hphantom{c_3=}{}
 +\big(E^ {2,4}+E^ {4,2}\big)\otimes \big(\widehat{E^ {2,3}+E^ {3,2}}\big)\wedge
 \big(\widehat{E^ {3,4}+E^ {4,3}}\big),\nonumber \\
 c_4= \big(E^ {1,2}+E^ {2,1}\big)\otimes \big(\widehat{E^ {2,3}+E^ {3,2}}\big)\wedge
 \big(\widehat{E^ {2,4}+E^ {4,2}}\big)\nonumber\\
\hphantom{c_4=}{}
 + \big(E^ {2,3}+E^ {3,2}\big) \otimes \big(\widehat{E^ {1,2}+E^ {2,1}}\big)\wedge
 \big(\widehat{E^ {2,4}+E^ {4,2}}\big)\nonumber \\
\hphantom{c_4=}{}
 +
 \big(E^ {2,4}+E^ {4,2}\big) \otimes \big(\widehat{E^ {1,2}+E^ {2,1}}\big)\wedge
 \big(\widehat{E^ {2,3}+E^ {3,2}}\big),\nonumber \\
 c_5= \big(E^ {1,2}+E^ {2,1}\big)\otimes \big(\widehat{E^ {1,3}+E^ {3,1}}\big)\wedge
 \big(\widehat{E^ {1,4}+E^ {4,1}}\big) \nonumber\\
\hphantom{c_5=}{}
 + \big(E^ {1,3}+E^ {3,1}\big)\otimes \big(\widehat{E^ {1,2}+E^ {2,1}}\big)\wedge
 \big(\widehat{E^ {1,4}+E^ {4,1}}\big)
 \nonumber \\
\hphantom{c_5=}{}
 +\big(E^ {1,4}+E^ {4,1}\big)\otimes \big(\widehat{E^ {1,2}+E^ {2,1}}\big)\wedge
 \big(\widehat{E^ {1,3}+E^ {3,1}}\big),
\nonumber \\
 c_6= \big(E^ {1,2}+E^ {2,1}\big)\otimes \big(\widehat{E^ {1,4}+E^ {4,1}}\big)\wedge
 \big(\widehat{E^ {2,4}+E^ {4,2}}\big)\nonumber\\
\hphantom{c_6=}{}
 +
 \big(E^ {1,3}+E^ {3,1}\big)\otimes \big(\widehat{E^ {1,4}+E^ {4,1}}\big)\wedge
 \big(\widehat{E^ {3,4}+E^ {4,3}}\big)
\nonumber \\
\hphantom{c_6=}{}
 +\big(E^ {2,3}+E^ {3,2}\big)\otimes \big(\widehat{E^ {2,4}+E^ {4,2}}\big)\wedge
 \big(\widehat{E^ {3,4}+E^ {4,3}}\big).%\label{eqo4}
 \end{gather*}
\end{Lemma}

\subsection[A relation between $C^{\bcdot}(\fg)$ and $C^{\bcdot}(\textbf{F}(\fg))$ (written with the help of A.~Krutov)]{A relation between $\boldsymbol{C^{\bcdot}(\fg)}$ and $\boldsymbol{C^{\bcdot}(\textbf{F}(\fg))}$\\ (written with the help of A.~Krutov)}\label{F(C)}

Let $\fg$ be a~Lie superalgebra and~$\textbf{F}(\fg)$ its desuperization, i.e., \textbf{F} is the functor that forgets squaring. If $p=2$, then $E^{\bcdot}(V)\subset S^{\bcdot}(V)$, so there are an injective map $i\colon C^{\bcdot}(\textbf{F}(\fg))\tto C^{\bcdot}(\fg)$. Note that this is an embedding of vector spaces, not of algebras: the fact $ab=0$ for some $a,b\in C^{\bcdot}(\textbf{F}(\fg))$ does not imply $i(a)i(b)=0$.

Here are two pairs of examples: $\fwk(3; \alpha)$ vs.\ $\fbgl(3; \alpha)$ and $\fwk(4; \alpha)$ vs. $\fbgl(4; \alpha)$.

\begin{Lemma}\label{SSwk3a} For $\fg=\fwk(3; \alpha)$, where $\alpha\neq 0, 1$,
we take the Cartan matrix
\[
\begin{pmatrix}
\ev&\alpha&0\\
\alpha&\ev&1\\
0&1&\ev
\end{pmatrix}
\]
and the basis
\begin{gather}
x_1,\ x_2,\
x_3,\
x_4=[x_1,x_2],\ x_5=[x_2,x_3],\
x_6=[x_3,[x_1,x_2]],\
x_7=[[x_1,x_2],[x_2,x_3]].\label{wk3a}
\end{gather}
Then, $H^2(\fg; \fg)$ is spanned by the cocycles:
\begin{gather}
c_{-6}= \alpha(1+\alpha)y_1\otimes \widehat x_3\wedge \widehat x_7 +\alpha
y_1\otimes \widehat x_5\wedge \widehat x_6 + \alpha(1+\alpha) y_3\otimes \widehat x_1\wedge
\widehat x_7 +\alpha ^2y_3\otimes \widehat x_4\wedge
\widehat x_6\nonumber\\
\hphantom{c_{-6}=}{}+ \alpha y_4\otimes \widehat x_3\wedge \widehat x_6 +\alpha y_5\otimes \widehat x_1\wedge
\widehat x_6 +y_6\otimes \widehat x_1\wedge \widehat x_5 +\alpha y_6\otimes
\widehat x_3\wedge \widehat x_4+y_7\otimes \widehat x_1\wedge \widehat x_3,\nonumber\\
c_{-4}^1= \alpha (1+\alpha) x_1\otimes \widehat x_3\wedge \widehat x_7 +\alpha
x_1\otimes \widehat x_5\wedge \widehat x_6 +\alpha y_2\otimes \widehat x_3\wedge \widehat x_5
+\alpha y_3\otimes \widehat x_2\wedge \widehat x_5\nonumber\\
\hphantom{c_{-4}^1=}{}
+ \alpha (1+\alpha) y_3\otimes \widehat x_7\wedge \widehat y_1+\alpha y_5\otimes
\widehat x_2\wedge \widehat x_3 +\alpha y_5\otimes \widehat x_6\wedge \widehat y_1 +y_6\otimes
\widehat x_5\wedge \widehat y_1 \nonumber\\
\hphantom{c_{-4}^1=}{}
+ y_7\otimes \widehat x_3\wedge
\widehat y_1,\nonumber\\
c_{-4}^2= \alpha (1+\alpha) x_3\otimes \widehat x_1\wedge \widehat x_7 +\alpha
^2x_3\otimes \widehat x_4\wedge \widehat x_6 +\alpha y_1\otimes \widehat x_2\wedge \widehat x_4
+\alpha (1+\alpha) y_1\otimes \widehat x_7\wedge \widehat y_3\nonumber\\
\hphantom{c_{-4}^2=}{}
+ \alpha y_2\otimes \widehat x_1\wedge \widehat x_4 +y_4\otimes \widehat x_1\wedge \widehat x_2
+\alpha y_4\otimes \widehat x_6\wedge \widehat y_3 +\alpha y_6\otimes \widehat x_4\wedge
\widehat y_3 +y_7\otimes \widehat x_1\wedge
\widehat y_3,\nonumber\\
c_{-2}= \alpha x_1\otimes \widehat x_2\wedge \widehat x_4 +\alpha (1+\alpha)
x_1\otimes \widehat x_7\wedge \widehat y_3 +\alpha x_3\otimes \widehat x_2\wedge \widehat x_5
+\alpha (1+\alpha) x_3\otimes \widehat x_7\wedge \widehat y_1\nonumber\\
\hphantom{c_{-2}=}{}
+ \alpha y_2\otimes \widehat x_4\wedge \widehat y_1 +\alpha y_2\otimes \widehat x_5\wedge
\widehat y_3 +y_4\otimes \widehat x_2\wedge \widehat y_1
+\alpha y_5\otimes \widehat x_2\wedge \widehat y_3 +y_7\otimes \widehat y_1\wedge \widehat y_3,\nonumber\\
c_{0}= h_1\otimes \widehat x_4 \wedge \widehat y_4 +h_1\otimes \widehat x_6 \wedge \widehat y_6
 +h_1\otimes \widehat x_7 \wedge \widehat y_7 +\alpha h_3\otimes
\widehat x_7 \wedge \widehat y_7
+h_4\otimes \widehat x_2 \wedge \widehat y_2\nonumber\\
\hphantom{c_{0}=}{}
+ \alpha h_4\otimes \widehat x_4 \wedge \widehat y_4 +h_4\otimes \widehat x_5 \wedge \widehat
y_5 +\alpha h_4\otimes \widehat x_6 \wedge \widehat y_6+x_1\otimes \widehat x_4 \wedge
\widehat y_2 +x_1\otimes \widehat x_6 \wedge \widehat y_5\nonumber\\
\hphantom{c_{0}=}{}
+ \alpha x_2\otimes \widehat x_7 \wedge \widehat y_6+x_4\otimes \widehat x_7 \wedge \widehat y_5
+\alpha x_5\otimes \widehat x_7 \wedge \widehat y_4 +x_6\otimes \widehat x_7
\wedge \widehat y_2 +y_1\otimes \widehat x_2 \wedge \widehat y_4\nonumber\\
\hphantom{c_{0}=}{}
+ y_1\otimes \widehat x_5 \wedge \widehat y_6 +\alpha y_2\otimes \widehat x_6 \wedge \widehat
y_7+y_4\otimes \widehat x_5 \wedge \widehat y_7 +\alpha y_5\otimes \widehat x_4 \wedge \widehat
y_7+y_6\otimes \widehat x_2 \wedge \widehat y_7.\label{eqwk3a}
\end{gather}

In the case of $\fwk^{(1)}(3; \alpha)$ all cocycles are the same as above except $c_0$ which takes the form
\begin{gather*}
c_0=
(\alpha +1) h_1\otimes \widehat x_4\wedge \widehat y_4+(\alpha
+1) h_1\otimes \widehat x_6\wedge \widehat y_6+(\alpha
 +1) h_1\otimes \widehat x_7\wedge \widehat y_7+h_2\otimes \widehat x_2\wedge \widehat y_2\\
 \hphantom{c_0=}{}
 +\alpha
 h_2\otimes \widehat x_4\wedge \widehat y_4+h_2\otimes \widehat x_5\wedge \widehat y_5 +
 \alpha
h_2\otimes \widehat x_6\wedge \widehat y_6 +\big(\alpha ^2+\alpha \big)
 h_3\otimes \widehat x_7\wedge \widehat y_7\\
\hphantom{c_0=}{}
 +x_1\otimes \widehat h_2\wedge \widehat x_1 +(\alpha
+1) x_1\otimes \widehat x_4\wedge \widehat y_2 +(\alpha+1) x_1\otimes \widehat x_6\wedge \widehat y_5 +
 \big(\alpha^2+\alpha \big)
 x_2\otimes \widehat x_7\wedge \widehat y_6 \\
\hphantom{c_0=}{}
 +x_3\otimes \widehat h_2\wedge \widehat x_3 +x_4\otimes \widehat h_2\wedge \widehat x_4 +(\alpha
 +1) x_4\otimes \widehat x_7\wedge \widehat y_5 +x_5\otimes \widehat h_2\wedge \widehat x_5\\
\hphantom{c_0=}{}
 +\big(\alpha^2+\alpha \big)
 x_5\otimes \widehat x_7\wedge \widehat y_4
+ (\alpha
 +1) x_6\otimes \widehat x_7\wedge \widehat y_2 +y_1\otimes \widehat h_2\wedge \widehat y_1 +(\alpha
 +1) y_1\otimes \widehat x_2\wedge \widehat y_4 \\
\hphantom{c_0=}{}
 +(\alpha+1) y_1\otimes \widehat x_5\wedge \widehat y_6+\big(\alpha^2+\alpha \big)
 y_2\otimes \widehat x_6\wedge \widehat y_7 +y_3\otimes \widehat h_2\wedge \widehat y_3 + y_4\otimes \widehat h_2\wedge \widehat y_4\\
\hphantom{c_0=}{}
 +(\alpha
 +1) y_4\otimes \widehat x_5\wedge \widehat y_7 +y_5\otimes \widehat h_2\wedge \widehat y_5+\big(\alpha
 ^2+\alpha \big)
 y_5\otimes \widehat x_4\wedge \widehat y_7+(\alpha
 +1) y_6\otimes \widehat x_2\wedge \widehat y_7.
\end{gather*}

In the case of $\fwk^{(1)}(3;\alpha)/\fc$ all the cocycles are the same except $c_0$ which takes the form
\begin{gather*}
c_0= h_2\otimes \widehat x_2\wedge \widehat y_2+\alpha h_2\otimes \widehat x_4\wedge \widehat y_4+h_2\otimes \widehat x_5\wedge \widehat y_5+\alpha h_2\otimes \widehat x_6\wedge \widehat y_6+\big(\alpha ^2+\alpha \big)
 h_3\otimes \widehat x_4\wedge \widehat y_4\\
\hphantom{c_0=}{}
 +\big(\alpha ^2+\alpha \big)
 h_3\otimes \widehat x_6\wedge \widehat y_6 +
 x_1\otimes \widehat h_2\wedge \widehat x_1 +(\alpha
 +1) x_1\otimes \widehat x_4\wedge \widehat y_2 +(\alpha
 +1) x_1\otimes \widehat x_6\wedge \widehat y_5 \\
\hphantom{c_0=}{}
 +\big(\alpha ^2+\alpha \big)
 x_2\otimes \widehat x_7\wedge \widehat y_6+x_3\otimes \widehat h_2\wedge \widehat x_3+x_4\otimes \widehat h_2\wedge \widehat x_4+
 (\alpha
 +1) x_4\otimes \widehat x_7\wedge \widehat y_5\\
\hphantom{c_0=}{}
 +x_5\otimes \widehat h_2\wedge \widehat x_5+\big(\alpha ^2+\alpha \big)
 x_5\otimes \widehat x_7\wedge \widehat y_4+(\alpha
 +1) x_6\otimes \widehat x_7\wedge \widehat y_2+y_1\otimes \widehat h_2\wedge \widehat y_1\\
\hphantom{c_0=}{}
 + (\alpha
 +1) y_1\otimes \widehat x_2\wedge \widehat y_4+(\alpha
 +1) y_1\otimes \widehat x_5\wedge \widehat y_6 +\big(\alpha ^2+\alpha \big)
 y_2\otimes \widehat x_6\wedge \widehat y_7 +y_3\otimes \widehat h_2\wedge \widehat y_3\\
\hphantom{c_0=}{}
 +y_4\otimes \widehat h_2\wedge \widehat y_4 +(\alpha
 +1) y_4\otimes \widehat x_5\wedge \widehat y_7 +
 y_5\otimes \widehat h_2\wedge \widehat y_5 +\big(\alpha ^2+\alpha \big)
 y_5\otimes \widehat x_4\wedge \widehat y_7 \\
\hphantom{c_0=}{}
 +(\alpha
 +1) y_6\otimes \widehat x_2\wedge \widehat y_7.
\end{gather*}
\end{Lemma}

\begin{Lemma}\label{bgl3a} For $\fg=\fbgl(3;\alpha)$, where $\alpha\neq 0, 1$,
$($the super analog of $\fwk(3;\alpha)$ and a~nonexisting for $p=2$ analog of
$\mathfrak{osp}(4|2; \alpha))$, we take the Cartan matrix and the
basis as for ${\fg=\fwk(3;\alpha)}$, see~\eqref{wk3a} with all
Chevalley generators even, except $x_1$, $y_1$. Then, $H^2(\fg; \fg)$
is spanned by the following cocycles~\eqref{eqbgl3a} $($compare with the
cocycles \eqref{eqwk3a} of the same weight and superscript for the
desuperization of $\fbgl(3;\alpha))$. We also give the expressions for the deformed squaring~$s_{{\rm new}}(\cdot)$ defined by the global with parameter $t$ deforms of these infinitesimal cocycles $(h_1+\alpha h_3$ spans the center$)$.

Let $c_i$ be a~cocycle
 of~$\fwk(3;\alpha)$ and $\overline c_i$ the corresponding cocycle of $\fbgl(3;\alpha)$. Let ``\textbf{new!}'' mark the cocycle that does not exist for~$\fwk(3;\alpha)$. We have
\begin{gather}
 \overline c_{-8} = (h_1+\alpha h_3)\otimes \widehat x_7^2,\quad \textbf{new!}\nonumber\\
 \overline c_{-6}^1 = i(c_{-6}) + \alpha^2(h_2 + (1+\alpha)h_4)\otimes \widehat x_6,\nonumber\\
 \overline c_{-6}^2 = (h_1 + \alpha h_3)\otimes \widehat x_6^2,\quad \textbf{new!}\nonumber\\
 \overline c_{-4}^1 = i(c_{-4}^1) + \alpha^2(h_2+h_3+(1+\alpha)h_4)\otimes \widehat x_4^2,\nonumber\\
 \overline c_{-4}^2 = i(c_{-4}^2),\nonumber\\
 \overline c_{-4}^3 = (h_1 + \alpha h_3)\otimes \widehat x_4^2,\quad \textbf{new!}\nonumber\\
 \overline c_{-2}^1 = i(c_{-2}),\nonumber\\
 \overline c_{-2}^2 = (h_1 + \alpha h_3)\otimes x_1^2,\quad \textbf{new!}\nonumber\\
 \overline c_0 = i(c_0).\label{eqbgl3a}
\end{gather}
The new cocycles deform only squaring $($modulo center spanned by $h_1 + \alpha h_3$ they vanish$)$.
\end{Lemma}

\begin{Remarks}\quad
\begin{enumerate}\itemsep=0pt
\item[1)] It would
 be interesting to prove that these new cocycles in \eqref{eqbgl3a} define semi-trivial deforms. Better answer a~more general
 question: let $\fg$ be a~Lie superalgebra, and $\fc$ its center. Let $c = \sum_i c_i \otimes \widehat e_i^{\wedge 2}$,
 where $c_i\in\fc$ and $e_i\in \fg_\od$, be a~nontrivial cocycle. When the corresponding deform is semi-trivial?

\item[2)]Observe that, contrary to a~possible hasty conjecture, the cohomology of Lie superalgebra and its desuperization do not necessary
coincide even if the squaring in Lie superalgebra is zero on each element of the basis. This does not always imply that the squaring
is the zero mapping either since $(x+y)^2 = x^2+y^2 + [x,y]$.
\end{enumerate}

For example, $\fwk^{(1)}(3;\alpha)/\fc$ has more
cocycles than $\fbgl^{(1)}(3;\alpha)/\fc$. Indeed, the superization of the $\fg$-valued 2-form~$c_{-6}$ defined by~\eqref{eqwk3a}
is not even a~cocycle of $\fbgl^{(1)}(3;\alpha)/\fc$:
\begin{gather*}
 d_{\fbgl^{(1)}/\fc}(c_{-6}) =
 \alpha^2 \big(x_1\otimes\widehat x_1\wedge\widehat x_6 ^{\wedge 2} +x_3\otimes\widehat x_3\wedge\widehat x_6 ^{\wedge 2}
 + x_4\otimes\widehat x_4\wedge\widehat x_6 ^{\wedge 2} +x_5\otimes\widehat x_5\wedge\widehat x_6 ^{\wedge 2}\\
\hphantom{d_{\fbgl^{(1)}/\fc}(c_{-6}) =}{} + y_1\otimes\widehat x_6 ^{\wedge 2}\wedge\widehat y_1
 + y_3\otimes\widehat x_6 ^{\wedge 2}\wedge\widehat y_3
 + y_4\otimes\widehat x_6 ^{\wedge 2}\wedge\widehat y_4 + y_5\otimes\widehat x_6 ^{\wedge 2}\wedge \widehat y_5\big)
 \neq 0.
\end{gather*}
However, the deformed bracket $[a,b]_{\fbgl;c_{-6}} = [a,b]_{\fbgl} + t c_{-6}(a,b)$ still satisfies the part of Jacobi
identity not involving squaring.

Since $a^2 = 0$ for all $a\in\big(\fbgl^{(1)}(3;\alpha)/\fc\big)_\od$, the additional axiom
\[
 [x^2,y] = [x,[x,y]] \qquad \text{for any $x\in\fg_\od$, $y\in\fg$,}
\]
reduces to the form
\[
 [x,[x,y]] = 0\qquad \text{for any $x\in\fg_\od$, $y\in\fg$.}
\]
But for $x_1$ and $x_6$ (both odd) we have
\[
 [x_6,[x_6,x_1]_{\fbgl;c_{-6}}]_{\fbgl;c_{-6}} = [x_6, \alpha y_5]_{\fbgl;c_{-6}} = \alpha^2 x_1 \neq 0.
\]
Therefore, the deformed bracket $[\cdot,\cdot]_{\fbgl,c_{-6}}$ does not define a~Lie
superalgebra.

Summary for $\fbgl^{(1)}(3;\alpha)/\fc$ and $\fwk^{(1)}(3;\alpha)/\fc$:
 \begin{gather*}
 \overline c_{-4} = i(c_{-4}^2),\qquad
 \overline c_{-2} = i(c_{-2}),\qquad
 \overline c_0 = i(c_0),\\
 d_{\fbgl^{(1)}/\fc}(c_{-6}) \neq 0,\qquad
 d_{\fbgl^{(1)}/\fc}(c_{-4}^1) \neq 0.
 \end{gather*}
\end{Remarks}

\begin{Lemma}\label{L5.7} For $\fg=\fwk(4;\alpha)$, where $\alpha\neq 0, 1$,
we take the Cartan matrix
\begin{gather*}
\begin{pmatrix}
 \ev &\alpha&1 &0 \\
 \alpha&\ev&0 &0 \\
 1 &0 &\ev &1 \\
 0 &0 &1 &\ev
\end{pmatrix},
\end{gather*}
 and the basis
\begin{gather*}
 x_1,\ x_2,\ x_3,\ x_4,\
x_5=[x_1,x_2],\ x_6=[x_1,x_3], \ x_7=[x_3,x_4], \
x_8=[x_3,[x_1,x_2]],\nonumber\\
x_9=[x_4,[x_1,x_3]], \
x_{10}=[[x_1,x_2],[x_1,x_3]],\
x_{11}=[[x_1,x_2],[x_3,x_4]],\nonumber\\
x_{12}=[[x_1,x_2],[x_4,[x_1,x_3]]],\
x_{13}=[[x_3,[x_1,x_2]],[x_4,[x_1,x_3]]],\nonumber\\
x_{14}=[[x_4,[x_1,x_3]],[[x_1,x_2],[x_1,x_3]]],\
x_{15}=[[[x_1,x_2],[x_1,x_3]],[[x_1,x_2],[x_3,x_4]]].%\label{wk4cm}
\end{gather*}
Then, $H^2(\fg; \fg)$ is spanned by the cocycles which for
polygraphical reasons are divided into several groups $($see \eqref{eqpsl4c3}, \eqref{eqwk4ab}$)$:
\begin{gather}
c_{-12}= (1+\alpha)\big(\alpha^3y_3\otimes \widehat x_9\wedge \widehat x_{15}
+\alpha^3y_3\otimes \widehat x_{11}\wedge \widehat x_{14} +\alpha ^2y_3\otimes
\widehat x_{12}\wedge \widehat x_{13} +\alpha ^3y_6\otimes \widehat x_7\wedge \widehat x_{15}
\nonumber\\
\hphantom{c_{-12}=}{}
+ \alpha ^2y_6\otimes \widehat x_{11}\wedge \widehat x_{13} +\alpha ^3 y_7\otimes
\widehat x_6\wedge \widehat x_{15} +\alpha ^3y_7\otimes \widehat x_8\wedge \widehat x_{14} +\alpha
^2y_8\otimes \widehat x_7\wedge \widehat x_{14}\nonumber\\
\hphantom{c_{-12}=}{}
+ \alpha y_8\otimes \widehat x_9\wedge \widehat x_{13} +\alpha ^3y_9\otimes
\widehat x_3\wedge \widehat x_{15} +\alpha ^2 y_9\otimes \widehat x_8\wedge \widehat x_{13} +\alpha
y_{10}\otimes \widehat x_7\wedge \widehat x_{13} \nonumber\\
\hphantom{c_{-12}=}{}
+\alpha ^2y_{11}\otimes \widehat x_3\wedge
\widehat x_{14}
+ \alpha y_{11}\otimes \widehat x_6\wedge \widehat x_{13} +\alpha y_{12}\otimes
\widehat x_3\wedge \widehat x_{13} +\alpha y_{13}\otimes \widehat x_3\wedge \widehat x_{12} \nonumber\\
\hphantom{c_{-12}=}{}
+\alpha y_{13}\otimes \widehat x_7\wedge \widehat x_{10}\big) +
 \alpha y_{13}\otimes \widehat x_8\wedge \widehat x_9 +\alpha y_{14}\otimes
\widehat x_3\wedge \widehat x_{11} +\alpha y_{14}\otimes \widehat x_7\wedge \widehat x_8\nonumber\\
\hphantom{c_{-12}=}{}
 +
y_{15}\otimes \widehat x_3\wedge
\widehat x_9 +y_{15}\otimes \widehat x_6\wedge \widehat x_7+
 \alpha
^2\big(1+\alpha^2\big) y_7\otimes \widehat x_{10}\wedge \widehat x_{13} +\alpha
 y_{13}\otimes \widehat x_6\wedge \widehat x_{11},\nonumber\\
c_{-10}= (1+\alpha)\big(\alpha ^3x_3\otimes \widehat x_9\wedge \widehat x_{15}
+\alpha ^3x_3\otimes \widehat x_{11}\wedge \widehat x_{14} + \alpha ^2x_3\otimes
\widehat x_{12}\wedge \widehat x_{13} +\alpha ^3y_1\otimes \widehat x_4\wedge \widehat x_{15}
\nonumber\\
\hphantom{c_{-10}=}{}
+ \alpha ^2y_1\otimes \widehat x_{11}\wedge \widehat x_{12} +\alpha ^3y_4\otimes
\widehat x_1\wedge \widehat x_{15} +\alpha ^3y_4\otimes \widehat x_5\wedge \widehat x_{14} +\alpha
^2 y_4\otimes \widehat x_{10}\wedge \widehat x_{12} \nonumber\\
\hphantom{c_{-10}=}{}
+\alpha ^2
y_5\otimes \widehat x_4\wedge \widehat x_{14} +
 \alpha y_5\otimes \widehat x_9\wedge \widehat x_{12} +\alpha ^2y_9\otimes
\widehat x_5\wedge \widehat x_{12} +\alpha ^3 y_9\otimes \widehat x_{15}\wedge \widehat y_3\nonumber\\
\hphantom{c_{-10}=}{}
+\alpha y_{10}\otimes \widehat x_4\wedge \widehat x_{12} +\alpha
y_{11}\otimes \widehat x_1\wedge \widehat x_{12} +
 \alpha ^2 y_{11}\otimes \widehat x_{14}\wedge \widehat y_3 \nonumber\\
\hphantom{c_{-10}=}{}
 +\alpha y_{12}\otimes
\widehat x_4\wedge \widehat x_{10} +\alpha y_{12}\otimes \widehat x_{13}\wedge \widehat y_3 +\alpha
y_{13}\otimes \widehat x_{12}\wedge \widehat y_3 \big)
 +\alpha y_{12}\otimes \widehat x_1\wedge \widehat x_{11} \nonumber\\
\hphantom{c_{-10}=}{}
 +\alpha y_{12}\otimes
\widehat x_5\wedge \widehat x_9 + \alpha y_{14}\otimes \widehat x_4\wedge \widehat x_5 +\alpha y_{14}\otimes \widehat x_{11}\wedge \widehat y_3\nonumber\\
\hphantom{c_{-10}=}{}
 +y_{15}\otimes \widehat x_1\wedge \widehat x_4 + y_{15}\otimes
\widehat x_9\wedge \widehat y_3,\nonumber\\
c_{-8}^1= (1+\alpha)\big(\alpha ^3x_4\otimes \widehat x_1\wedge \widehat x_{15}
+\alpha ^3 x_4\otimes \widehat x_5\wedge \widehat x_{14} +\alpha ^2 x_4\otimes
\widehat x_{10}\wedge \widehat x_{12} +\alpha ^3x_7\otimes \widehat x_6\wedge \widehat x_{15}
\nonumber\\
\hphantom{c_{-8}^1=}{}
+ \alpha^3x_7\otimes \widehat x_8\wedge \widehat x_{14} +\alpha ^2y_1\otimes \widehat x_8\wedge \widehat x_{10} +\alpha
^3 y_1\otimes \widehat x_{15}\wedge \widehat y_4 +\alpha
y_5\otimes \widehat x_6\wedge \widehat x_{10}\nonumber\\
\hphantom{c_{-8}^1=}{}
+ \alpha ^2y_5\otimes \widehat x_{14}\wedge \widehat y_4 + \alpha^2y_6\otimes
\widehat x_5\wedge \widehat x_{10} + \alpha ^3y_6\otimes \widehat x_{15}\wedge
\widehat y_7+\alpha y_8\otimes \widehat x_1\wedge \widehat x_{10} \nonumber\\
\hphantom{c_{-8}^1=}{}
 +\alpha ^2y_8\otimes
\widehat x_{14}\wedge \widehat y_7 +
 \alpha y_{10}\otimes \widehat x_{12}\wedge \widehat y_4 +\alpha y_{10}\otimes
\widehat x_{13}\wedge \widehat y_7 +\alpha y_{12}\otimes \widehat x_{10}\wedge \widehat y_4 \nonumber\\
\hphantom{c_{-8}^1=}{}
 +
\alpha y_{13}\otimes \widehat x_{10}\wedge \widehat y_7\big)+
 \alpha y_{10}\otimes \widehat x_1\wedge \widehat x_8 +\alpha y_{10}\otimes
\widehat x_5\wedge \widehat x_6 +\alpha y_{14}\otimes \widehat x_5\wedge \widehat y_4
\nonumber\\
\hphantom{c_{-8}^1=}{}
 +\alpha y_{14}\otimes \widehat x_8\wedge \widehat y_7 +
 y_{15}\otimes \widehat x_1\wedge
\widehat y_4 + y_{15}\otimes \widehat x_6\wedge \widehat y_7 + \alpha ^2\big(1+\alpha^2\big)x_7\otimes
\widehat x_{10}\wedge \widehat x_{13},\nonumber\\
c_{-8}^2= (1+\alpha) \big ( \alpha ^3x_1\otimes \widehat x_4\wedge \widehat x_{15} +\alpha
^2x_1\otimes \widehat x_{11}\wedge \widehat x_{12} +\alpha ^3 x_6\otimes
\widehat x_7\wedge \widehat x_{15} +\alpha ^2x_6\otimes \widehat x_{11}\wedge \widehat x_{13}\nonumber\\
\hphantom{c_{-8}^2=}{}
 + \alpha ^2y_2\otimes \widehat x_4\wedge \widehat x_{13} +\alpha^2y_2\otimes
\widehat x_7\wedge \widehat x_{12}+\alpha ^2 y_4\otimes \widehat x_2\wedge \widehat x_{13} +\alpha^2
y_4\otimes \widehat x_8\wedge \widehat x_{11} \nonumber\\
\hphantom{c_{-8}^2=}{}
+
\alpha ^3 y_4\otimes \widehat x_{15}\wedge \widehat y_1 +
 \alpha y_5\otimes \widehat x_7\wedge \widehat x_{11} +\alpha ^2y_7\otimes
\widehat x_2\wedge \widehat x_{12} + \alpha ^2y_7\otimes \widehat x_5\wedge \widehat x_{11}
\nonumber\\
\hphantom{c_{-8}^2=}{}
+\alpha
^3y_7\otimes \widehat x_{15}\wedge \widehat y_6 +\alpha y_8\otimes \widehat x_4\wedge
\widehat x_{11} +
 \alpha y_{11}\otimes \widehat x_4\wedge \widehat x_8 +\alpha y_{11}\otimes
\widehat x_5\wedge \widehat x_7\nonumber\\
\hphantom{c_{-8}^2=}{}
 +\alpha y_{11}\otimes \widehat x_{12}\wedge \widehat y_1 +\alpha
y_{11}\otimes \widehat x_{13}\wedge \widehat y_6\big) +
 \alpha y_{12}\otimes \widehat x_2\wedge \widehat x_7 +\alpha y_{12}\otimes
\widehat x_{11}\wedge \widehat y_1 \nonumber\\
\hphantom{c_{-8}^2=}{}
+\alpha y_{13}\otimes \widehat x_2\wedge \widehat x_4 +\alpha
y_{13}\otimes \widehat x_{11}\wedge \widehat y_6 +
 y_{15}\otimes
\widehat x_4\wedge \widehat y_1 + y_{15}\otimes \widehat x_7\wedge \widehat y_6,\nonumber\\
c_{-6}^1= (1+\alpha)\big(\alpha x_2\otimes \widehat x_4\wedge \widehat x_{13}
+\alpha x_2\otimes \widehat x_7\wedge \widehat x_{12} +\alpha x_5\otimes
\widehat x_4\wedge \widehat x_{14}+x_5\otimes \widehat x_9\wedge \widehat x_{12}\nonumber\\
\hphantom{c_{-6}^1=}{}
+\alpha
x_8\otimes \widehat x_7\wedge \widehat x_{14}+ x_8\otimes \widehat x_9\wedge \widehat x_{13} +y_1\otimes \widehat x_7\wedge \widehat x_9
+y_4\otimes \widehat x_6\wedge \widehat x_9 +\alpha y_4\otimes \widehat x_{13}\wedge \widehat y_2\nonumber\\
\hphantom{c_{-6}^1=}{}
+\alpha ^2 y_4\otimes \widehat x_{14}\wedge \widehat y_5
+ y_6\otimes \widehat x_4\wedge \widehat x_9 + y_7\otimes \widehat x_1\wedge \widehat x_9 +\alpha
y_7\otimes \widehat x_{12}\wedge \widehat y_2 +\alpha
^2 y_7\otimes \widehat x_{14}\wedge \widehat y_8\nonumber\\
\hphantom{c_{-6}^1=}{}
 + y_9\otimes \widehat x_1\wedge \widehat x_7 +
 y_9\otimes \widehat x_4\wedge \widehat x_6 +\alpha y_9\otimes \widehat x_{12}\wedge \widehat y_5
+\alpha y_9\otimes \widehat x_{13}\wedge \widehat y_8\big) +y_{12}\otimes
\widehat x_7\wedge \widehat y_2\nonumber\\
\hphantom{c_{-6}^1=}{}
 + y_{12}\otimes \widehat x_9\wedge \widehat y_5
+
 y_{13}\otimes \widehat x_4\wedge \widehat y_2 +y_{13}\otimes \widehat x_9\wedge \widehat y_8 +
y_{14}\otimes \widehat x_4\wedge \widehat y_5 + y_{14}\otimes \widehat x_7\wedge \widehat y_8,\nonumber\\
c_{-6}^2= (1+\alpha)\big(\alpha ^2x_1\otimes \widehat x_8\wedge \widehat x_{10}
+\alpha ^3x_1\otimes \widehat x_{15}\wedge \widehat y_4 +\alpha ^2 x_4\otimes
\widehat x_2\wedge \widehat x_{13} +\alpha ^2 x_4\otimes \widehat x_8\wedge \widehat x_{11}\nonumber\\
\hphantom{c_{-6}^2=}{}
+ \alpha ^3x_4\otimes \widehat x_{15}\wedge \widehat y_1 + \alpha^3 x_9\otimes
\widehat x_3\wedge \widehat x_{15} +\alpha ^2 x_9\otimes \widehat x_8\wedge \widehat x_{13} +\alpha
^2 y_2\otimes \widehat x_3\wedge \widehat x_{10} \nonumber\\
\hphantom{c_{-6}^2=}{}
+\alpha ^2
y_2\otimes \widehat x_{13}\wedge \widehat y_4 +
 \alpha ^2y_3\otimes \widehat x_2\wedge \widehat x_{10} +\alpha^2 y_3\otimes
\widehat x_5\wedge \widehat x_8 +\alpha^3 y_3\otimes \widehat x_{15}\wedge \widehat y_9
 \nonumber\\
\hphantom{c_{-6}^2=}{}
 +\alpha
y_5\otimes \widehat x_3\wedge \widehat x_8 +\alpha y_8\otimes \widehat x_3\wedge \widehat x_5
+
 \alpha y_8\otimes \widehat x_{10}\wedge \widehat y_1 +\alpha y_8\otimes
\widehat x_{11}\wedge \widehat y_4
 \nonumber\\
\hphantom{c_{-6}^2=}{}
+\alpha y_8\otimes \widehat x_{13}\wedge \widehat y_9 +\alpha
y_{11}\otimes \widehat x_8\wedge \widehat y_4\big) +\alpha y_{13}\otimes \widehat x_2\wedge \widehat y_4 +
 \alpha y_{10}\otimes \widehat x_2\wedge \widehat x_3
 \nonumber\\
\hphantom{c_{-6}^2=}{}
 +\alpha y_{10}\otimes
\widehat x_8\wedge \widehat y_1 +\alpha y_{13}\otimes \widehat x_8\wedge \widehat y_9 +y_{15}\otimes
\widehat x_3\wedge
\widehat y_9 + y_{15}\otimes \widehat y_1\wedge \widehat y_4,\nonumber\\
c_{-4}^1= (1+\alpha)\big(\alpha x_2\otimes \widehat x_3\wedge \widehat x_{10}
+\alpha x_2\otimes \widehat x_{13}\wedge \widehat y_4 +x_4\otimes \widehat x_6\wedge
\widehat x_9 +\alpha x_4\otimes \widehat x_{13}\wedge \widehat y_2 \nonumber\\
\hphantom{c_{-4}^1=}{}
+ \alpha ^2 x_4\otimes \widehat x_{14}\wedge \widehat y_5 +x_5\otimes \widehat x_6\wedge
\widehat x_{10} +\alpha x_5\otimes \widehat x_{14}\wedge \widehat y_4 + \alpha x_{11}\otimes
\widehat x_3\wedge \widehat x_{14}\nonumber\\
\hphantom{c_{-4}^1=}{}
+x_{11}\otimes \widehat x_6\wedge \widehat x_{13} +
 y_1\otimes \widehat x_3\wedge \widehat x_6 +y_3\otimes \widehat x_1\wedge \widehat x_6 +\alpha
y_3\otimes \widehat x_{10}\wedge \widehat y_2 +\alpha ^2y_3\otimes \widehat x_{14}\wedge
\widehat y_{11}\nonumber\\
\hphantom{c_{-4}^1=}{}
 +y_6\otimes \widehat x_1\wedge \widehat x_3 +
 y_6\otimes \widehat x_9\wedge \widehat y_4 +\alpha y_6\otimes \widehat x_{10}\wedge \widehat y_5
+\alpha y_6\otimes \widehat x_{13}\wedge \widehat y_{11} + y_9\otimes \widehat x_6\wedge
\widehat y_4\big) \nonumber\\
\hphantom{c_{-4}^1=}{}
+y_{10}\otimes \widehat x_3\wedge \widehat y_2 +
 y_{10}\otimes \widehat x_6\wedge \widehat y_5 +y_{13}\otimes \widehat x_6\wedge \widehat y_{11}
+y_{13}\otimes \widehat y_2\wedge \widehat y_4 +y_{14}\otimes \widehat x_3\wedge \widehat y_{11}
 \nonumber\\
\hphantom{c_{-4}^1=}{}
+y_{14}\otimes \widehat y_4\wedge
\widehat y_5,\nonumber\\
c_{-4}^2= (1+\alpha)\big(\alpha ^2 x_3\otimes \widehat x_2\wedge \widehat x_{10}
+\alpha ^2 x_3\otimes \widehat x_5\wedge \widehat x_8 +\alpha ^3 x_3\otimes
\widehat x_{15}\wedge \widehat y_9 +\alpha ^2 x_6\otimes \widehat x_5\wedge \widehat x_{10}\nonumber\\
\hphantom{c_{-4}^2=}{}
+ \alpha ^3 x_6\otimes \widehat x_{15}\wedge \widehat y_7 + \alpha ^2x_7\otimes
\widehat x_2\wedge \widehat x_{12} +\alpha ^2 x_7\otimes \widehat x_5\wedge \widehat x_{11} + \alpha
^3 x_7\otimes \widehat x_{15}\wedge \widehat y_6\nonumber\\
\hphantom{c_{-4}^2=}{}
 +\alpha ^2 x_9\otimes \widehat x_5\wedge
\widehat x_{12} +
 \alpha ^3 x_9\otimes \widehat x_{15}\wedge \widehat y_3 + \alpha ^2 y_2\otimes
\widehat x_{10}\wedge \widehat y_3 +\alpha ^2 y_2\otimes \widehat x_{12}\wedge \widehat y_7
\nonumber\\
\hphantom{c_{-4}^2=}{}
 +\alpha
y_5\otimes \widehat x_8\wedge \widehat y_3 +\alpha y_5\otimes \widehat x_{10}\wedge \widehat y_6
+ \alpha y_5\otimes \widehat x_{11}\wedge \widehat y_7 +\alpha y_5\otimes
\widehat x_{12}\wedge \widehat y_9 +\alpha y_8\otimes \widehat x_5\wedge \widehat y_3
\nonumber\\
\hphantom{c_{-4}^2=}{}
+\alpha
y_{11}\otimes \widehat x_5\wedge \widehat y_7\big) +\alpha y_{10}\otimes
\widehat x_2\wedge \widehat y_3
+ \alpha y_{10}\otimes \widehat x_5\wedge \widehat y_6 +\alpha y_{12}\otimes
\widehat x_2\wedge \widehat y_7
\nonumber\\
\hphantom{c_{-4}^2=}{}
 +\alpha y_{12}\otimes \widehat x_5\wedge \widehat y_9 +
y_{15}\otimes \widehat y_3\wedge \widehat y_9 +y_{15}\otimes \widehat y_6\wedge
\widehat y_7,\nonumber\\
c_{-2}= (1+\alpha)\big(\alpha x_2\otimes \widehat x_{10}\wedge \widehat y_3 +\alpha
x_2\otimes \widehat x_{12}\wedge \widehat y_7 + x_3\otimes \widehat x_1\wedge \widehat x_6 +\alpha
x_3\otimes \widehat x_{10}\wedge \widehat y_2\nonumber\\
\hphantom{c_{-2}=}{}
 +\alpha ^2x_3\otimes
\widehat x_{14}\wedge \widehat y_{11} +
 x_7\otimes \widehat x_1\wedge \widehat x_9 +\alpha x_7\otimes \widehat x_{12}\wedge \widehat y_2
+\alpha ^2 x_7\otimes \widehat x_{14}\wedge \widehat y_8
\nonumber\\
\hphantom{c_{-2}=}{}
 +x_8\otimes \widehat x_1\wedge
\widehat x_{10} +\alpha x_8\otimes \widehat x_{14}\wedge \widehat y_7 +
 x_{11}\otimes \widehat x_1\wedge \widehat x_{12} +\alpha x_{11}\otimes
\widehat x_{14}\wedge \widehat y_3 +y_1\otimes \widehat x_6\wedge \widehat y_3
\nonumber\\
\hphantom{c_{-2}=}{}
 +y_1\otimes
\widehat x_9\wedge \widehat y_7 +\alpha
y_1\otimes \widehat x_{10}\wedge \widehat y_8 +
 \alpha y_1\otimes \widehat x_{12}\wedge \widehat y_{11} +y_6\otimes \widehat x_1\wedge \widehat y_3
+
y_9\otimes \widehat x_1\wedge \widehat y_7\big)\nonumber\\
\hphantom{c_{-2}=}{}
+ y_{10}\otimes \widehat x_1\wedge \widehat y_8 + y_{10}\otimes \widehat y_2\wedge \widehat y_3
+y_{12}\otimes \widehat x_1\wedge \widehat y_{11} + y_{12}\otimes \widehat y_2\wedge \widehat y_7
+y_{14}\otimes \widehat y_3\wedge \widehat y_{11}
\nonumber\\
\hphantom{c_{-2}=}{}
+ y_{14}\otimes \widehat y_7\wedge
\widehat y_8,\label{eqpsl4c3}\\
c_{0}= \alpha h_1\otimes \widehat x_{14}\wedge \widehat y_{14} +\alpha ^2h_1\otimes
\widehat x_{15}\wedge \widehat y_{15} +h_2\otimes \widehat x_5\wedge \widehat y_5 +h_2\otimes
\widehat x_8\wedge \widehat y_8 +h_2\otimes \widehat x_{10}\wedge \widehat y_{10}
\nonumber\\
\hphantom{c_{0}=}{}
+ h_2\otimes \widehat x_{11}\wedge \widehat y_{11} +h_2\otimes \widehat x_{12}\wedge \widehat y_{12}
+h_2\otimes \widehat x_{13}\wedge \widehat y_{13} +\alpha ^2h_2\otimes
\widehat x_{14}\wedge \widehat y_{14} +h_3\otimes \widehat x_1\wedge \widehat y_1\nonumber\\
\hphantom{c_{0}=}{}
+ \alpha h_3\otimes \widehat x_5\wedge \widehat y_5 +h_3\otimes \widehat x_6\wedge \widehat y_6
+\alpha h_3\otimes \widehat x_8\wedge \widehat y_8 +h_3\otimes \widehat x_9\wedge \widehat y_9
+\alpha h_3\otimes \widehat x_{10}\wedge \widehat y_{10}\nonumber\\
\hphantom{c_{0}=}{}
 + \alpha h_3\otimes \widehat x_{11}\wedge \widehat y_{11} +\alpha h_3\otimes
\widehat x_{12}\wedge \widehat y_{12} +\alpha ^2(1+\alpha) h_3\otimes \widehat x_{14}\wedge
\widehat y_{14}
\nonumber\\
\hphantom{c_{0}=}{}
+\alpha ^3(1+\alpha)
h_3\otimes \widehat x_{15}\wedge \widehat y_{15} +
 \alpha h_4\otimes \widehat x_{12}\wedge \widehat y_{12} + \alpha h_4\otimes
\widehat x_{13}\wedge \widehat y_{13} +\alpha h_4\otimes \widehat x_{14}\wedge \widehat y_{14}
\nonumber\\
\hphantom{c_{0}=}{}
+\alpha ^2 h_4\otimes \widehat x_{15}\wedge \widehat y_{15} +x_1\otimes \widehat h_1\wedge
\widehat x_1 +
 \alpha x_1\otimes \widehat x_{10}\wedge \widehat y_8 +\alpha x_1\otimes
\widehat x_{12}\wedge \widehat y_{11}
\nonumber\\
\hphantom{c_{0}=}{}
 +\alpha x_1\otimes \widehat x_{14}\wedge \widehat y_{13}
+x_2\otimes \widehat h_1\wedge \widehat x_2 +x_2\otimes \widehat x_5\wedge \widehat y_1 +
 x_2\otimes \widehat x_8\wedge \widehat y_6 +x_2\otimes \widehat x_{11}\wedge \widehat y_9
\nonumber\\
\hphantom{c_{0}=}{}
+\alpha^3x_2 \otimes \widehat x_{15}\wedge \widehat y_{14} +\alpha x_3\otimes
\widehat x_{13}\wedge \widehat y_{12} +x_4\otimes \widehat h_1\wedge \widehat x_4 +
 \alpha x_4\otimes \widehat x_{12}\wedge \widehat y_{10}\nonumber\\
\hphantom{c_{0}=}{}
 +x_5\otimes \widehat x_{10}\wedge
\widehat y_6+x_5\otimes \widehat x_{12}\wedge \widehat y_9 +\alpha x_5\otimes \widehat x_{15}\wedge \widehat y_{13} +x_6\otimes \widehat h_1\wedge \widehat x_6
+
 \alpha x_6\otimes \widehat x_{10}\wedge \widehat y_5 \nonumber\\
\hphantom{c_{0}=}{}
+\alpha x_6\otimes
\widehat x_{13}\wedge \widehat y_{11} +\alpha x_6\otimes \widehat x_{14}\wedge \widehat y_{12}
+x_7\otimes \widehat h_1\wedge \widehat x_7 +\alpha x_7\otimes \widehat x_{13}\wedge \widehat y_{10}\nonumber\\
\hphantom{c_{0}=}{}
+ x_8\otimes \widehat x_{10}\wedge \widehat y_1 +x_8\otimes \widehat x_{13}\wedge \widehat y_9
+\alpha x_8\otimes \widehat x_{15}\wedge \widehat y_{12} +\alpha x_9\otimes
\widehat x_{12}\wedge \widehat y_5 +\alpha x_9\otimes
\widehat x_{13}\wedge \widehat y_8\nonumber\\
\hphantom{c_{0}=}{}
+ \alpha x_9\otimes \widehat x_{14}\wedge \widehat y_{10} +x_{10}\otimes \widehat h_1\wedge
\widehat x_{10} +x_{10}\otimes \widehat x_{14}\wedge \widehat y_9 +\alpha x_{10}\otimes
\widehat x_{15}\wedge \widehat y_{11}\nonumber\\
\hphantom{c_{0}=}{}
+x_{11}\otimes \widehat h_1\wedge \widehat x_{11}
+ x_{11}\otimes \widehat x_{12}\wedge \widehat y_1 + x_{11}\otimes \widehat x_{13}\wedge \widehat y_6
+\alpha x_{11}\otimes \widehat x_{15}\wedge \widehat y_{10}
\nonumber\\
\hphantom{c_{0}=}{}
 +x_{12}\otimes
\widehat x_{14}\wedge \widehat y_6 + \alpha x_{12}\otimes \widehat x_{15}\wedge \widehat y_8
+
 x_{13}\otimes \widehat x_{14}\wedge \widehat y_1 + \alpha x_{13}\otimes
\widehat x_{15}\wedge \widehat y_5
\nonumber\\
\hphantom{c_{0}=}{}
+x_{14}\otimes \widehat h_1\wedge \widehat x_{14} + y_1\otimes
\widehat h_1\wedge \widehat y_1 +\alpha y_1\otimes \widehat x_8\wedge \widehat y_{10} +
 \alpha y_1\otimes \widehat x_{11}\wedge \widehat y_{12} +\alpha y_1\otimes
\widehat x_{13}\wedge \widehat y_{14}
\nonumber\\
\hphantom{c_{0}=}{}
+y_2\otimes \widehat h_1\wedge
\widehat y_2 +y_2\otimes \widehat x_1\wedge \widehat y_5 +y_2\otimes \widehat x_6\wedge \widehat y_8 +
 y_2\otimes \widehat x_9\wedge \widehat y_{11} +\alpha^3y_2\otimes \widehat x_{14}\wedge
\widehat y_{15}
\nonumber\\
\hphantom{c_{0}=}{}
+ \alpha y_3\otimes \widehat x_{12}\wedge \widehat y_{13} +y_4\otimes
\widehat h_1\wedge \widehat y_4 +\alpha y_4\otimes \widehat x_{10}\wedge \widehat y_{12}+
 y_5\otimes \widehat x_6\wedge \widehat y_{10} +y_5\otimes \widehat x_9\wedge \widehat y_{12}
\nonumber\\
\hphantom{c_{0}=}{}
+\alpha y_5\otimes \widehat x_{13}\wedge \widehat y_{15} +y_6\otimes \widehat h_1\wedge
\widehat y_6 +\alpha y_6\otimes \widehat x_5\wedge \widehat y_{10} +
 \alpha y_6\otimes \widehat x_{11}\wedge \widehat y_{13}
\nonumber\\
\hphantom{c_{0}=}{}
 +\alpha y_6\otimes
\widehat x_{12}\wedge \widehat y_{14} +y_7\otimes \widehat h_1\wedge \widehat y_7 +\alpha y_7\otimes
\widehat x_{10}\wedge \widehat y_{13} +
y_8\otimes \widehat x_1\wedge \widehat y_{10} +
 y_8\otimes \widehat x_9\wedge \widehat y_{13}
\nonumber\\
\hphantom{c_{0}=}{}
 +\alpha y_8\otimes \widehat x_{12}\wedge
\widehat y_{15} +\alpha y_9\otimes \widehat x_5\wedge \widehat y_{12} +\alpha y_9\otimes
\widehat x_8\wedge \widehat y_{13} +\alpha
y_9\otimes \widehat x_{10}\wedge \widehat y_{14}
\nonumber\\
\hphantom{c_{0}=}{}
+
 y_{10}\otimes \widehat h_1\wedge \widehat y_{10} +y_{10}\otimes \widehat x_9\wedge \widehat y_{14}
+\alpha y_{10}\otimes \widehat x_{11}\wedge \widehat y_{15} + y_{11}\otimes
\widehat h_1\wedge \widehat y_{11} +y_{11}\otimes \widehat x_1\wedge \widehat y_{12}
\nonumber\\
\hphantom{c_{0}=}{}
+
 y_{11}\otimes \widehat x_6\wedge \widehat y_{13} +\alpha y_{11}\otimes
\widehat x_{10}\wedge
\widehat y_{15} +y_{12}\otimes \widehat x_6\wedge \widehat y_{14} +\alpha y_{12}\otimes \widehat x_8\wedge \widehat y_{15}\nonumber\\
\hphantom{c_{0}=}{}
+ y_{13}\otimes \widehat x_1\wedge \widehat y_{14} +\alpha y_{13}\otimes \widehat x_5\wedge
\widehat y_{15} + y_{14}\otimes \widehat h_1\wedge \widehat y_{14}.\label{eqwk4ab}
\end{gather}
\end{Lemma}

\begin{Lemma}\label{L5.7bg} For $\fbgl(4;\alpha)$, we take the same Cartan matrix and the
basis as for $\fwk(4;\alpha)$ with all
Chevalley generators even, except $x_1, y_1$.
Let $c_i$ be a~cocycle of~$\fwk(4;\alpha)$, see~\eqref{eqpsl4c3},
\eqref{eqwk4ab}; let $\overline c_i$ be the
 corresponding cocycle of~$\fbgl(4;\alpha)$. We see that
\begin{gather*}
 \overline c_{-12} = i(c_{-12}),\qquad \overline c_{-10} = i(c_{-10}), \\
 \overline c_{-8}^{1} = i\big(c_{-8}^{1}\big),\qquad
 \overline c_{-8}^{2} = i\big(c_{-8}^{2}\big) + \alpha^2(1+\alpha) h_3 \otimes (\widehat x_{11})^{\wedge 2},\\
 \overline c_{-6}^{1} = i\big(c_{-6}^{1}\big) + (1+\alpha) h_3 \otimes (\widehat x_9)^{\wedge 2},\qquad
 \overline c_{-6}^{2} = i\big(c_{-6}^{2}\big) + \alpha^2(1+\alpha)(h_3 + h_4)\otimes (\widehat x_8)^{\wedge 2},\\
 \overline c_{-4}^{1} = i\big(c_{-4}^{1}\big) + (1+\alpha)(h_3 + h_4)\otimes (\widehat x_6)^{\wedge 2},\qquad
 \overline c_{-4}^{2} = i\big(c_{-4}^{2}\big) + (1+\alpha)(h_3 + h_4)\otimes (\widehat x_6)^{\wedge 2},\\
 \overline c_{-2}= i(c_{-2}) + (1+\alpha) h_4 \otimes (\widehat x_1)^{\wedge 2},\qquad
 \overline c_0= i(c_0).
 \end{gather*}
\end{Lemma}

\subsection[$\mathfrak{gl}(a\vert a+pk)$ and its simple relatives (by A.~Krutov)]{$\boldsymbol{\mathfrak{gl}(a\vert a+pk)}$ and its simple relatives (by A.~Krutov)}\label{ssGL}

The following results are incomplete for $p=2$: The cases
$\mathfrak{gl}(1|7)$ and $\mathfrak{psl}(1|7)$; $\mathfrak{gl}(2|4)$ and $\mathfrak{psl}(2|4)$;
$\mathfrak{gl}(2|6)$ and $\mathfrak{psl}(2|6)$; $\mathfrak{gl}(3|5)$ and $\mathfrak{psl}(3|5)$; $\mathfrak{gl}(4|4)$ and $\mathfrak{psl}(4|4)$ will be studied elsewhere,
together with the superizations of Chebochko's results on the other Lie algebras with $ADE$ root systems.

Here we list what A.~Krutov computed at our request. For $\mathfrak{gl}(n)$ with any $n>1$, the element~$h_n$ acts on $x_i$ as $[h_n, x_i] = \delta_{i,1}x_i$.

$\bullet$ $\mathfrak{gl}(5)$, positive elements of Chevalley basis are
\begin{gather*}
  x_1, \
  x_2, \
  x_3, \
  x_4, \
  x_5=[x_1, x_2], \
  x_6=[x_2, x_3], \
  x_7=[x_3, x_4], \
  x_8=[x_3, [x_1,x_2]], \\
  x_9=[x_4, [x_2,x_3]], \
  x_{10}=[[x_1,x_2], [x_3,x_4]].
\end{gather*}

$\bullet$ $\mathfrak{gl}(6)$, positive elements of Chevalley basis are
\begin{gather*}
  x_1, \
  x_2, \
  x_3, \
  x_4, \
  x_5, \
  x_6=[x_1, x_2], \
  x_7=[x_2, x_3], \
  x_8=[x_3, x_4], \
  x_9=[x_4, x_5],\\
  x_{10}=[x_3, [x_1,x_2]], \
  x_{11}=[x_4, [x_2,x_3]], \
  x_{12}=[x_5, [x_3,x_4]], \
  x_{13}=[[x_1,x_2], [x_3,x_4]], \\
  x_{14}=[[x_2,x_3], [x_4,x_5]], \
  x_{15}=[[x_4,x_5], [x_3,[x_1,x_2]]].
\end{gather*}

$\bullet$ $p=5$:
\begin{itemize}\itemsep=0pt
\item[--]$\fg=\mathfrak{gl}(5)$ and $\fg=\mathfrak{psl}(5)$: in both cases $H^2(\fg;\fg)=0$.
\item[--] $\fg=\mathfrak{gl}(1|6)$ and $\fg=\mathfrak{psl}(1|6)$: in both cases $H^2(\fg;\fg)=0$.
\item[--] $\fg=\mathfrak{gl}(2|7)$ and $\fg=\mathfrak{psl}(2|7)$: in both %${}^*$
 cases $H^2(\fg;\fg)=0$.
\end{itemize}

$\bullet$ $p=3$:
\begin{itemize}\itemsep=0pt
\item[--] $\fg=\mathfrak{gl}(3)$ and $\fg=\mathfrak{psl}(3)$: in both %${}^*$
 cases $H^2(\fg;\fg)=0$.

\item[--]$\fg=\mathfrak{gl}(6)$ and for $\fg=\mathfrak{psl}(6)$: in both %${}^*$
 cases $H^2(\fg;\fg)=0$.

\item[--]$\fg=\mathfrak{gl}(1|4)$: $H^2(\fg; \fg)=0$; for $\fg:=\mathfrak{psl}(1|4)$ the space
$H^2(\fg;\fg)$ is spanned by (observe that
indices of the cocycles are the degrees with respect to the standard $\Zee$-grading)
\begin{gather*}
 c_{-7} =
 2 y_2\otimes\widehat x_6\wedge \widehat x_{10}+y_2\otimes\widehat x_8\wedge \widehat x_9+2 y_5\otimes\widehat x_6\wedge \widehat x_9+2 y_6\otimes\widehat x_2\wedge \widehat x_{10}\\
 \hphantom{c_{-7} =}{}
+2 y_6\otimes\widehat x_5\wedge \widehat x_9 +y_8\otimes\widehat x_2\wedge \widehat x_9
+2 y_9\otimes\widehat x_2\wedge \widehat x_8+2 y_9\otimes\widehat x_5\wedge \widehat x_6\\
\hphantom{c_{-7} =}{}
+y_{10}\otimes\widehat x_2\wedge \widehat x_6,\\
 c_{-4} =
2 y_1\otimes\widehat x_3\wedge \widehat x_7+2 y_3\otimes\widehat x_1\wedge \widehat x_7+y_3\otimes\widehat x_{10}\wedge \widehat y_2+2 y_7\otimes\widehat x_1\wedge \widehat x_3+2 y_7\otimes\widehat x_8\wedge \widehat y_2\\
\hphantom{c_{-4} =}{}
+y_8\otimes\widehat x_7\wedge \widehat y_2
+2 y_{10}\otimes\widehat x_3\wedge \widehat y_2+2 x_2\otimes\widehat x_3\wedge \widehat x_{10}+2 x_2\otimes\widehat x_7\wedge \widehat x_8,\\
 c_{-2} =
 2 y_4\otimes\widehat x_7\wedge \widehat y_1+y_4\otimes\widehat x_9\wedge \widehat y_5+2 y_7\otimes\widehat x_4\wedge \widehat y_1+y_7\otimes\widehat x_9\wedge \widehat y_8+2 y_9\otimes\widehat x_4\wedge \widehat y_5\\
 \hphantom{c_{-2} =}{}
+y_9\otimes\widehat x_7\wedge \widehat y_8+x_1\otimes\widehat x_4\wedge \widehat x_7+x_5\otimes\widehat x_4\wedge \widehat x_9+x_8\otimes\widehat x_7\wedge \widehat x_9,\\
 c_{-1} =
 x_3\otimes\widehat x_{10}\wedge \widehat y_6+x_6\otimes\widehat x_{10}\wedge \widehat y_3+2 y_1\otimes\widehat x_4\wedge \widehat y_3+2 y_4\otimes\widehat x_1\wedge \widehat y_3+2 y_4\otimes\widehat x_5\wedge \widehat y_6\\
 \hphantom{c_{-1} =}{}
 +y_5\otimes\widehat x_4\wedge \widehat y_6
+x_3\otimes\widehat x_1\wedge \widehat x_4+2 x_6\otimes\widehat x_4\wedge \widehat x_5+2 y_{10}\otimes\widehat y_3\wedge \widehat y_6.
\end{gather*}
\item[--] $\fg=\mathfrak{gl}(1|7)$ and $\fg=\mathfrak{psl}(1|7)$: in both cases $H^2(\fg;\fg)=0$.
\item[--] $\fg=\mathfrak{gl}(2|5)$ and $\fg=\mathfrak{psl}(2|5)$: in both cases $H^2(\fg;\fg)=0$.
For $\fg=\mathfrak{gl}(3|3)$, we have $H^2(\fg;\fg)=0$; for $\fg=\mathfrak{psl}(3|3)$, the space $H^2(\fg;\fg)$ is spanned by (observe that
indices of the cocycles are the degrees in the standard $\Zee$-grading)
\begin{gather*}
 c_{-9} =
 2 y_3\otimes\widehat x_{11}\wedge \widehat x_{15}+2 y_3\otimes\widehat x_{13}\wedge \widehat x_{14}+y_7\otimes\widehat x_8\wedge \widehat x_{15}+y_7\otimes\widehat x_{12}\wedge \widehat x_{13}\\
\hphantom{c_{-9} =}{}
+y_8\otimes\widehat x_7\wedge \widehat x_{15} +y_8\otimes\widehat x_{10}\wedge \widehat x_{14}
+2 y_{10}\otimes\widehat x_8\wedge \widehat x_{14}+y_{10}\otimes\widehat x_{11}\wedge \widehat x_{12}\\
\hphantom{c_{-9} =}{}
+2 y_{11}\otimes\widehat x_3\wedge \widehat x_{15}+y_{11}\otimes\widehat x_{10}\wedge \widehat x_{12}+2 y_{12}\otimes\widehat x_7\wedge \widehat x_{13}+y_{12}\otimes\widehat x_{10}\wedge \widehat x_{11}\\
\hphantom{c_{-9} =}{} +y_{13}\otimes\widehat x_3\wedge \widehat x_{14}+y_{13}\otimes\widehat x_7\wedge \widehat x_{12}+y_{14}\otimes\widehat x_3\wedge \widehat x_{13}+y_{14}\otimes\widehat x_8\wedge \widehat x_{10}\\
\hphantom{c_{-9} =}{}
+2 y_{15}\otimes\widehat x_3\wedge \widehat x_{11}+2 y_{15}\otimes\widehat x_7\wedge \widehat x_8.
\end{gather*}
\item[--]
$\fg=\mathfrak{gl}(3|6)$ and $\fg=\mathfrak{psl}(3|6)$: in both cases $H^2(\fg;\fg)=0$.
\end{itemize}

$\bullet$ $p=2$:
\begin{itemize}
\item[--]
$\fg=\mathfrak{gl}(4)$ and $\fg=\mathfrak{psl}(4)$, see Chebochko's Theorem~\ref{NC2}.
\end{itemize}

Recently Chebochko and Kuznetsov published a~realization of the two inequivalent deforms
of $\mathfrak{psl}(4)$ (they keep calling it ``$G_2$''; it is $\overline{A_3}$ in \eqref{NC2}),
see \cite{ChKu} which is an expansion of
\cite{Ch1}, where Chebochko proved that the 20-dimensional, see \eqref{NC2},
space of cocycles constitutes
two (apart from the origin) orbits relative $\operatorname{Aut}(\mathfrak{psl}(4))$, and
described the two respective Lie algebras. A~posteriori this result
is natural because
the deforms of the Hamiltonian Lie algebra $\mathfrak{h}_\Pi(2k;\un)$ (classified by Skryabin, see \cite{Sk0, SkH}) exist for $p=2$, the derived algebra of $\mathfrak{h}_\Pi(2k;\un)$
can inherit these deforms, and $\mathfrak{psl}(4)\simeq\mathfrak{h}_\Pi^{(1)}(4;(1111)|0)$. Here are explicit cocycles in the Chevalley basis given by formula~\eqref{gl4cm}:
\begin{itemize}\itemsep=0pt
\item[--]$\mathfrak{psl}(4)\simeq \textbf{F}(\mathfrak{h}_\Pi^{(1)}(0|4))$ (an exceptional case due to this occasional isomorphism; should be compared with quantization of $\mathfrak{po}_\Pi(4;\One|0)$ and the result in \cite{KST}).
For $\fg=\mathfrak{psl}(4)$, the space $H^2(\fg;\fg)$ is spanned by
\begin{gather*}
 c_{-6}^{1} =
 y_1\otimes \widehat x_4\wedge \widehat x_6+y_4\otimes \widehat x_1\wedge \widehat x_6+y_6\otimes \widehat x_1\wedge \widehat x_4,\\
 c_{-6}^{2} =
 y_3\otimes \widehat x_5\wedge \widehat x_6+y_5\otimes \widehat x_3\wedge \widehat x_6+y_6\otimes \widehat x_3\wedge \widehat x_5,\\
 c_{-4}^{1} =
 y_2\otimes \widehat h_3\wedge \widehat x_6+y_4\otimes \widehat h_3\wedge \widehat x_5+y_5\otimes \widehat h_3\wedge \widehat x_4+y_6\otimes \widehat h_3\wedge \widehat x_2+h_2\otimes \widehat x_4\wedge \widehat x_5\\
 \hphantom{c_{-4}^{1} =}{}
 +y_2\otimes \widehat x_1\wedge \widehat x_5+y_2\otimes \widehat x_3\wedge \widehat x_4
 +y_6\otimes \widehat x_4\wedge \widehat y_1+y_6\otimes \widehat x_5\wedge \widehat y_3,\\
 c_{-4}^{2} =
 y_2\otimes \widehat h_2\wedge \widehat x_6+y_2\otimes \widehat h_3\wedge \widehat x_6+y_4\otimes \widehat h_2\wedge \widehat x_5+y_4\otimes \widehat h_3\wedge \widehat x_5+y_5\otimes \widehat h_2\wedge \widehat x_4\\
 \hphantom{c_{-4}^{2} =}{}
 +y_5\otimes \widehat h_3\wedge \widehat x_4+y_6\otimes \widehat h_2\wedge \widehat x_2
 +y_6\otimes \widehat h_3\wedge \widehat x_2+h_3\otimes \widehat x_2\wedge \widehat x_6+y_2\otimes \widehat x_1\wedge \widehat x_5\\
 \hphantom{c_{-4}^{2} =}{}
 +y_2\otimes \widehat x_3\wedge \widehat x_4+y_4\otimes \widehat x_2\wedge \widehat x_3+y_5\otimes \widehat x_1\wedge \widehat x_2+x_1\otimes \widehat x_4\wedge \widehat x_6+x_3\otimes \widehat x_5\wedge \widehat x_6,\\
c_{-2}^{1} =
 y_2\otimes \widehat x_4\wedge \widehat y_3+y_4\otimes \widehat x_2\wedge \widehat y_3+x_3\otimes \widehat x_2\wedge \widehat x_4,\\
c_{-2}^{2} =
 y_2\otimes \widehat x_5\wedge \widehat y_1+y_5\otimes \widehat x_2\wedge \widehat y_1+x_1\otimes \widehat x_2\wedge \widehat x_5,\\
c_{-2}^{3} =
 h_3\otimes \widehat x_6\wedge \widehat y_2+y_1\otimes \widehat h_2\wedge \widehat x_3+y_3\otimes \widehat h_2\wedge \widehat x_1+x_2\otimes \widehat h_2\wedge \widehat x_6+y_6\otimes \widehat h_2\wedge \widehat y_2\\
 \hphantom{c_{-2}^{3} =}{}
 +y_1\otimes \widehat x_5\wedge \widehat y_2+y_1\otimes \widehat x_6\wedge \widehat y_4+y_3\otimes \widehat x_4\wedge \widehat y_2+y_3\otimes \widehat x_6\wedge \widehat y_5,\\
c_{-2}^{4} =
 y_1\otimes \widehat h_3\wedge \widehat x_3+y_3\otimes \widehat h_3\wedge \widehat x_1+h_2\otimes \widehat x_1\wedge \widehat x_3+x_2\otimes \widehat h_3\wedge \widehat x_6+y_6\otimes \widehat h_3\wedge \widehat y_2\\
 \hphantom{c_{-2}^{4} =}{}
 +y_1\otimes \widehat x_5\wedge \widehat y_2+y_1\otimes \widehat x_6\wedge \widehat y_4+y_3\otimes \widehat x_4\wedge \widehat y_2
 +y_3\otimes \widehat x_6\wedge \widehat y_5+y_4\otimes \widehat x_3\wedge \widehat y_2\\
\hphantom{c_{-2}^{4} =}{}
 +y_5\otimes \widehat x_1\wedge \widehat y_2+x_4\otimes \widehat x_1\wedge \widehat x_6+x_5\otimes \widehat x_3\wedge \widehat x_6,\\
c_{0}^{1} =
 h_3\otimes \widehat x_4\wedge \widehat y_5+x_3\otimes \widehat h_2\wedge \widehat x_1+x_5\otimes \widehat h_2\wedge \widehat x_4+y_1\otimes \widehat h_2\wedge \widehat y_3+y_4\otimes \widehat h_2\wedge \widehat y_5\\
\hphantom{c_{0}^{1} =}{}
 +x_3\otimes \widehat x_4\wedge \widehat y_2+x_3\otimes \widehat x_6\wedge \widehat y_5+y_1\otimes \widehat x_2\wedge \widehat y_5+y_1\otimes \widehat x_4\wedge \widehat y_6,\\
c_{0}^{2} =
 h_3\otimes \widehat x_5\wedge \widehat y_4+x_1\otimes \widehat h_2\wedge \widehat x_3+x_4\otimes \widehat h_2\wedge \widehat x_5+y_3\otimes \widehat h_2\wedge \widehat y_1+y_5\otimes \widehat h_2\wedge \widehat y_4\\
 \hphantom{c_{0}^{2} =}{}
 +x_1\otimes \widehat x_5\wedge \widehat y_2+x_1\otimes \widehat x_6\wedge \widehat y_4+y_3\otimes \widehat x_2\wedge \widehat y_4+y_3\otimes \widehat x_5\wedge \widehat y_6,\\
c_{0}^{3} =
 h_2\otimes \widehat x_1\wedge \widehat y_3+x_3\otimes \widehat h_2\wedge \widehat x_1+x_3\otimes \widehat h_3\wedge \widehat x_1+x_5\otimes \widehat h_2\wedge \widehat x_4+x_5\otimes \widehat h_3\wedge \widehat x_4\\
 \hphantom{c_{0}^{3} =}{}
 +y_1\otimes \widehat h_2\wedge \widehat y_3+y_1\otimes \widehat h_3\wedge \widehat y_3+y_4\otimes \widehat h_2\wedge \widehat y_5
 +y_4\otimes \widehat h_3\wedge \widehat y_5+x_2\otimes \widehat x_4\wedge \widehat y_3\\
\hphantom{c_{0}^{3} =}{}
 +y_1\otimes \widehat x_2\wedge \widehat y_5+y_1\otimes \widehat x_4\wedge \widehat y_6+y_4\otimes \widehat x_1\wedge \widehat y_6+x_5\otimes \widehat x_1\wedge \widehat x_2+y_6\otimes \widehat y_3\wedge \widehat y_5,\\
c_{0}^{4} =
 h_2\otimes \widehat x_3\wedge \widehat y_1+x_1\otimes \widehat h_2\wedge \widehat x_3+x_1\otimes \widehat h_3\wedge \widehat x_3+x_4\otimes \widehat h_2\wedge \widehat x_5+x_4\otimes \widehat h_3\wedge \widehat x_5\\
\hphantom{c_{0}^{4} =}{}
 +y_3\otimes \widehat h_2\wedge \widehat y_1+y_3\otimes \widehat h_3\wedge \widehat y_1
 +y_5\otimes \widehat h_2\wedge \widehat y_4+y_5\otimes \widehat h_3\wedge \widehat y_4+x_2\otimes \widehat x_5\wedge \widehat y_1\\
\hphantom{c_{0}^{4} =}{}
 +y_3\otimes \widehat x_2\wedge \widehat y_4+y_3\otimes \widehat x_5\wedge \widehat y_6+y_5\otimes \widehat x_3\wedge \widehat y_6+x_4\otimes \widehat x_2\wedge \widehat x_3+y_6\otimes \widehat y_1\wedge \widehat y_4.
\end{gather*}

\item[--] $\mathfrak{gl}(2|2)$ and $\mathfrak{psl}(2|2)\simeq\mathfrak{h}_\Pi^{(1)}(0|4)$ (an exceptional case due to this occasional isomorphism; compare with quantization of $\mathfrak{po}(0|4)$, see \cite{LSh3} and Tyutin's result \cite{Tyut} over $\Cee$.)

\item[--] $\fg=\mathfrak{gl}(6)$: $H^2(\fg;\fg)=0$ and for $\fg:=\mathfrak{psl}(6)$ the space $H^2(\fg;\fg)$ is spanned, in agreement with Chebochko's Theorem~\ref{NC2}, by (observe thatindices of the cocycles are the degrees in the standard $\Zee$-grading)
\begin{gather*}
 c_{-9} =
 y_3\otimes\widehat x_{11}\wedge \widehat x_{15}+y_3\otimes\widehat x_{13}\wedge \widehat x_{14}+y_7\otimes\widehat x_8\wedge \widehat x_{15}+y_7\otimes\widehat x_{12}\wedge \widehat x_{13}\\
\hphantom{c_{-9} =}{}
+y_8\otimes\widehat x_7\wedge \widehat x_{15} +y_8\otimes\widehat x_{10}\wedge \widehat x_{14}+y_{10}\otimes\widehat x_8\wedge \widehat x_{14}
 +y_{10}\otimes\widehat x_{11}\wedge \widehat x_{12}\\
\hphantom{c_{-9} =}{}
+y_{11}\otimes\widehat x_3\wedge \widehat x_{15} +y_{11}\otimes\widehat x_{10}\wedge \widehat x_{12}+y_{12}\otimes\widehat x_7\wedge \widehat x_{13}+y_{12}\otimes\widehat x_{10}\wedge \widehat x_{11}\\
\hphantom{c_{-9} =}{}
+y_{13}\otimes\widehat x_3\wedge \widehat x_{14} +y_{13}\otimes\widehat x_7\wedge \widehat x_{12}+y_{14}\otimes\widehat x_3\wedge \widehat x_{13}+y_{14}\otimes\widehat x_8\wedge \widehat x_{10}\\
\hphantom{c_{-9} =}{}
+y_{15}\otimes\widehat x_3\wedge \widehat x_{11} +y_{15}\otimes\widehat x_7\wedge \widehat x_8,\\
 c_{-7} =
 y_2\otimes\widehat x_4\wedge \widehat x_{15}+y_2\otimes\widehat x_9\wedge \widehat x_{13}+y_4\otimes\widehat x_2\wedge \widehat x_{15}+y_4\otimes\widehat x_6\wedge \widehat x_{14}+y_6\otimes\widehat x_4\wedge \widehat x_{14}\\
 \hphantom{c_{-7} =}{}
 +y_6\otimes\widehat x_9\wedge \widehat x_{11}+y_9\otimes\widehat x_2\wedge \widehat x_{13}
 +y_9\otimes\widehat x_6\wedge \widehat x_{11}+y_{11}\otimes\widehat x_6\wedge \widehat x_9\\
\hphantom{c_{-7} =}{}
 +y_{11}\otimes\widehat x_{15}\wedge \widehat y_3+y_{13}\otimes\widehat x_2\wedge \widehat x_9+y_{13}\otimes\widehat x_{14}\wedge \widehat y_3+y_{14}\otimes\widehat x_4\wedge \widehat x_6\\
\hphantom{c_{-7} =}{}
 +y_{14}\otimes\widehat x_{13}\wedge \widehat y_3
 +y_{15}\otimes\widehat x_2\wedge \widehat x_4+y_{15}\otimes\widehat x_{11}\wedge \widehat y_3+x_3\otimes\widehat x_{11}\wedge \widehat x_{15}\\
\hphantom{c_{-7} =}{}
 +x_3\otimes\widehat x_{13}\wedge \widehat x_{14},\\
 c_{-5}^{1} =
 y_2\otimes\widehat x_5\wedge \widehat x_{10}+y_2\otimes\widehat x_{15}\wedge \widehat y_4+y_5\otimes\widehat x_2\wedge \widehat x_{10}+y_5\otimes\widehat x_6\wedge \widehat x_7+y_6\otimes\widehat x_5\wedge \widehat x_7\\
 \hphantom{c_{-5}^{1} =}{}
 +y_6\otimes\widehat x_{14}\wedge \widehat y_4+y_7\otimes\widehat x_5\wedge \widehat x_6
 +y_7\otimes\widehat x_{15}\wedge \widehat y_8+y_{10}\otimes\widehat x_2\wedge \widehat x_5\\
\hphantom{c_{-5}^{1} =}{}
 +y_{10}\otimes\widehat x_{14}\wedge \widehat y_8+y_{14}\otimes\widehat x_6\wedge \widehat y_4+y_{14}\otimes\widehat x_{10}\wedge \widehat y_8+y_{15}\otimes\widehat x_2\wedge \widehat y_4\\
\hphantom{c_{-5}^{1} =}{}
 +y_{15}\otimes\widehat x_7\wedge \widehat y_8
 +x_4\otimes\widehat x_2\wedge \widehat x_{15}+x_4\otimes\widehat x_6\wedge \widehat x_{14}+x_8\otimes\widehat x_7\wedge \widehat x_{15}\\
\hphantom{c_{-5}^{1} =}{}
 +x_8\otimes\widehat x_{10}\wedge \widehat x_{14},\\
 c_{-5}^{2} =
 y_1\otimes\widehat x_4\wedge \widehat x_{12}+y_1\otimes\widehat x_8\wedge \widehat x_9+y_4\otimes\widehat x_1\wedge \widehat x_{12}+y_4\otimes\widehat x_{15}\wedge \widehat y_2+y_8\otimes\widehat x_1\wedge \widehat x_9\\
 \hphantom{c_{-5}^{2} =}{}
 +y_8\otimes\widehat x_{15}\wedge \widehat y_7+y_9\otimes\widehat x_1\wedge \widehat x_8
 +y_9\otimes\widehat x_{13}\wedge \widehat y_2+y_{12}\otimes\widehat x_1\wedge \widehat x_4\\
\hphantom{c_{-5}^{2} =}{}
 +y_{12}\otimes\widehat x_{13}\wedge \widehat y_7+y_{13}\otimes\widehat x_9\wedge \widehat y_2+y_{13}\otimes\widehat x_{12}\wedge \widehat y_7+y_{15}\otimes\widehat x_4\wedge \widehat y_2\\
\hphantom{c_{-5}^{2} =}{}
 +y_{15}\otimes\widehat x_8\wedge \widehat y_7
 +x_2\otimes\widehat x_4\wedge \widehat x_{15}+x_2\otimes\widehat x_9\wedge \widehat x_{13}+x_7\otimes\widehat x_8\wedge \widehat x_{15}\\
\hphantom{c_{-5}^{2} =}{}
 +x_7\otimes\widehat x_{12}\wedge \widehat x_{13},\\
 c_{-3}^{1} = y_2\otimes\widehat x_{10}\wedge \widehat y_5+y_2\otimes\widehat x_{13}\wedge \widehat y_9+y_6\otimes\widehat x_7\wedge \widehat y_5+y_6\otimes\widehat x_{11}\wedge \widehat y_9+y_7\otimes\widehat x_6\wedge \widehat y_5\\
 \hphantom{c_{-3}^{1} =}{}
 +y_7\otimes\widehat x_{13}\wedge \widehat y_{12}+y_{10}\otimes\widehat x_2\wedge \widehat y_5
+y_{10}\otimes\widehat x_{11}\wedge \widehat y_{12}+y_{11}\otimes\widehat x_6\wedge \widehat y_9\\
\hphantom{c_{-3}^{1} =}{}
+y_{11}\otimes\widehat x_{10}\wedge \widehat y_{12}+y_{13}\otimes\widehat x_2\wedge \widehat y_9+y_{13}\otimes\widehat x_7\wedge \widehat y_{12}+x_5\otimes\widehat x_2\wedge \widehat x_{10}\\
\hphantom{c_{-3}^{1} =}{}
+x_5\otimes\widehat x_6\wedge \widehat x_7
+x_9\otimes\widehat x_2\wedge \widehat x_{13}+x_9\otimes\widehat x_6\wedge \widehat x_{11}+x_{12}\otimes\widehat x_7\wedge \widehat x_{13}\\
\hphantom{c_{-3}^{1} =}{}
+x_{12}\otimes\widehat x_{10}\wedge \widehat x_{11},\\
 c_{-3}^{2} = y_4\otimes\widehat x_{12}\wedge \widehat y_1+y_4\otimes\widehat x_{14}\wedge \widehat y_6+y_8\otimes\widehat x_9\wedge \widehat y_1+y_8\otimes\widehat x_{14}\wedge \widehat y_{10}+y_9\otimes\widehat x_8\wedge \widehat y_1\\
\hphantom{c_{-3}^{2} =}{}
 +y_9\otimes\widehat x_{11}\wedge \widehat y_6+y_{11}\otimes\widehat x_9\wedge \widehat y_6
+y_{11}\otimes\widehat x_{12}\wedge \widehat y_{10}+y_{12}\otimes\widehat x_4\wedge \widehat y_1\\
\hphantom{c_{-3}^{2} =}{}
+y_{12}\otimes\widehat x_{11}\wedge \widehat y_{10}+y_{14}\otimes\widehat x_4\wedge \widehat y_6+y_{14}\otimes\widehat x_8\wedge \widehat y_{10}+x_1\otimes\widehat x_4\wedge \widehat x_{12}\\
\hphantom{c_{-3}^{2} =}{}
+x_1\otimes\widehat x_8\wedge \widehat x_9
+x_6\otimes\widehat x_4\wedge \widehat x_{14}+x_6\otimes\widehat x_9\wedge \widehat x_{11}+x_{10}\otimes\widehat x_8\wedge \widehat x_{14}\\
\hphantom{c_{-3}^{2} =}{}
+x_{10}\otimes\widehat x_{11}\wedge \widehat x_{12},\\
 c_{-3}^{3} = x_2\otimes\widehat x_{15}\wedge \widehat y_4+x_4\otimes\widehat x_{15}\wedge \widehat y_2+y_1\otimes\widehat x_3\wedge \widehat x_5+y_1\otimes\widehat x_{12}\wedge \widehat y_4+y_3\otimes\widehat x_1\wedge \widehat x_5\\
\hphantom{c_{-3}^{3} =}{}
 +y_3\otimes\widehat x_{15}\wedge \widehat y_{11}+y_5\otimes\widehat x_1\wedge \widehat x_3
+y_5\otimes\widehat x_{10}\wedge \widehat y_2+y_{10}\otimes\widehat x_5\wedge \widehat y_2\\
\hphantom{c_{-3}^{3} =}{}
+y_{10}\otimes\widehat x_{12}\wedge \widehat y_{11}+y_{12}\otimes\widehat x_1\wedge \widehat y_4+y_{12}\otimes\widehat x_{10}\wedge \widehat y_{11}+y_{15}\otimes\widehat x_3\wedge \widehat y_{11}\\
\hphantom{c_{-3}^{3} =}{}
+x_2\otimes\widehat x_5\wedge \widehat x_{10}
+x_4\otimes\widehat x_1\wedge \widehat x_{12}+x_{11}\otimes\widehat x_3\wedge \widehat x_{15}+x_{11}\otimes\widehat x_{10}\wedge \widehat x_{12}\\
\hphantom{c_{-3}^{3} =}{}
+y_{15}\otimes\widehat y_2\wedge \widehat y_4,\\
 c_{-1}^{1} = x_2\otimes\widehat x_{10}\wedge \widehat y_5+x_2\otimes\widehat x_{13}\wedge \widehat y_9+x_5\otimes\widehat x_{10}\wedge \widehat y_2+x_9\otimes\widehat x_{13}\wedge \widehat y_2+y_1\otimes\widehat x_3\wedge \widehat y_5\\
\hphantom{c_{-1}^{1} =}{}
 +y_1\otimes\widehat x_8\wedge \widehat y_9+y_3\otimes\widehat x_1\wedge \widehat y_5
+y_3\otimes\widehat x_{13}\wedge \widehat y_{14}+y_8\otimes\widehat x_1\wedge \widehat y_9+y_8\otimes\widehat x_{10}\wedge \widehat y_{14}\\
\hphantom{c_{-1}^{1} =}{}
+y_{10}\otimes\widehat x_8\wedge \widehat y_{14}+y_{13}\otimes\widehat x_3\wedge \widehat y_{14}+x_5\otimes\widehat x_1\wedge \widehat x_3+x_9\otimes\widehat x_1\wedge \widehat x_8
\\
\hphantom{c_{-1}^{1} =}{}
+x_{14}\otimes\widehat x_3\wedge \widehat x_{13}+x_{14}\otimes\widehat x_8\wedge \widehat x_{10}+y_{10}\otimes\widehat y_2\wedge \widehat y_5+y_{13}\otimes\widehat y_2\wedge \widehat y_9,\\
 c_{-1}^{2} =
 x_1\otimes\widehat x_{12}\wedge \widehat y_4+x_4\otimes\widehat x_{12}\wedge \widehat y_1+x_4\otimes\widehat x_{14}\wedge \widehat y_6+x_6\otimes\widehat x_{14}\wedge \widehat y_4+y_3\otimes\widehat x_5\wedge \widehat y_1\\
 \hphantom{c_{-1}^{2} =}{}
 +y_3\otimes\widehat x_{14}\wedge \widehat y_{13}+y_5\otimes\widehat x_3\wedge \widehat y_1
+y_5\otimes\widehat x_7\wedge \widehat y_6+y_7\otimes\widehat x_5\wedge \widehat y_6+y_7\otimes\widehat x_{12}\wedge \widehat y_{13}\\
 \hphantom{c_{-1}^{2} =}{}
+y_{12}\otimes\widehat x_7\wedge \widehat y_{13}+y_{14}\otimes\widehat x_3\wedge \widehat y_{13}+x_1\otimes\widehat x_3\wedge \widehat x_5+x_6\otimes\widehat x_5\wedge \widehat x_7\\
\hphantom{c_{-1}^{2} =}{}
+x_{13}\otimes\widehat x_3\wedge \widehat x_{14}+x_{13}\otimes\widehat x_7\wedge \widehat x_{12}+y_{12}\otimes\widehat y_1\wedge \widehat y_4+y_{14}\otimes\widehat y_4\wedge \widehat y_6,\\
 c_{-1}^{3} =
 x_3\otimes\widehat x_{15}\wedge \widehat y_{11}+x_7\otimes\widehat x_{15}\wedge \widehat y_8+x_8\otimes\widehat x_{15}\wedge \widehat y_7+x_{11}\otimes\widehat x_{15}\wedge \widehat y_3+y_1\otimes\widehat x_5\wedge \widehat y_3\\
\hphantom{c_{-1}^{3} =}{}
 +y_1\otimes\widehat x_9\wedge \widehat y_8+y_5\otimes\widehat x_1\wedge \widehat y_3
+y_5\otimes\widehat x_6\wedge \widehat y_7+y_6\otimes\widehat x_5\wedge \widehat y_7+y_6\otimes\widehat x_9\wedge \widehat y_{11}\\
\hphantom{c_{-1}^{3} =}{}
+y_9\otimes\widehat x_1\wedge \widehat y_8+y_9\otimes\widehat x_6\wedge \widehat y_{11}+x_3\otimes\widehat x_1\wedge \widehat x_5+x_7\otimes\widehat x_5\wedge \widehat x_6
+x_8\otimes\widehat x_1\wedge \widehat x_9\\
\hphantom{c_{-1}^{3} =}{}
+x_{11}\otimes\widehat x_6\wedge \widehat x_9+y_{15}\otimes\widehat y_3\wedge \widehat y_{11}+y_{15}\otimes\widehat y_7\wedge \widehat y_8.
\end{gather*}

\item[--] $\mathfrak{gl}(8)$ and $\mathfrak{psl}(8)$: in both cases $H^2(\fg;\fg)=0$, see Chebochko's Theorem~\ref{NC2}.
\end{itemize}

For $\fg=\mathfrak{gl}(1|3)$ in Chevalley basis given by formula~\eqref{gl4cm}, the space $H^2(\fg;\fg)$ is spanned by
\begin{gather*}
 c_{-6} =
 h_2\otimes (\widehat x_6)^{\wedge2}+y_1\otimes \widehat x_5\wedge \widehat x_6+y_3\otimes \widehat x_4\wedge \widehat x_6+y_4\otimes \widehat x_3\wedge \widehat x_6+y_5\otimes \widehat x_1\wedge \widehat x_6\\
 \hphantom{c_{-6} =}{}
 +y_6\otimes \widehat x_1\wedge \widehat x_5+y_6\otimes \widehat x_3\wedge \widehat x_4,\\
 c_{-4} =
 h_1\otimes (\widehat x_4)^{\wedge2}+h_3\otimes (\widehat x_4)^{\wedge2},\\
 c_{-2} =
 h_1\otimes (\widehat x_1)^{\wedge2}+h_3\otimes (\widehat x_1)^{\wedge2}.
\end{gather*}

For $\fg=\mathfrak{psl}(1|3)$, we have
$H^2(\fg;\fg) =0$.
For $\fg=\mathfrak{gl}(1|5)$, the space $H^2(\fg;\fg)$ is spanned by
\begin{gather*}
c_{-10} =
 h_1\otimes\widehat x_{15}{}^{\wedge 2}+ h_3\otimes\widehat x_{15}{}^{\wedge 2}+h_5\otimes\widehat x_{15}{}^{\wedge 2},\\
c_{-8} =
 h_1\otimes\widehat x_{13}{}^{\wedge 2}+ h_3\otimes\widehat x_{13}{}^{\wedge 2}+h_5\otimes\widehat x_{13}{}^{\wedge 2},\\
c_{-6} =
 h_1\otimes\widehat x_{10}{}^{\wedge 2}+h_3\otimes\widehat x_{10}{}^{\wedge 2}+h_5\otimes\widehat x_{10}{}^{\wedge 2},\\
c_{-4} =
 h_1\otimes\widehat x_6{}^{\wedge 2}+h_3\otimes\widehat x_6{}^{\wedge 2}+h_5\otimes\widehat x_6{}^{\wedge 2},\\
c_{-2} =
 h_1\otimes\widehat x_1{}^{\wedge 2}+h_3\otimes\widehat x_1{}^{\wedge 2}+h_5\otimes\widehat x_1{}^{\wedge 2}.
\end{gather*}

For $\fg=\mathfrak{psl}(1|5)$, the space $H^2(\fg;\fg)$ is spanned by (observe that
indices of the cocycles are the degrees in the standard $\Zee$-grading)
\begin{gather*}
c_{-9} = y_3\otimes\widehat x_{11}\wedge \widehat x_{15}+y_3\otimes\widehat x_{13}\wedge \widehat x_{14}+y_7\otimes\widehat x_8\wedge \widehat x_{15}+y_7\otimes\widehat x_{12}\wedge \widehat x_{13}+y_8\otimes\widehat x_7\wedge \widehat x_{15}\\
\hphantom{c_{-9} =}{}
+y_8\otimes\widehat x_{10}\wedge \widehat x_{14}
 +y_{10}\otimes\widehat x_8\wedge \widehat x_{14}+y_{10}\otimes\widehat x_{11}\wedge \widehat x_{12}+y_{11}\otimes\widehat x_3\wedge \widehat x_{15}\\
\hphantom{c_{-9} =}{}
 +y_{11}\otimes\widehat x_{10}\wedge \widehat x_{12}+y_{12}\otimes\widehat x_7\wedge \widehat x_{13}
 +y_{12}\otimes\widehat x_{10}\wedge \widehat x_{11}+y_{13}\otimes\widehat x_3\wedge \widehat x_{14}\\
\hphantom{c_{-9} =}{}
 +y_{13}\otimes\widehat x_7\wedge \widehat x_{12}+y_{14}\otimes\widehat x_3\wedge \widehat x_{13}
 +y_{14}\otimes\widehat x_8\wedge \widehat x_{10}+y_{15}\otimes\widehat x_3\wedge \widehat x_{11}+y_{15}\otimes\widehat x_7\wedge \widehat x_8,\\
c_{-7} =
 y_2\otimes\widehat x_4\wedge \widehat x_{15}+y_2\otimes\widehat x_9\wedge \widehat x_{13}+y_4\otimes\widehat x_2\wedge \widehat x_{15}+y_4\otimes\widehat x_6\wedge \widehat x_{14}+y_6\otimes\widehat x_4\wedge \widehat x_{14}\\
\hphantom{c_{-7} =}{}
 +y_6\otimes\widehat x_9\wedge \widehat x_{11}
 +y_9\otimes\widehat x_2\wedge \widehat x_{13}+y_9\otimes\widehat x_6\wedge \widehat x_{11}+y_{11}\otimes\widehat x_6\wedge \widehat x_9+y_{11}\otimes\widehat x_{15}\wedge \widehat y_3\\
\hphantom{c_{-7} =}{}
 +y_{13}\otimes\widehat x_2\wedge \widehat x_9+y_{13}\otimes\widehat x_{14}\wedge \widehat y_3
 +y_{14}\otimes\widehat x_4\wedge \widehat x_6+y_{14}\otimes\widehat x_{13}\wedge \widehat y_3+y_{15}\otimes\widehat x_2\wedge \widehat x_4\\
\hphantom{c_{-7} =}{}
 +y_{15}\otimes\widehat x_{11}\wedge \widehat y_3+x_3\otimes\widehat x_{11}\wedge \widehat x_{15}+x_3\otimes\widehat x_{13}\wedge \widehat x_{14},\\
c_{-5}^1 =
 y_2\otimes\widehat x_5\wedge \widehat x_{10}+y_2\otimes\widehat x_{15}\wedge \widehat y_4+y_5\otimes\widehat x_2\wedge \widehat x_{10}+y_5\otimes\widehat x_6\wedge \widehat x_7+y_6\otimes\widehat x_5\wedge \widehat x_7\\
 \hphantom{c_{-5}^1 =}{}
 +y_6\otimes\widehat x_{14}\wedge \widehat y_4
 +y_7\otimes\widehat x_5\wedge \widehat x_6+y_7\otimes\widehat x_{15}\wedge \widehat y_8+y_{10}\otimes\widehat x_2\wedge \widehat x_5+y_{10}\otimes\widehat x_{14}\wedge \widehat y_8\\
 \hphantom{c_{-5}^1 =}{}
 +y_{14}\otimes\widehat x_6\wedge \widehat y_4+y_{14}\otimes\widehat x_{10}\wedge \widehat y_8
 +y_{15}\otimes\widehat x_2\wedge \widehat y_4+y_{15}\otimes\widehat x_7\wedge \widehat y_8+x_4\otimes\widehat x_2\wedge \widehat x_{15}\\
 \hphantom{c_{-5}^1 =}{}
 +x_4\otimes\widehat x_6\wedge \widehat x_{14}+x_8\otimes\widehat x_7\wedge \widehat x_{15}+x_8\otimes\widehat x_{10}\wedge \widehat x_{14},\\
c_{-5}^2 =
 y_1\otimes\widehat x_4\wedge \widehat x_{12}+y_1\otimes\widehat x_8\wedge \widehat x_9+y_4\otimes\widehat x_1\wedge \widehat x_{12}+y_4\otimes\widehat x_{15}\wedge \widehat y_2+y_8\otimes\widehat x_1\wedge \widehat x_9\\
\hphantom{c_{-5}^2 =}{}
 +y_8\otimes\widehat x_{15}\wedge \widehat y_7
 +y_9\otimes\widehat x_1\wedge \widehat x_8+y_9\otimes\widehat x_{13}\wedge \widehat y_2+y_{12}\otimes\widehat x_1\wedge \widehat x_4+y_{12}\otimes\widehat x_{13}\wedge \widehat y_7\\
\hphantom{c_{-5}^2 =}{}
 +y_{13}\otimes\widehat x_9\wedge \widehat y_2+y_{13}\otimes\widehat x_{12}\wedge \widehat y_7
 +y_{15}\otimes\widehat x_4\wedge \widehat y_2+y_{15}\otimes\widehat x_8\wedge \widehat y_7+x_2\otimes\widehat x_4\wedge \widehat x_{15}\\
\hphantom{c_{-5}^2 =}{}
 +x_2\otimes\widehat x_9\wedge \widehat x_{13}+x_7\otimes\widehat x_8\wedge \widehat x_{15}+x_7\otimes\widehat x_{12}\wedge \widehat x_{13},\\
c_{-3}^1 =
 y_2\otimes\widehat x_{10}\wedge \widehat y_5+y_2\otimes\widehat x_{13}\wedge \widehat y_9+y_6\otimes\widehat x_7\wedge \widehat y_5+y_6\otimes\widehat x_{11}\wedge \widehat y_9+y_7\otimes\widehat x_6\wedge \widehat y_5\\
 \hphantom{c_{-3}^1 =}{}
 +y_7\otimes\widehat x_{13}\wedge \widehat y_{12}
 +y_{10}\otimes\widehat x_2\wedge \widehat y_5+y_{10}\otimes\widehat x_{11}\wedge \widehat y_{12}+y_{11}\otimes\widehat x_6\wedge \widehat y_9+y_{11}\otimes\widehat x_{10}\wedge \widehat y_{12}\\
\hphantom{c_{-3}^1 =}{}
+y_{13}\otimes\widehat x_2\wedge \widehat y_9 +y_{13}\otimes\widehat x_7\wedge \widehat y_{12}
 +x_5\otimes\widehat x_2\wedge \widehat x_{10}+x_5\otimes\widehat x_6\wedge \widehat x_7+x_9\otimes\widehat x_2\wedge \widehat x_{13}\\
\hphantom{c_{-3}^1 =}{}
+x_9\otimes\widehat x_6\wedge \widehat x_{11} +x_{12}\otimes\widehat x_7\wedge \widehat x_{13}+x_{12}\otimes\widehat x_{10}\wedge \widehat x_{11},\\
c_{-3}^2 =
 y_4\otimes\widehat x_{12}\wedge \widehat y_1+y_4\otimes\widehat x_{14}\wedge \widehat y_6+y_8\otimes\widehat x_9\wedge \widehat y_1+y_8\otimes\widehat x_{14}\wedge \widehat y_{10}+y_9\otimes\widehat x_8\wedge \widehat y_1\\
\hphantom{c_{-3}^2 =}{}
 +y_9\otimes\widehat x_{11}\wedge \widehat y_6
 +y_{11}\otimes\widehat x_9\wedge \widehat y_6+y_{11}\otimes\widehat x_{12}\wedge \widehat y_{10}+y_{12}\otimes\widehat x_4\wedge \widehat y_1+y_{12}\otimes\widehat x_{11}\wedge \widehat y_{10}\\
\hphantom{c_{-3}^2 =}{}
 +y_{14}\otimes\widehat x_4\wedge \widehat y_6+y_{14}\otimes\widehat x_8\wedge \widehat y_{10}
 +x_1\otimes\widehat x_4\wedge \widehat x_{12}+x_1\otimes\widehat x_8\wedge \widehat x_9+x_6\otimes\widehat x_4\wedge \widehat x_{14}\\
\hphantom{c_{-3}^2 =}{}
 +x_6\otimes\widehat x_9\wedge \widehat x_{11}+x_{10}\otimes\widehat x_8\wedge \widehat x_{14}+x_{10}\otimes\widehat x_{11}\wedge \widehat x_{12},\\
c_{-3}^3 =
 x_2\otimes\widehat x_{15}\wedge \widehat y_4+x_4\otimes\widehat x_{15}\wedge \widehat y_2+y_1\otimes\widehat x_3\wedge \widehat x_5+y_1\otimes\widehat x_{12}\wedge \widehat y_4+y_3\otimes\widehat x_1\wedge \widehat x_5\\
\hphantom{c_{-3}^3 =}{}
 +y_3\otimes\widehat x_{15}\wedge \widehat y_{11}
 +y_5\otimes\widehat x_1\wedge \widehat x_3+y_5\otimes\widehat x_{10}\wedge \widehat y_2+y_{10}\otimes\widehat x_5\wedge \widehat y_2+y_{10}\otimes\widehat x_{12}\wedge \widehat y_{11}\\
\hphantom{c_{-3}^3 =}{}
 +y_{12}\otimes\widehat x_1\wedge \widehat y_4+y_{12}\otimes\widehat x_{10}\wedge \widehat y_{11}
 +y_{15}\otimes\widehat x_3\wedge \widehat y_{11}+x_2\otimes\widehat x_5\wedge \widehat x_{10}+x_4\otimes\widehat x_1\wedge \widehat x_{12}\\
\hphantom{c_{-3}^3 =}{}
 +x_{11}\otimes\widehat x_3\wedge \widehat x_{15}+x_{11}\otimes\widehat x_{10}\wedge \widehat x_{12}+y_{15}\otimes\widehat y_2\wedge \widehat y_4,\\
c_{-1}^1 =
 x_2\otimes\widehat x_{10}\wedge \widehat y_5+x_2\otimes\widehat x_{13}\wedge \widehat y_9+x_5\otimes\widehat x_{10}\wedge \widehat y_2+x_9\otimes\widehat x_{13}\wedge \widehat y_2+y_1\otimes\widehat x_3\wedge \widehat y_5\\
\hphantom{c_{-1}^1 =}{}
 +y_1\otimes\widehat x_8\wedge \widehat y_9
 +y_3\otimes\widehat x_1\wedge \widehat y_5+y_3\otimes\widehat x_{13}\wedge \widehat y_{14}+y_8\otimes\widehat x_1\wedge \widehat y_9+y_8\otimes\widehat x_{10}\wedge \widehat y_{14}\\
\hphantom{c_{-1}^1 =}{}
 +y_{10}\otimes\widehat x_8\wedge \widehat y_{14}+y_{13}\otimes\widehat x_3\wedge \widehat y_{14}
 +x_5\otimes\widehat x_1\wedge \widehat x_3+x_9\otimes\widehat x_1\wedge \widehat x_8+x_{14}\otimes\widehat x_3\wedge \widehat x_{13}\\
\hphantom{c_{-1}^1 =}{}
 +x_{14}\otimes\widehat x_8\wedge \widehat x_{10}+y_{10}\otimes\widehat y_2\wedge \widehat y_5+y_{13}\otimes\widehat y_2\wedge \widehat y_9,\\
c_{-1}^2 =
 x_1\otimes\widehat x_{12}\wedge \widehat y_4+x_4\otimes\widehat x_{12}\wedge \widehat y_1+x_4\otimes\widehat x_{14}\wedge \widehat y_6+x_6\otimes\widehat x_{14}\wedge \widehat y_4+y_3\otimes\widehat x_5\wedge \widehat y_1\\
 \hphantom{c_{-1}^2 =}{}
 +y_3\otimes\widehat x_{14}\wedge \widehat y_{13}
 +y_5\otimes\widehat x_3\wedge \widehat y_1+y_5\otimes\widehat x_7\wedge \widehat y_6+y_7\otimes\widehat x_5\wedge \widehat y_6+y_7\otimes\widehat x_{12}\wedge \widehat y_{13}\\
 \hphantom{c_{-1}^2 =}{}
 +y_{12}\otimes\widehat x_7\wedge \widehat y_{13}+y_{14}\otimes\widehat x_3\wedge \widehat y_{13}
 +x_1\otimes\widehat x_3\wedge \widehat x_5+x_6\otimes\widehat x_5\wedge \widehat x_7+x_{13}\otimes\widehat x_3\wedge \widehat x_{14}\\
 \hphantom{c_{-1}^2 =}{}
 +x_{13}\otimes\widehat x_7\wedge \widehat x_{12}+y_{12}\otimes\widehat y_1\wedge \widehat y_4+y_{14}\otimes\widehat y_4\wedge \widehat y_6,\\
c_{-1}^3 =
 x_3\otimes\widehat x_{15}\wedge \widehat y_{11}+x_7\otimes\widehat x_{15}\wedge \widehat y_8+x_8\otimes\widehat x_{15}\wedge \widehat y_7+x_{11}\otimes\widehat x_{15}\wedge \widehat y_3+y_1\otimes\widehat x_5\wedge \widehat y_3\\
 \hphantom{c_{-1}^3 =}{}
 +y_1\otimes\widehat x_9\wedge \widehat y_8
 +y_5\otimes\widehat x_1\wedge \widehat y_3+y_5\otimes\widehat x_6\wedge \widehat y_7+y_6\otimes\widehat x_5\wedge \widehat y_7+y_6\otimes\widehat x_9\wedge \widehat y_{11}\\
 \hphantom{c_{-1}^3 =}{}
 +y_9\otimes\widehat x_1\wedge \widehat y_8+y_9\otimes\widehat x_6\wedge \widehat y_{11}
 +x_3\otimes\widehat x_1\wedge \widehat x_5+x_7\otimes\widehat x_5\wedge \widehat x_6+x_8\otimes\widehat x_1\wedge \widehat x_9\\
 \hphantom{c_{-1}^3 =}{}
 +x_{11}\otimes\widehat x_6\wedge \widehat x_9+y_{15}\otimes\widehat y_3\wedge \widehat y_{11}+y_{15}\otimes\widehat y_7\wedge \widehat y_8.
\end{gather*}

\subsection{Comparison with Shen's ``variations"}\label{L5.9.1} The Lie algebra\footnote{In
\cite{She1}, Shen mentioned that $\mathfrak{psl}(4):=\fg(2)$ in the basis
given in \cite[p.~346]{FH}. Hence, the property of the root systems
to have all roots of ``equal length'' is meaningless if $p=2$, cf.~\cite{Ch1}.} $\fg=\mathfrak{psl}(4)$ has no Cartan matrix; its relative that
has a~Cartan matrix given by~\eqref{gl4cm} is $\mathfrak{gl}(4)$.
In order to compare with Shen's ``variations of
$\fg(2)$ theme'' we give, in parentheses, together with the
generators of~$\mathfrak{sl}(4)$, their weights in terms of the root systems
of~$\mathfrak{gl}(4)$ (we designate them~$\beta$s) and~$\fg(2)$ (we designate
them~$\alpha$s):
\begin{gather}
\begin{pmatrix}
\ev &1 &0\\
1 &\ev &1\\
0&1 &\ev
\end{pmatrix},
\nonumber\\
 x_1=E_{12}\
(\beta_1=\alpha_2),\qquad x_2=E_{23} \ (\beta_2=\alpha_1),
\qquad x_3=E_{34}\ (\beta_3=2\alpha_1+\alpha_2),\nonumber\\
x_4:=[x_1,x_2]=E_{13}\ (\beta_1+\beta_2=\alpha_1+\alpha_2),\qquad
x_5:=[x_2,x_3]=E_{24}\ (\beta_2+\beta_3=3\alpha_1+\alpha_2),\nonumber\\
{}x_6:=[x_3,[x_1,x_2]]=E_{14}\
(\beta_1+\beta_2+\beta_3=3\alpha_1+2\alpha_2).\label{gl4cm}
\end{gather}
Then, $H^2(\fg; \fg)$ is spanned by the cocycles \eqref{eqpsl4}. The
parameters $(a)$, $(b)$ in the first column correspond to Shen's
parameters $(a,0)$, $(0,b)$ of his $V_4G(2; a, b)$. The index of the cocycle is its degree under the standard grading:
\begin{alignat}{4}
&c_{-6}^1\quad &&-4(\alpha_1+\alpha_2)\quad &&y_1\otimes \widehat x_4\wedge \widehat
x_6+y_4\otimes \widehat x_1\wedge \widehat x_6+y_6\otimes \widehat x_1\wedge
\widehat x_4&\nonumber\\
&c_{-6}^2\quad &&-4(2\alpha_1+\alpha_2)\quad &&y_3\otimes \widehat x_5\wedge \widehat
x_6+y_5\otimes \widehat x_3\wedge \widehat x_6+y_6\otimes \widehat x_3\wedge
\widehat x_5& \nonumber\\
&c_{-4}^1 \ (b)\quad &&-2(2\alpha_1+\alpha_2)\quad &&h_2\otimes \widehat x_4\wedge \widehat
x_5+y_2\otimes \widehat h_3\wedge \widehat x_6+y_2\otimes \widehat x_1\wedge
x_5+y_2\otimes \widehat x_3\wedge \widehat x_4& \nonumber\\
&\quad &&\quad &&+y_4\otimes \widehat h_3\wedge \widehat
x_5+y_5\otimes \widehat h_3\wedge \widehat x_4+y_6\otimes \widehat h_3\wedge \widehat
x_2+y_6\otimes \widehat x_4\wedge \widehat y_1& \nonumber\\
&\quad &&\quad && +y_6\otimes \widehat x_5\wedge \widehat y_3& \nonumber\\
&c_{-4}^2\quad &&-3(\alpha_1+\alpha_2)\quad &&h_3\otimes \widehat x_2\wedge \widehat
x_6+x_1\otimes \widehat x_4\wedge \widehat x_6+x_3\otimes \widehat x_5\wedge \widehat
x_6+y_2\otimes \widehat h_2\wedge \widehat
x_6& \nonumber\\
&\quad &&\quad &&+y_2\otimes \widehat h_3\wedge \widehat x_6+y_2\otimes \widehat x_1\wedge \widehat
x_5+y_2\otimes \widehat x_3\wedge \widehat x_4+y_4\otimes \widehat h_2\wedge \widehat
x_5& \nonumber\\
&\quad &&\quad &&+y_4\otimes \widehat h_3\wedge \widehat x_5+y_4\otimes \widehat x_2\wedge \widehat
x_3+y_5\otimes \widehat h_2\wedge \widehat x_4+y_5\otimes \widehat h_3\wedge \widehat
x_4& \nonumber\\
&\quad &&\quad &&+y_5\otimes \widehat x_1\wedge \widehat x_2+y_6\otimes \widehat h_2\wedge \widehat
x_2+y_6\otimes \widehat h_3\wedge \widehat x_2 & \nonumber\\
&c_{-2}^1\quad &&0 \quad &&x_3\otimes \widehat x_2\wedge \widehat x_4 +y_2\otimes \widehat x_4\wedge
\widehat y_3 +y_4\otimes \widehat x_2\wedge \widehat y_3& \nonumber\\
&c_{-2}^2\quad &&-4\alpha_1 \quad &&x_1\otimes \widehat x_2\wedge \widehat x_5 +y_2\otimes \widehat
x_5\wedge \widehat y_1 +y_5\otimes \widehat x_2\wedge
\widehat y_1& \nonumber\\
&c_{-2}^3\quad &&-2(\alpha_1+\alpha_2)\quad &&h_3\otimes \widehat x_6\wedge \widehat y_2
+x_2\otimes \widehat h_2\wedge \widehat x_6 +y_1\otimes \widehat h_2\wedge \widehat x_3
+y_1\otimes \widehat x_5\wedge \widehat y_2& \nonumber\\
&\quad &&\quad &&+y_1\otimes \widehat x_6\wedge \widehat y_4 +y_3\otimes \widehat h_2\wedge \widehat x_1 +y_3\otimes \widehat x_4\wedge \widehat y_2
+y_3\otimes \widehat x_6\wedge \widehat y_5 & \nonumber\\
&\quad &&\quad && +y_6\otimes \widehat h_2\wedge \widehat y_2& \nonumber\\
&c_{-2}^4\quad &&-2(\alpha_1+\alpha_2)\quad &&h_2\otimes \widehat x_1\wedge \widehat x_3
+x_2\otimes \widehat h_3\wedge \widehat x_6 +x_4\otimes \widehat x_1\wedge \widehat x_6
+x_5\otimes \widehat x_3\wedge \widehat x_6 & \nonumber\\
&\quad &&\quad &&+y_1\otimes \widehat h_3\wedge
\widehat x_3+y_1\otimes \widehat x_5\wedge \widehat y_2 +y_1\otimes \widehat x_6\wedge \widehat y_4
+y_3\otimes \widehat h_3\wedge
\widehat x_1 & \nonumber\\
&\quad &&\quad && +y_3\otimes \widehat x_4\wedge \widehat y_2 +y_3\otimes \widehat x_6\wedge \widehat y_5 +
y_4\otimes \widehat x_3\wedge \widehat y_2 +y_5\otimes \widehat x_1\wedge \widehat y_2& \nonumber\\
&\quad &&\quad && +y_6\otimes \widehat h_3\wedge
\widehat y_2& \nonumber\\
&c_{0}^1\quad &&2\alpha_1\quad &&h_3\otimes \widehat x_4\wedge \widehat y_5 +x_3\otimes \widehat
h_2\wedge \widehat x_1 +x_3\otimes \widehat x_4\wedge \widehat y_2 +x_3\otimes \widehat
x_6\wedge \widehat y_5& \nonumber\\
&\quad &&\quad &&+ x_5\otimes \widehat h_2\wedge
\widehat x_4+y_1\otimes \widehat h_2\wedge \widehat y_3 +y_1\otimes \widehat x_2\wedge \widehat y_5
+y_1\otimes \widehat x_4\wedge \widehat y_6 & \nonumber\\
&\quad &&\quad && +y_4\otimes \widehat h_2\wedge
\widehat y_5& \nonumber\\
&c_{0}^2\quad &&-2\alpha_1\quad &&h_3\otimes \widehat x_5\wedge \widehat y_4 +x_1\otimes \widehat
h_2\wedge \widehat x_3 +x_1\otimes \widehat x_5\wedge \widehat y_2 +x_1\otimes \widehat
x_6\wedge \widehat y_4& \nonumber\\
&\quad &&\quad &&+x_4\otimes \widehat h_2\wedge
\widehat x_5+y_3\otimes \widehat h_2\wedge \widehat y_1 +y_3\otimes \widehat x_2\wedge \widehat y_4
+y_3\otimes \widehat x_5\wedge \widehat y_6 & \nonumber\\
&\quad &&\quad && +y_5\otimes \widehat h_2\wedge
\widehat y_4& \nonumber\\
&c_{0}^3\ (a)\quad &&2\alpha_1\quad &&h_2\otimes \widehat x_1\wedge \widehat y_3 +x_2\otimes \widehat
x_4\wedge \widehat y_3 +x_3\otimes \widehat h_2\wedge \widehat x_1 +x_3\otimes \widehat
h_3\wedge \widehat x_1& \nonumber\\
&\quad &&\quad &&+x_5\otimes \widehat h_2\wedge
\widehat x_4+x_5\otimes \widehat h_3\wedge \widehat x_4 +x_5\otimes \widehat x_1\wedge \widehat x_2
+y_1\otimes \widehat h_2\wedge
\widehat y_3 & \nonumber\\
&\quad &&\quad && +y_1\otimes \widehat h_3\wedge \widehat y_3 +y_1\otimes \widehat x_2\wedge \widehat y_5 +
y_1\otimes \widehat x_4\wedge \widehat y_6 +y_4\otimes \widehat h_2\wedge \widehat y_5& \nonumber\\
&\quad &&\quad &&
+y_4\otimes \widehat h_3\wedge \widehat y_5 +y_4\otimes \widehat x_1\wedge \widehat y_6
+y_6\otimes \widehat y_3\wedge
\widehat y_5& \nonumber\\
&c_{0}^4\quad &&-2\alpha_1\quad &&h_2\otimes \widehat x_3\wedge \widehat y_1 +x_1\otimes \widehat
h_2\wedge \widehat x_3 +x_1\otimes \widehat h_3\wedge \widehat x_3 +x_2\otimes \widehat
x_5\wedge \widehat y_1& \nonumber\\
&\quad &&\quad &&+x_4\otimes \widehat h_2\wedge
\widehat x_5+x_4\otimes \widehat h_3\wedge \widehat x_5 +x_4\otimes \widehat x_2\wedge \widehat x_3
+y_3\otimes \widehat h_2\wedge
\widehat y_1 & \nonumber\\
&\quad &&\quad && +y_3\otimes \widehat h_3\wedge \widehat y_1 +y_3\otimes \widehat x_2\wedge \widehat y_4 +
y_3\otimes \widehat x_5\wedge \widehat y_6 +y_5\otimes \widehat h_2\wedge \widehat y_4& \nonumber\\
&\quad &&\quad &&
+y_5\otimes \widehat h_3\wedge \widehat y_4 +y_5\otimes \widehat x_3\wedge \widehat y_6
+y_6\otimes \widehat y_1\wedge \widehat y_4&
\label{eqpsl4}
\end{alignat}

%\subsection{On isomorphisms}\label{ssOnIso} See footnote to Open problems \ref{OP}.

\subsection[$\mathfrak{gs}(2)$: Shen's analog of $\fg(2)$ for $p=2$]{$\boldsymbol{\mathfrak{gs}(2)}$: Shen's analog of $\boldsymbol{\fg(2)}$ for $\boldsymbol{p=2}$}\label{ssShen} For a~description of this
Lie algebra, which we call $\mathfrak{gs}(2)$ in honor of Shen Guangyu, see \cite{BGLL1}. Let $\fg=\mathfrak{gs}(2)$.

The space of
$H^2(\fg;\fg)$ is spanned by the following cocycles:
\begin{gather*}%\label{deforShen}
c_{-8}= \partial_1\otimes (\widehat{u_1u_2})\wedge
(\widehat{u_1u_2\partial_1})+\partial_2\otimes (\widehat{u_1u_2})\wedge
(\widehat{u_1u_2\partial_2})+v\otimes (\widehat{u_1u_2\partial_1})\wedge (\widehat{u_1u_2\partial_2}),\\
c_{-4}= \partial_1\otimes (\widehat{v})\wedge (\widehat{u_1u_2\partial_1})
+\partial_1\otimes (\widehat{u_1\partial_1})\wedge (\widehat{u_2v}) +\partial_1\otimes
(\widehat{u_2\partial_2})\wedge (\widehat{u_2v})\\
\hphantom{c_{-4}=}{}
+\partial_1\otimes (\widehat{u_2})\wedge
(\widehat{u_1u_2})+
 \partial_1\otimes (\widehat{u_1 v})\wedge (\widehat{u_2\partial_1})+\partial_2\otimes
(\widehat{v})\wedge (\widehat{u_1u_2\partial_2})\\
\hphantom{c_{-4}=}{}
+ \partial_2\otimes
(\widehat{u_1\partial_1})\wedge (\widehat{u_1v}) +\partial_2\otimes
(\widehat{u_2\partial_2})\wedge
(\widehat{u_1v}) +
 \partial_2\otimes (\widehat{u_2v}r)\wedge (\widehat{u_1\partial_2})\\
\hphantom{c_{-4}=}{}
+v\otimes
(\widehat{u_1\partial_1})\wedge (\widehat{u_1u_2})+v\otimes (\widehat{u_2\partial_2})\wedge
(\widehat{u_1u_2}) +v\otimes (\widehat{u_2})\wedge
(\widehat{u_1u_2\partial_2}) \\
\hphantom{c_{-4}=}{}
+ u_1\partial_1\otimes (\widehat{u_2v})\wedge (\widehat{u_1u_2\partial_2})
+u_2\partial_2\otimes (\widehat{u_1v})\wedge (\widehat{u_1u_2\partial_1}) +u_2\otimes
(\widehat{u_1\partial_1})\wedge (\widehat{u_1u_2\partial_1}) \\
\hphantom{c_{-4}=}{}
+ u_2\otimes (\widehat{u_2\partial_2})\wedge (\widehat{u_1u_2\partial_1}) +u_2\otimes
(\widehat{u_2v})\wedge (\widehat{u_1u_2}) +u_2\otimes (\widehat{u_2\partial_1})\wedge
(\widehat{u_1u_2\partial_2})\\
\hphantom{c_{-4}=}{}
+u_1\otimes (\widehat{u_1\partial_1})\wedge
(\widehat{u_1u_2\partial_2}) +
 u_1\otimes (\widehat{u_2\partial_2})\wedge (\widehat{u_1u_2\partial_2}) +u_1\otimes
(\widehat{u_1v})\wedge (\widehat{u_1u_2}) \\
\hphantom{c_{-4}=}{}
+ u_1\otimes (\widehat{u_1\partial_2})\wedge
(\widehat{u_1u_2\partial_1}) +u_2v\otimes (\widehat{u_1u_2})\wedge (\widehat{u_1u_2\partial_1}),\\
%\label{brr}
c_{-2}^1= \partial_2\otimes (\widehat{v})\wedge (\widehat{u_2v})+\partial_2\otimes
(\widehat{u_1\partial_1})\wedge (\widehat{u_2})+ \partial_2\otimes
(\widehat{u_2\partial_2})\wedge (\widehat{u_2}) +v\otimes (\widehat{u_1\partial_1})\wedge
(\widehat{u_2\partial_1})\\
\hphantom{c_{-2}^1=}{}
+ v\otimes (\widehat{u_2\partial_2})\wedge (\widehat{u_2\partial_1})
+u_2\partial_2\otimes (\widehat{u_2})\wedge (\widehat{u_1u_2\partial_1}) +u_2\otimes
(\widehat{u_2v})\wedge (\widehat{u_2\partial_1})\\
\hphantom{c_{-2}^1=}{}
 +u_1\otimes (\widehat{v})\wedge
(\widehat{u_1u_2\partial_1})
+ u_1\otimes (\widehat{u_2})\wedge (\widehat{u_1u_2})+u_1\otimes (\widehat{u_1v})\wedge
(\widehat{u_2\partial_1}) \\
\hphantom{c_{-2}^1=}{}
+u_2v\otimes (\widehat{u_2\partial_1})\wedge
(\widehat{u_1u_2\partial_1})+u_1v\otimes (\widehat{u_1\partial_1})\wedge
(\widehat{u_1u_2\partial_1})+ u_1v\otimes (\widehat{u_2\partial_2})\wedge (\widehat{u_1u_2\partial_1})\\
\hphantom{c_{-2}^1=}{}
+u_1v\otimes
(\widehat{u_2\partial_1})\wedge (\widehat{u_1u_2\partial_2}) +
u_1\partial_2\otimes (\widehat{u_1\partial_1})\wedge (\widehat{u_1u_2})
+ u_1\partial_2\otimes (\widehat{u_2\partial_2})\wedge
(\widehat{u_1u_2})\\
\hphantom{c_{-2}^1=}{}
 +u_1\partial_2\otimes (\widehat{u_2})\wedge (\widehat{u_1u_2\partial_2})+
 u_1\partial_2\otimes (\widehat{u_2 v})\wedge (\widehat{u_1v}) +u_1u_2\otimes (\widehat{u_2
v})\wedge (\widehat{u_1u_2\partial_1})\\
\hphantom{c_{-2}^1=}{}
+u_1u_2\partial_2\otimes (\widehat{u_1u_2})\wedge
(\widehat{u_1u_2\partial_1}),\\
c_{-2}^2= \partial_1\otimes (\widehat{v})\wedge (\widehat{u_1 v})+\partial_1\otimes
(\widehat{u_1\partial_1})\wedge (\widehat{u_1})+\partial_1\otimes
(\widehat{u_2\partial_2})\wedge (\widehat{u_1})+v\otimes (\widehat{u_1\partial_1})\wedge
(\widehat{u_1\partial_2})\\
\hphantom{c_{-2}^2=}{}
+ v\otimes (\widehat{u_2\partial_2})\wedge
(\widehat{u_1\partial_2})+u_1\partial_1\otimes (\widehat{u_1})\wedge
(\widehat{u_1u_2\partial_2})+u_2\otimes (\widehat{v})\wedge
(\widehat{u_1u_2\partial_2})\\
\hphantom{c_{-2}^2=}{}
+u_2\otimes (\widehat{u_1})\wedge (\widehat{u_1u_2})+
 u_2\otimes (\widehat{u_2v})\wedge (\widehat{u_1\partial_2})+u_1\otimes (\widehat{u_1v})\wedge
(\widehat{u_1\partial_2})\\
\hphantom{c_{-2}^2=}{}
+u_2v\otimes (\widehat{u_1\partial_1})\wedge
(\widehat{u_1u_2\partial_2})+
u_2v\otimes (\widehat{u_2\partial_2})\wedge (\widehat{u_1u_2\partial_2})+
 u_2 v\otimes (\widehat{u_1\partial_2})\wedge
(\widehat{u_1u_2\partial_1})\\
\hphantom{c_{-2}^2=}{}
+(u_1v)\otimes (\widehat{u_1\partial_2})\wedge
(\widehat{u_1u_2\partial_2})+u_2\partial_1\otimes (\widehat{u_1\partial_1})\wedge
(\widehat{u_1 u_2})+
 u_2\partial_1\otimes (\widehat{u_2\partial_2})\wedge (\widehat{u_1
u_2})\\
\hphantom{c_{-2}^2=}{}
+u_2\partial_1\otimes (\widehat{u_1})\wedge
(\widehat{u_1u_2\partial_1})+u_2\partial_1\otimes (\widehat{u_2v})\wedge (\widehat{u_1v})+
 u_1 u_2\otimes (\widehat{u_1v})\wedge (\widehat{u_1u_2\partial_2})\\
\hphantom{c_{-2}^2=}{}
 +u_1u_2\partial_1
\otimes (\widehat{u_1 u_2})\wedge
(\widehat{u_1u_2\partial_2}),\\
c_{-2}^3= \partial_2\otimes (\widehat{\partial_1})\wedge (\widehat{u_1
u_2})+\partial_2\otimes (\widehat{u_2\partial_2})\wedge
(\widehat{u_2})+\partial_2\otimes (\widehat{u_1})\wedge (\widehat{u_2\partial_1})+v\otimes
(\widehat{u_2\partial_2})\wedge (\widehat{u_2\partial_1})\\
\hphantom{c_{-2}^3=}{}
+ u_1\partial_1 \otimes (\widehat{u_2\partial_1})\wedge (\widehat{u_1 u_2})+u_1v\otimes
(\widehat{u_2\partial_2})\wedge (\widehat{u_1u_2\partial_1})+u_1v\otimes (\widehat{u_2
v})\wedge (\widehat{u_1 u_2})\\
\hphantom{c_{-2}^3=}{}
+ (u_1v)\otimes (\widehat{u_2\partial_1})\wedge
(\widehat{u_1u_2\partial_2})+u_1\partial_2\otimes (\widehat{u_2\partial_2})\wedge
(\widehat{u_1 u_2}),\\
 c_{2}^1= \partial_1\otimes (\widehat{\partial_1})\wedge
(\widehat{v})+u_1\partial_1\otimes (\widehat{v})\wedge (\widehat{u_1\partial_2 })+u_2\partial_2
\otimes (\widehat{\partial_1})\wedge (\widehat{u_1v})+u_2\otimes (\widehat{\partial_1})\wedge
(\widehat{u_1\partial_1})\\
\hphantom{c_{2}^1=}{}
+ u_1\otimes (\widehat{\partial_1})\wedge (\widehat{u_1\partial_2})+(u_2v)\otimes
(\widehat{\partial_1})\wedge (\widehat{u_1 u_2})+u_2v\otimes (\widehat{u_1\partial_1})\wedge
(\widehat{u_1}) \\
\hphantom{c_{2}^1=}{}
+u_2 v\otimes (\widehat{u_2})\wedge (\widehat{u_1\partial_2})+
 u_2\partial_1 \otimes (\widehat{\partial_1})\wedge
(\widehat{u_1v})+u_2\partial_1\otimes (\widehat{v})\wedge (\widehat{u_1\partial_1
})\\
\hphantom{c_{2}^1=}{}
+u_2\partial_1 \otimes (\widehat{u_2})\wedge (\widehat{u_1})+u_1 u_2\otimes
(\widehat{u_1\partial_1})\wedge
(\widehat{u_1\partial_2})+
 u_1 u_2\otimes (\widehat{u_1})\wedge (\widehat{u_1v})\\
\hphantom{c_{2}^1=}{}
 +u_1u_2\partial_1\otimes
(\widehat{v})\wedge (\widehat{u_1u_2\partial_2})+u_1u_2\partial_1\otimes
(\widehat{u_1\partial_1})\wedge
(\widehat{u_1v}) +
 u_1u_2\partial_1\otimes (\widehat{u_2v})\wedge (\widehat{u_1\partial_2}),\\
c_{2}^2= u_1\partial_1\otimes (\widehat{\partial_1})\wedge (\widehat{u_2v})
+u_1\partial_1\otimes (\widehat{v})\wedge (\widehat{u_2\partial_1})+u_2\otimes
(\widehat{\partial_1})\wedge (\widehat{u_2\partial_1})+u_1\otimes (\widehat{\partial_1})\wedge
(\widehat{u_2\partial_2})\\
\hphantom{c_{2}^2=}{}
 +u_1\otimes (\widehat{v})\wedge (\widehat{u_2})+u_1v\otimes (\widehat{\partial_1})\wedge (\widehat{u_1
u_2})+u_1v\otimes (\widehat{v})\wedge (\widehat{u_2v})+u_1v\otimes
(\widehat{u_2\partial_2})\wedge (\widehat{u_2})\\
\hphantom{c_{2}^2=}{}
+ u_1v\otimes (\widehat{u_1})\wedge (\widehat{u_2\partial_1})+u_1\partial_2\otimes
(\widehat{\partial_1})\wedge (\widehat{u_1v}) +u_1\partial_2 \otimes (\widehat{v})\wedge
(\widehat{u_2\partial_2})\\
\hphantom{c_{2}^2=}{}
+u_1 u_2\otimes (\widehat{u_2\partial_2})\wedge
(\widehat{u_2\partial_1}) +
 u_1 u_2\otimes (\widehat{u_2})\wedge (\widehat{u_2v})+u_1u_2\partial_2\otimes
(\widehat{u_2\partial_2})\wedge
(\widehat{u_2v}) \\
\hphantom{c_{2}^2=}{}
+u_1u_2\partial_2\otimes (\widehat{u_1v})\wedge (\widehat{u_2\partial_1}),\\
c_{4}^1= u_2\partial_2 \otimes (\widehat{\partial_1})\wedge (\widehat{u_2})+u_2
v\otimes (\widehat{\partial_1})\wedge (\widehat{u_1\partial_1})+u_2v\otimes
(\widehat{\partial_1})\wedge (\widehat{u_2\partial_1})+u_2v\otimes (\widehat{v})\wedge
(\widehat{u_2})\\
\hphantom{c_{4}^1=}{}
+ u_1v\otimes (\widehat{\partial_1})\wedge
(\widehat{u_1\partial_1})+u_1u_2\partial_1\otimes (\widehat{u_1\partial_1})\wedge
(\widehat{u_2})+u_1u_2\partial_2\otimes (\widehat{\partial_1})\wedge (\widehat{u_1
u_2})\\
\hphantom{c_{4}^1=}{}
+ u_1u_2\partial_2\otimes (\widehat{u_1\partial_1})\wedge
(\widehat{u_1})+u_1u_2\partial_2\otimes (\widehat{u_2})\wedge
(\widehat{u_1\partial_2}),\\
c_4^2= u_1\partial_1\otimes (\widehat{\partial_1})\wedge (\widehat{u_1})+u_2v\otimes
(\widehat{\partial_1})\wedge (\widehat{u_2\partial_2})+u_1v\otimes (\widehat{\partial_1})\wedge
(\widehat{u_1\partial_2})\\
\hphantom{c_4^2=}{}
+(u_1v)\otimes (\widehat{\partial_1})\wedge
(\widehat{u_2\partial_2})+
 u_1v\otimes (\widehat{v})\wedge (\widehat{u_1})+u_1u_2\partial_1\otimes
(\widehat{\partial_1})\wedge (\widehat{u_1 u_2})\\
\hphantom{c_4^2=}{}
+u_1u_2\partial_1\otimes
(\widehat{u_2\partial_2})\wedge (\widehat{u_2})+
 u_1u_2\partial_1\otimes (\widehat{u_1})\wedge
(\widehat{u_2\partial_1})+u_1u_2\partial_2\otimes (\widehat{u_2\partial_2})\wedge
(\widehat{u_1}),\\
c_4^3= u_2v\otimes (\widehat{\partial_1})\wedge
(\widehat{u_1\partial_2})+u_2\partial_1\otimes (\widehat{\partial_1})\wedge
(\widehat{u_1})+u_1u_2\partial_1 \otimes (\widehat{u_1})\wedge
(\widehat{u_1\partial_2}),\\
c_4^4= u_1v\otimes (\widehat{\partial_1})\wedge
(\widehat{u_2\partial_1})+u_1\partial_2\otimes (\widehat{\partial_1})\wedge
(\widehat{u_2})+u_1u_2\partial_2 \otimes (\widehat{u_2})\wedge (\widehat{u_2\partial_1}).
\end{gather*}

\newpage

\section{Queer Lie superalgebras and queerified Lie algebras}

In the textbook \cite{Ls}, it is demonstrated that there are at
least two\footnote{Both these Lie superalgebras $\mathfrak{gl}(V)$ (respectively, $\mathfrak{gl}(V)$) can be considered---if the superdimension of the superspace $V$ is equal to $2^{s-1}|2^{s-1}$---as ``quantized'' versions of one more analog of $\mathfrak{gl}(n)$, namely of the
Poisson Lie superalgebra $\mathfrak{po}(0|2s)$ (respectively, $\mathfrak{po}(0|2s+1)$), which corresponds to the value 0 of the parameter in the 1-parametric family of its deforms.} super versions of
$\mathfrak{gl}(n)$: a~naive one, $\mathfrak{gl}(n|m)$, and the ``queer" one,
\begin{align*}
\mathfrak{q}(n):={} & \left\{X\in \mathfrak{gl}(n|n)\mid [X, \Pi_{2n}]=0 \ \text{for
$\Pi_{2n}=\begin{pmatrix}
0 &1_n\\
-1_n &0
\end{pmatrix}$} \right\}\\
={}& \left\{X=(A,B):=
\begin{pmatrix}
A &B\\
-B &A
\end{pmatrix}\right\},
\end{align*}
which preserves a~complex structure given by the odd operator $\Pi_{2n}$. We
set
\[
\mathfrak{sq}(n):= \{X=(A,B)\in \mathfrak{q}(n)\mid \mathrm{tr}(B)=0 \}.
\]
Let $\mathfrak{psq}(n):=\mathfrak{sq}(n)/\Kee 1_{2n}$ be the
projectivization of $\mathfrak{sq}(n)$. We denote the images of the
$A_{ij}$-elements of $A\in\mathfrak{sq}(n)_\ev$ in
$\mathfrak{psq}(n)$ by $a_{ij}$ and the images of the
$B_{ij}$-elements of $B\in\mathfrak{sq}(n)_\od$ by~$b_{ij}$; let $b_{i}:=b_{ii}-b_{i+1, i+1}$.

If $\operatorname{char}\Kee=p>2$, then the derived superalgebra
$\mathfrak{psq}^{(1)}(n)=[\mathfrak{psq}(n), \mathfrak{psq}(n)]$ is
simple and is only different from $\mathfrak{psq}(n)$ for $n=pk$,
when $\mathfrak{sq}(n)$ contains $\Pi_{2n}$.

\begin{Lemma}\label{L6.1} For $\fg=\mathfrak{psq}(3)$, we have $H^2(\fg;\fg)=0$ for
$p\geq 5$.

For $p=3$, when $\fg=\mathfrak{psq}(3)$ is not simple, we take the simple
algebra $\fg'=\mathfrak{psq}^{(1)}(3)$ and for a~basis of $H^2(\fg';\fg')$ we
take the $2$-cocycle
\begin{gather*}%\label{eqpsq3}
\underline{c_0}= 2a_{2,2}\otimes \widehat a_{2,2}\wedge \widehat b_{2}
+ 2a_{2,2}\otimes \widehat a_{3,3}\wedge \widehat b_{2}+
2a_{2,2}\otimes \widehat a_{1,2}\wedge \widehat b_{2,1}+
a_{2,2}\otimes \widehat a_{1,3}\wedge
\widehat b_{3,1} \\
\hphantom{\underline{c_0}=}{}
+ a_{2,2}\otimes \widehat a_{3,2}\wedge \widehat b_{2,3}+
 a_{2,2}\otimes \widehat a_{3,1}\wedge \widehat b_{1,3} +
2a_{3,3}\otimes \widehat a_{2,2}\wedge \widehat b_{1}+ 2a_{3,3}\otimes
\widehat a_{2,2}\wedge \widehat b_{2}\\
\hphantom{\underline{c_0}=}{}
+ 2a_{3,3}\otimes
\widehat a_{3,3}\wedge \widehat b_{1}+ 2a_{3,3}\otimes
\widehat a_{3,3}\wedge
\widehat b_{2}+
 a_{3,3}\otimes \widehat a_{2,1}\wedge \widehat b_{1,2}+
2a_{3,3}\otimes \widehat a_{3,2}\wedge \widehat b_{2,3}\\
\hphantom{\underline{c_0}=}{}
 +
a_{3,3}\otimes \widehat a_{3,1}\wedge \widehat b_{1,3}+
2a_{1,2}\otimes \widehat a_{1,2}\wedge \widehat b_{2}+
2a_{1,2}\otimes \widehat a_{3,2}\wedge
\widehat b_{1,3}+
 a_{2,3}\otimes \widehat a_{2,3}\wedge \widehat b_{2}\\
\hphantom{\underline{c_0}=}{}
 + a_{1,3}\otimes
\widehat a_{2,2}\wedge \widehat b_{1,3}+ a_{1,3}\otimes
\widehat a_{3,3}\wedge \widehat b_{1,3}+ a_{1,3}\otimes
\widehat a_{2,3}\wedge \widehat b_{1,2}+ a_{1,3}\otimes
\widehat a_{1,3}\wedge \widehat b_{2} \\
\hphantom{\underline{c_0}=}{}
+ 2a_{2,1}\otimes \widehat a_{2,1}\wedge \widehat b_{1}+ a_{3,2}\otimes
\widehat a_{3,2}\wedge \widehat b_{1}+ a_{3,1}\otimes
\widehat a_{2,2}\wedge \widehat b_{3,1}+ a_{3,1}\otimes
\widehat a_{3,3}\wedge \widehat b_{3,1}
\\
\hphantom{\underline{c_0}=}{}
+ 2a_{3,1}\otimes \widehat a_{3,1}\wedge \widehat b_{1}+ a_{3,1}\otimes \widehat a_{3,1}\wedge \widehat b_{2}+ b_{1}\otimes
\widehat b_{1,3}\wedge \widehat b_{3,1}+ 2b_{1}\otimes
\widehat b_{1}\wedge \widehat b_{1} + b_{2}\otimes
\widehat b_{2}\wedge \widehat b_{2}\\
\hphantom{\underline{c_0}=}{}
+ b_{1,2}\otimes \widehat a_{2,2}\wedge \widehat a_{1,2}+
 b_{1,2}\otimes \widehat a_{3,3}\wedge \widehat a_{1,2}+
b_{1,2}\otimes \widehat a_{1,3}\wedge \widehat a_{3,2} +
b_{1,2}\otimes \widehat b_{1,3}\wedge \widehat b_{3,2}\\
\hphantom{\underline{c_0}=}{}
+
2b_{2,3}\otimes \widehat a_{2,2}\wedge \widehat a_{2,3}+ 2b_{2,3}\otimes \widehat a_{3,3}\wedge \widehat a_{2,3}+
 2b_{2,3}\otimes \widehat a_{1,3}\wedge \widehat a_{2,1} +
b_{2,3}\otimes \widehat b_{1}\wedge \widehat b_{2,3}\\
\hphantom{\underline{c_0}=}{}
+ 2b_{2,3}\otimes
\widehat b_{2}\wedge \widehat b_{2,3} + 2b_{1,3}\otimes
\widehat a_{1,2}\wedge \widehat a_{2,3}+ 2b_{1,3}\otimes
\widehat b_{1}\wedge \widehat b_{1,3}
+
 2b_{2,1}\otimes \widehat a_{2,2}\wedge \widehat a_{2,1}\\
\hphantom{\underline{c_0}=}{}
 +
2b_{2,1}\otimes \widehat a_{3,3}\wedge \widehat a_{2,1} +
2b_{2,1}\otimes \widehat a_{2,3}\wedge \widehat a_{3,1} +
b_{2,1}\otimes \widehat b_{1}\wedge \widehat b_{2,1}+ 2b_{2,1}\otimes \widehat b_{2}\wedge \widehat b_{2,1}\\
\hphantom{\underline{c_0}=}{}
+
 b_{2,1}\otimes \widehat b_{2,3}\wedge \widehat b_{3,1}+
b_{3,2}\otimes \widehat a_{2,2}\wedge \widehat a_{3,2}+ b_{3,2}\otimes
\widehat a_{3,3}\wedge \widehat a_{3,2} + b_{3,2}\otimes
\widehat a_{1,2}\wedge \widehat a_{3,1} \\
\hphantom{\underline{c_0}=}{}
+ b_{3,2}\otimes
\widehat b_{1,2}\wedge \widehat b_{3,1}+
 2b_{3,1}\otimes \widehat a_{2,1}\wedge \widehat a_{3,2}.
\end{gather*}
\end{Lemma}

\begin{Lemma}\label{L6.3} For $\fg=\mathfrak{psq}(4)$, we have $H^2(\fg;\fg)=0$ for
$p=3$ and for $p=5$.
\end{Lemma}

\begin{Lemma}\label{L6.2} For $\fg=\mathfrak{q}(\mathfrak{sl}(3))$, for a~basis of
$H^2(\fg;\fg)$ we take the $2$-cocycles
\begin{gather}
\underline{A}= h_1\otimes (\widehat y_1\wedge \Pi(\widehat x_1))+h_1\otimes
(\widehat y_3\wedge \Pi(\widehat x_3))+ h_2\otimes (\widehat x_1\wedge \Pi(\widehat y_1))+h_2\otimes
(\widehat y_1\wedge \Pi(\widehat x_1))
\nonumber\\
\hphantom{\underline{A}=}{}
+ h_2\otimes (\widehat y_3\wedge \Pi(\widehat x_3))+x_1\otimes (\widehat x_1\wedge
\Pi(\widehat h_1))+ x_2\otimes (\widehat x_2\wedge \Pi(\widehat h_1))+ x_2\otimes
(\widehat x_3\wedge \Pi(\widehat y_1))
\nonumber\\
\hphantom{\underline{A}=}{}
+ y_1\otimes (\widehat y_1\wedge \Pi(\widehat h_1))+y_1\otimes (\widehat y_1\wedge
\Pi(\widehat h_2))+y_1\otimes (\widehat y_3\wedge \Pi(\widehat x_2))+ y_2\otimes (\widehat x_1\wedge
\Pi(\widehat y_3))
\nonumber\\
\hphantom{\underline{A}=}{}
+ y_2\otimes (\widehat y_2\wedge \Pi(\widehat h_1))+ y_3\otimes (\widehat y_1\wedge
\Pi(\widehat y_2))+y_3\otimes (\widehat y_3\wedge \Pi(\widehat h_1))+y_3\otimes (\widehat y_3\wedge
\Pi(\widehat h_2))
\nonumber\\
\hphantom{\underline{A}=}{}
+ \Pi(h_1)\otimes (\Pi(\widehat h_1))^{\wedge^2} +\Pi(h_2)\otimes
(\Pi(\widehat x_2)\wedge \Pi(\widehat y_2))+\Pi(h_2)\otimes (\Pi(\widehat h_2))^{\wedge^2}
\nonumber\\
\hphantom{\underline{A}=}{}
+\Pi(x_1)\otimes (\Pi(\widehat h_1)\wedge \Pi(\widehat x_1))+ \Pi(x_1)\otimes (\Pi(\widehat h_2)\wedge \Pi(\widehat x_1))+ \Pi(x_1)\otimes
(\Pi(\widehat x_3)\wedge \Pi(\widehat y_2))
\nonumber\\
\hphantom{\underline{A}=}{}
+ \Pi(x_3)\otimes (\Pi(\widehat h_1)\wedge
\Pi(\widehat x_3))+
 \Pi(x_3)\otimes (\Pi(\widehat h_2)\wedge \Pi(\widehat x_3))+ \Pi(x_3)\otimes
(\Pi(\widehat x_1)\wedge \Pi(\widehat x_2))
\nonumber\\
\hphantom{\underline{A}=}{}
+
\Pi(y_1)\otimes (\Pi(\widehat h_1)\wedge \Pi(\widehat y_1)),
\nonumber\\
 B=
 h_1\otimes(\Pi (\widehat h_1))^{\wedge 2}+h_2\otimes(\Pi (\widehat h_2))^{\wedge 2}
 +\Pi (x_1)\otimes(\widehat x_1\wedge \Pi (\widehat h_2))+\Pi (x_2)\otimes(\widehat x_2\wedge \Pi (\widehat h_1))
\nonumber\\
\hphantom{B=}{}
 +\Pi (x_3)\otimes(\widehat x_3\wedge \Pi (\widehat h_1))+\Pi (x_3)\otimes(\widehat x_3\wedge \Pi(\widehat h_2))
 +\Pi (h_1)\otimes(\widehat y_1\wedge \Pi(\widehat x_1))
\nonumber\\
\hphantom{B=}{}
 +\Pi (h_2)\otimes(\widehat y_2\wedge \Pi(\widehat x_2)) + h_1\otimes(\widehat x_3\wedge \widehat y_3)
 +h_2\otimes(\widehat x_3\wedge \widehat y_3)+y_1\otimes(\Pi(\widehat h_2)\wedge \Pi(\widehat y_1))
\nonumber\\
\hphantom{B=}{}
 +y_2\otimes(\Pi(\widehat h_1)\wedge \Pi(\widehat y_2))+\Pi (y_3)\otimes(\widehat y_3\wedge \Pi(\widehat h_1))
 +\Pi (y_3)\otimes(\widehat y_3\wedge \Pi (\widehat h_2))
\nonumber\\
\hphantom{B=}{}
 +x_3\otimes(\Pi (\widehat x_1)\wedge \Pi(\widehat x_2))+\Pi (x_1)\otimes(\widehat x_3\wedge \Pi(\widehat y_2))
 +\Pi (x_2)\otimes(\widehat x_3\wedge \Pi (\widehat y_1))
\nonumber\\
\hphantom{B=}{}
 +y_1\otimes(\widehat x_2\wedge \widehat y_3)+y_2\otimes(\widehat x_1 \wedge \widehat y_3)
 +y_3\otimes(\widehat y_1 \wedge \widehat y_2).\label{eqpsq311}
\end{gather}
\end{Lemma}

\section{On integrability of infinitesimal deformations}\label{Sintegra}

The main tool in proving integrability of infinitesimal deformations is the technique of Massey products, cf.~\cite{BLW}.
Superization of the relevant definitions is supposed to be straightforward (just apply the sign rule) but is not quite (as, for example, is the definition of the supertrace that differs on the even and odd supermatrices), so we give a~precise definition; compare with the non-super case in~\cite{BLW}.

Let $\mathfrak{g}_t$ be a~$1$-parameter deformation of a~Lie superalgebra
$\mathfrak{g}$, given by an infinitesimal cocycle $c=c^1$ and higher
degree terms $c^2$, $c^3$, \dots. The Jacobi identity modulo
$t^{n+1}$ reads
\begin{gather}
\sum\limits_{i+j=n,\, i, j\geq 0}c^i(c^j(x,y),z)+(-1)^{p(z)(p(y)+p(x))}c^i(c^j(z,x),y)\nonumber\\
\qquad{}+(-1)^{p(x)(p(y)+p(z))}c^i(c^j(y,z),x) +(c^i\leftrightarrow c^j)=0.\label{br}
\end{gather}
The expression \eqref{br} can be rewritten as
\[
0= \sum\limits_{0\leq i<j\leq n;\, i+j=n}[[c^i,c^j]]+\begin{cases}0&\text{if $n=2k+1$},\\
 c^k\circ c^k&\text{if $n=2k$},\end{cases}
\]
where the non-super versions of the brackets
\begin{gather*}
[[c^i,c^j]](x,y,z):= c^i(c^j(x,y),z)+
 (-1)^{p(z)(p(y)+p(x))}c^i(c^j(z,x),y)\\
\hphantom{[[c^i,c^j]](x,y,z):=}{}
+(-1)^{p(x)(p(y)+p(z))}c^i(c^j(y,z),x)+(c^i\leftrightarrow c^j)
\end{gather*}
are called \textit{Nijenhuis
brackets} (in differential geometry) or \textit{Massey brackets} (in
deformation theory) and where
\begin{gather*}
c^k\circ c^k (x,y,z):=c^k(c^k(x,y),z)+(-1)^{p(z)(p(y)+p(x))} c^k(c^k(z,x),y)\\
\hphantom{c^k\circ c^k (x,y,z):=}{} +(-1)^{p(x)(p(y)+p(z))} c^k(c^k(y,z),x)
\end{gather*}
are \textit{Massey squares}. The brackets and squares are similarly defined for
cocycles of higher degrees. The collection of all brackets
$[[\cdot,\cdot]]$, together with squares if $p=2$ (this Massey square exists even when $p>2$ and has nothing to do with the squaring $=$ ``half bracket'' of the odd element with itself) defines a~graded
Lie superalgebra structure on $H^{\bcdot}(\mathfrak{g};
\mathfrak{g})$. The formula~\eqref{br} can be expressed as a~\textit{Maurer--Cartan}
equation: \begin{gather*}%\label{MaCart}
dc^n=\sum_{0< i<j\leq n;\, i+j=n}[[c^i,c^j]]+\begin{cases}0&\text{if $n=2k+1$},\\
 c^k\circ c^k&\text{if $n=2k$}.
 \end{cases}
\end{gather*}
Let $\{a_1, \ldots, a_n\}$ be a~basis of $\fg$. Suppose that $c=a_l\otimes \widehat a_m\wedge \widehat a_n$ and $\tilde c=a_i\otimes \widehat a_j\wedge \widehat a_k$. We have
\begin{gather*}
c(\tilde c(x,y),z)+(-1)^{p(z)(p(y)+p(x))}c(\tilde c(x,y),z)+(-1)^{p(x)(p(y)+p(z))}c(\tilde c(y,z),x)\\
\qquad{}=
\big(\delta^i_m a_l \otimes \widehat a_j\wedge \widehat a_k\wedge \widehat a_n-(-1)^{p(a_m)p(a_n)} \delta^i_na_l\otimes \widehat a_j\wedge \widehat a_k\wedge \widehat a_m+\delta^j_l a_i \otimes \widehat a_m\wedge \widehat a_n\wedge \widehat a_k\\
\qquad\quad{} -(-1)^{p(a_j)p(a_k)} \delta^l_k a_i\otimes \widehat a_m\wedge \widehat a_n\wedge \widehat a_j\big) (x,y,z).
\end{gather*}

Because of the squaring, the definition of Massey brackets must be different in characteristic~2; for the correct definition, see \cite{BLW,FL}.

The above gives a~clear procedure for the prolongation of an
infinitesimal deformation (expressed here for simplicity only for a
$1$-parameter deformation): given a~first degree deformation via a
cocycle $c=c^1$, one must compute its \textit{Massey square} $[[c,c]]$. The following cases can be encountered:

If $[[c,c]]=0$, the infinitesimal deformation fulfills the Jacobi
identity, and therefore, is a~true deformation.

If $[[c,c]]\in
Z^3(\mathfrak{g},\mathfrak{g})$ is not a~coboundary, the
infinitesimal deformation is obstructed and cannot be prolonged. If
 $[[c,c]]=d\alpha$ with $\alpha\not=0$, then
$-\alpha t^2$ is the second degree term of the deformation. In order
to prolong to the third degree, one has to compute the next step---the Massey product $[[c,\alpha]]$. Once again, there are the three
possibilities
\[
1)\ [[c,\alpha]]=0,\qquad 2)
\ [[c,\alpha]]=d\beta \ \text{with} \ \beta\not=0, \qquad
3)\ [[c,\alpha]]\not=d\beta \ \text{for any $\beta$.}
\]

If
$[[c,\alpha]]=d\beta$, then $-\beta t^3$ gives the third degree
prolongation of the deformation. In order to go up to degree $4$
then, one has to be able to compensate
$[[\alpha,\alpha]]+[[c,\beta]]$ by a~coboundary $d\gamma$, and so
on. One must be careful to keep track of all terms coming in to
compensate low degree Massey brackets in a~multiparameter
deformation.

\begin{Lemma}\label{integra3} For $p=2$, the Lie algebra
$\mathfrak{o}^{(1)}(3)$ admits the global deformation given by
\[
[\cdot ,\cdot ]_{\alpha,\beta}=[\cdot ,\cdot ]+\alpha c_1 +\beta
c_2.
\]
Denote the deform by $\mathfrak{o}^{(1)}(3,\alpha,\beta)$. This Lie algebra
is simple if and only if $\alpha \beta \not =1$.

If $\alpha \beta=1$, then $\mathfrak{o}^{(1)}(3,\alpha,\beta)$ has an ideal
$I=\operatorname{Span} \{h, x+\alpha y \}$.

For $\alpha \beta\not = 1$, the Lie algebra
$\mathfrak{o}^{(1)}(3,\alpha,\beta)$ is simple with two outer derivations
given by
\[
\alpha h \otimes \widehat h + \alpha x\otimes \widehat x + x\otimes \widehat y, \qquad \beta
h\otimes \widehat h + \beta x\otimes \widehat x + y\otimes \widehat x.
\]
\end{Lemma}

Reducing the operator $\operatorname{ad}_h$ to the Jordan normal form, we
immediately see (in the eigenbasis) that the deform depends,
actually, on one parameter, not two. It is not difficult to show
that each of these $1$-parameter deforms is isomorphic to the
initial algebra; so these deforms are \emph{semi-trivial}, see~\cite{BLLS}.

Proof of Lemma \ref{integra3} is straightforward. The next two
lemmas were proved with the aid of the \textsc{SuperLie} code, cf.~\cite{Gr}.

Every cocycle of positive degree can be extended to a~global
deformation by means of Massey products as follows from the same argument
people use to prove the existence of the highest and lowest weight vector
of any finite-dimensional module over any finite-dimensional simple
Lie algebra over $\Cee$.

\begin{Lemma}\label{integra2} Let $p=2$, and $\alpha\neq 0, 1$.
The bracket
\eqref{defglob} satisfies the Jacobi identity in the following
cases:
\begin{enumerate}\itemsep=0pt
\item[$(a)$] For $\mathfrak{o}^{(1)}(3)$, $\mathfrak{o}^{(1)}(5)$, $\mathfrak{psl}(4)$, $\mathfrak{psl}(6)$, and
 $\mathfrak{o}\mathfrak{o}_{I \Pi}^{(1)}(1|2)$, for each of the homogeneous
cocycles $c$.
\item[$(b)$] For $\fwk(3;\alpha)$ and $\fwk(4;\alpha)$, and also for
$\fbgl(3;\alpha)$ and $\fbgl(4;\alpha)$, for each of the homogeneous
cocycles $c$, except the ones of degree~$0$.
\end{enumerate}
\end{Lemma}

\begin{Lemma}\label{integral3} Let $p=2$, and $\alpha\neq 0, 1$.
\begin{enumerate}\itemsep=0pt
\item[$1)$] For $\fwk(3;\alpha)$, the global
deform corresponding to $c_0$, see \eqref{eqwk3a}, is as follows:
\begin{gather}\label{wk3glob}
{}[\cdot,\cdot]_t=[\cdot,\cdot]+tc_0+t^2h_1\otimes \widehat x_7\wedge \widehat
y_7.
\end{gather}
For $\fbgl(3;\alpha)$, the global
deform corresponding to $c_0$, see~\eqref{eqbgl3a}, is of a~form similar to~\eqref{wk3glob}.

\item[$2)$] For $\fwk(4;\alpha)$, the global
deform corresponding to $c_0$, see \eqref{eqpsl4c3}, is as follows:
\begin{gather}\label{65}
[\cdot,\cdot]_t=[\cdot,\cdot]+tc_0+t^2A+ t^3B,
\end{gather}
where
\begin{gather*}
A = \alpha h_1\otimes \widehat x _{14} \wedge
\widehat y_{14}
+\alpha ^2
h_1\otimes \widehat x_{15} \wedge
\widehat y_{15}
+h_2\otimes \widehat x_{10} \wedge
\widehat y_{10}
+h_2\otimes \widehat x_{12} \wedge
 \widehat y_{12}\\
\hphantom{A=}{}
 +h_2\otimes \widehat x_{13} \wedge
 \widehat y_{13}
 +(\alpha +1)
 h_2\otimes \widehat x_{14} \wedge
 \widehat y_{14}
+\alpha
h_3\otimes \widehat x_{14} \wedge
\widehat y_{14}
+\alpha ^2
 h_3\otimes \widehat x_{15} \wedge
 \widehat y_{15}\\
\hphantom{A=}{}
 +\alpha
 h_4\otimes \widehat x_{14} \wedge
 \widehat y_{14}
 +\alpha ^2
 h_4\otimes \widehat x_{15} \wedge
 \widehat y_{15}
 +\alpha
 x_1\otimes \widehat x_{14} \wedge
 \widehat y_{13}\\
\hphantom{A=}{}
 +
\big(\alpha ^2+\alpha \big)
 x_2\otimes \widehat x_{15} \wedge
 \widehat y_{14} +\alpha
 x_5\otimes \widehat x_{15} \wedge
 \widehat y_{13}
 +\alpha
 x_6\otimes \widehat x_{14} \wedge
 \widehat y_{12}\\
\hphantom{A=}{}
 +\alpha
 x_8\otimes \widehat x_{15} \wedge
 \widehat y_{12}
 +\alpha
 x_9\otimes \widehat x_{14} \wedge
 \widehat y_{10}
 +
 \alpha
 x_{11}\otimes \widehat x_{15} \wedge
 \widehat y_{10}
 +\alpha
 y_1\otimes \widehat x_{13} \wedge
 \widehat y_{14}\\
\hphantom{A=}{}
 + \big(\alpha ^2+\alpha \big)
 y_2\otimes \widehat x_{14} \wedge
 \widehat y_{15}
 +\alpha
 y_5\otimes \widehat x_{13} \wedge
 \widehat y_{15}
 +\alpha
 y_6\otimes \widehat x_{12} \wedge
 \widehat y_{14}
 +
 \alpha
 y_8\otimes \widehat x_{12} \wedge
 \widehat y_{15}\\
\hphantom{A=}{}
 +\alpha
 y_9\otimes \widehat x_{10} \wedge
 \widehat y_{14}
 +\alpha
 y_{11}\otimes \widehat x_{10} \wedge
 \widehat y_{15}, \\
B = h_2\otimes \widehat x_{14} \wedge
 \widehat y_{14}
 +\alpha
 h_3\otimes \widehat x_{14} \wedge
 \widehat y_{14}
 +\alpha ^2
 h_3\otimes \widehat x_{15} \wedge
 \widehat y_{15}
 +\alpha
 x_2\otimes \widehat x_{15} \wedge
 \widehat y_{14}\\
 \hphantom{B=}{}
 +\alpha
 y_2\otimes \widehat x_{14} \wedge
 \widehat y_{15}.
\end{gather*}
For $\fbgl(4;\alpha)$, the answer, and its form, are similar to formula~\eqref{65}.

\item[$3)$] For $\mathfrak{q}(\mathfrak{sl}(3))$, the global deform corresponding to the even cocyle $B$, see \eqref{eqpsq311}, is
as follows:
$%\label{qsl3glob}
[\cdot,\cdot]_t=[\cdot,\cdot]+t B +t^2 C$,
where
\begin{gather*}
 C = h_1\otimes(\Pi (\widehat x_1)\wedge \Pi(\widehat y_1)) + h_2\otimes(\Pi(\widehat x_2)\wedge\Pi(\widehat y_2))
 +\Pi (h_1)\otimes(\widehat x_1\wedge \Pi(\widehat y_1))\\
 \hphantom{C =}{}
 +\Pi (h_2)\otimes(\widehat x_2\wedge\Pi(\widehat y_2))
 +h_1\otimes(\widehat x_1\wedge \widehat y_1)+h_2\otimes(\widehat x_2\wedge\widehat y_2)
\\
 \hphantom{C =}{}
+\Pi (x_3)\otimes(\widehat x_1\wedge \Pi(\widehat x_2)) +\Pi (x_3)\otimes(\widehat x_2\wedge\Pi(\widehat x_1))
 +x_3\otimes(\widehat x_1\wedge\widehat x_2)\\
 \hphantom{C =}{}
+y_3\otimes(\Pi(\widehat y_1)\wedge\Pi(\widehat y_2)) +\Pi (y_3)\otimes(\widehat y_1\wedge\Pi(\widehat y_2))+\Pi (y_3)\otimes(\widehat y_2\wedge\Pi(\widehat y_1)).
\end{gather*}
\end{enumerate}
\end{Lemma}

\begin{Lemma}\label{integra} For $p=3$, the bracket
\eqref{defglob} satisfies the Jacobi identity in the following cases
\emph{(for the definition of the respective algebras, see \cite{BGL2})}:
for $\fg=\mathfrak{br}(2; \varepsilon)$, where $\varepsilon\neq 0$, as well as
for $\mathfrak{brj}(2,3)$, $\mathfrak{br}(2;5)$, and $\mathfrak{br}(3)$, for each of the
homogeneous cocycles given in respective subsections above.
\end{Lemma}

\appendix
\section[Cohomology of simple Lie algebras with coefficients in the trivial module]{Cohomology of simple Lie algebras with coefficients\\ in the trivial module}\label{App}

\subsection{Useful facts}\label{ssHS}
Recall that over $\Cee$ the cohomology of simple Lie algebras with
coefficients in the trivial module, are known\footnote{The published
proofs known to us reduce the computation to compact forms of the
complex Lie groups with given Lie algebra, but using
Hochschild--Serre spectral sequence and induction on rank one can
easily obtain the answer algebraically; the base of induction being $\mathfrak{gl}(1)$ and
$\mathfrak{sl}(2)$ for which computations are easy by hand.} \cite{Car}:
Each space is naturally
endowed with a~Grassmann algebra structure generated by cocycles
$c_{2e_i+1}$ of degree\footnote{The numbers $2e_i+1$ are dimensions of
the irreducible $\mathfrak{sl}(2)$-modules in the adjoint action in $\fg$ of the
principally embedded $\mathfrak{sl}(2)\subset\fg$. Recall that the
\textit{principal} embedding is the one with the least (equal to $\operatorname{rk}\fg$) number
of irreducible $\mathfrak{sl}(2)$-modules into which $\fg$ decomposes, see \cite{GL2}.} $2e_i+1$, where the
\textit{exponents} $e_i$ are as follows:
\begin{equation}\label{CohSimple}\renewcommand{\arraystretch}{1.2}
\begin{tabular}{|l|l|}
\hline
Cartan's or Bourbaki's notation of $\fg$&the exponents $e_i$ \cite[Table 4]{OV} \cr
\hline
$A_n$ or $\mathfrak{sl}(n+1)$&$ 1, 2, 3,\dots, n$
\cr
$\begin{cases}B_n \ \text{or $\mathfrak{o} (2n+1)$}\cr
C_n \ \text{or $\mathfrak{sp}(2n)$}\end{cases}$
for $n\geq 2$&$1, 3, \dots,
2n-1$\cr
$D_n$ or $\mathfrak{o}(2n)$ &$1, 3, \dots, 2n-3; n-1$\cr
$G_2$ or $\fg(2)$ &$1, 5$\cr
$F_4$ or
$\mathfrak{f}(4)$&$1, 5, 7, 11$\cr
$E_6$ or $\mathfrak{e}(6)$&$1, 4, 5, 7, 8,
11$\cr
$E_7$ or $\mathfrak{e}(7)$&$1, 5, 7, 9, 11, 13, 17$\cr
$E_8$
or $\mathfrak{e}(8)$&$1, 7, 11, 13, 17, 19, 23, 29$\cr \hline
\end{tabular}
\end{equation}

Obviously, if $p$ is sufficiently big as compared with $\dim\fg$,
the fact \eqref{CohSimple} is true over $\Kee$ of characteristic
$p$, but for $p$ ``small'', the analog of the fact \eqref{CohSimple} is unknown.
To investigate this,
we will use the following facts (see~\cite{Fu}):
\begin{gather*}%\label{3facts}
H^{\bcdot}(\mathfrak{gl}(1))=\Lambda^{\bcdot}(c_1)\qquad \text{for the dual element $c_1:=1^*\in\mathfrak{gl}(1)^*$},\\
H^{n}(\fg_1\oplus\fg_2)\simeq
\mathop{\oplus}\limits_{a+b=n}H^{a}(\fg_1)\otimes H^{b}(\fg_2),\\
H^{n}(\fg; M_1\oplus M_2)\simeq H^{n}(\fg; M_1)\oplus H^{n}(\fg;M_2).
\end{gather*}

\subsection[Analogs of simple complex Lie algebras for $p$ much larger than $\dim\fg$]{Analogs of simple complex Lie algebras for $\boldsymbol{p}$ much larger than $\boldsymbol{\dim\fg}$}\label{ssPlarge}

Clearly, in this case, the cohomology $H^{\bcdot}(\fg)$ are described as
for $p=0$, see \eqref{CohSimple}. We did not know how does the
answer change for $p$ small: do degrees of generators of $H^{\bcdot}(\fg)$
vary, or $H^{\bcdot}(\fg)$ fails to be a~free supercommutative superalgebra,
or both? Here are the answers (only for the small values of $p$ for
which the answer differ from $p=0$), we give cocycles of nonpositive
degrees only (by symmetry this suffices). For the principal grading
of $\fg(A)$, we have the following answers obtained with the aid of \textsc{SuperLie}, see~\cite{Gr}.

For $\fg(2)$ and $p=5$, the dimensions and degrees of basis
cocycles are as follows:
\begin{alignat*}{3}%\label{q2}
&\deg = - 10\colon \quad && \text{one cocycle at each }H^5, \ H^6, \ H^8, \ H^9,&\\
&\deg = - 5\colon \quad &&\text{one cocycle at each }H^5, \ H^6, \ H^8, \ H^9,&\\
&\deg = 0\colon \quad &&\text{one cocycle at each }H^0, \ H^3, \ H^5, \ H^6, \ H^8, \ H^9, \ H^{11}, \ H^{14}.&
\end{alignat*}
Denote the cocycles of degree $5k$ in $H^i$ by $c_{5k}^i$. The
generators\footnote{The equality sign ``$=$'' should be understood as
equality of classes, rather than of cocycles themselves.} in $H^{\bcdot}$
are as follows, see~\eqref{g2genRel}. The relations are too numerous
and cumbrous to figure them out explicitly without motivations.
Below, the ``equalities'' of the type $a\times b=ab$ should
eventually be replaced by genuine relations, if there are any (in
addition to supercommutativity); for example, if
\begin{gather*}
c_{0}^{11}=gc_{0}^5c_{0}^6+ac_{5}^5c_{-5}^6+bc_{-5}^5c_{5}^6+
ec_{10}^5c_{-10}^6+fc_{-10}^5c_{10}^6\qquad \text{for some coefficients $a$, $b$, $e$, $f$, $g$,}
\end{gather*}
then
$c_{0}^{11}$ is not a~generator,
and there are 5 more relations. We did not investigate this.
\begin{equation}\label{g2genRel}
\renewcommand{\arraystretch}{1.2}
\begin{tabular}{|l|l|}
\hline generators&$c_{0}^3$, $c_{0}^5$, $c_{0}^6$, $c_{0}^{11}$?\\
&$c_{\pm5}^5$, $c_{\pm5}^6$;
$c_{\pm10}^5$, $c_{\pm10}^6$,\\
\hline
relations&these are equalities, not ``relations''\\
$c_{5}^ic_{10}^j=0$ and $c_{-5}^ic_{-10}^j=0$ for $i,j=5,6$&
$c_{0}^3c_{0}^5=c_{0}^8$, $c_{0}^3c_{0}^6=c_{0}^9$, $c_{0}^3
c_{0}^{11}=c_{0}^{14}$,\\
&$c_{0}^3c_{\pm5}^5=c_{\pm5}^8$, $c_{0}^3c_{\pm5}^6=c_{\pm5}^9$,\\
&$c_{0}^3c_{\pm10}^5=c_{\pm10}^8$, $c_{0}^3c_{\pm10}^6=c_{\pm10}^9$, $\dots$\\
\hline
\end{tabular}
\end{equation}

For $\mathfrak{sl}(3)$ and $p=2$, the dimensions and degrees of basis
cocycles are as follows (here, for brevity,
$kH^5$ means that $\dim H^5=k$):
\begin{alignat*}{3}%\label{sl3p=2}
& \deg = - 4\colon \quad && H^5, \ H^8,&\\
&\deg = - 2\colon \quad &&2H^3, \ 2H^5,\ 2H^8,&\\
&\deg = 0\colon \quad && H^0, \ 9H^3, \ 16H^4, \ 9H^5, \ H^8.&
\end{alignat*}
Our initial result for $\deg = 0$ was corrected in~\cite{IT}.
In the list of generators we numerate (arbitrarily)
the cocycles forming a~basis of $H_w^k$ by $c_{w}^{k,i}$, where $w$
is the degree (weight) if $\dim H_w^k>1$:
\begin{equation*}%\label{sl3genRel}
\renewcommand{\arraystretch}{1.4}
\begin{tabular}{|l|l|}
\hline generators&$c_{0}^{3,i}$ for $i=1,2$; $c_{0}^{4,i}$
for $i=1,2, 3, 4$;
$c_{0}^{5,i}$ for $i=1,2, 3$?\\
&$c_{\pm2}^{3,i}$ for $i=1,2$; $c_{\pm2}^{5,i}$ for $i=1,2$;
$c_{\pm4}^5$\\
\hline
relations&$c_{2}^Ic_{4}^5=0$ and $c_{-2}^Ic_{-4}^5=0$ for
$I=(3,i), (5,i)$, $\dots$\\
\hline
\end{tabular}
\end{equation*}

For $\mathfrak{sl}(3)$ and $p=3$, the dimensions and degrees of
basis cocycles are as follows:
\begin{gather*}
\deg =-3\colon \quad 2H^2, \ 2H^3, \ 2H^5, \ 2H^6,\nonumber\\
\deg =-0\colon \quad H^0, \ 2H^2, \ 3H^3, \ 3H^5, \ 2H^6, \ H^8.%\label{sl3p=3}
\end{gather*}

For $\mathfrak{sp}(4)\simeq\mathfrak{o}(5)$ and $p\geq5$, the dimensions and degrees of
basis cocycles are the same as for $p=0$, whereas for $p=3$, see
the following cases of its deformations.

For $\mathfrak{br}(2;\varepsilon)$, the dimensions and degrees of basis
cocycles are as follows, except for $\varepsilon=-1$ when $\mathfrak{br}(2;-1)=\mathfrak{sp}(4)$:
\begin{alignat*}{3}
&\deg = -6\colon \quad && H^3,& \nonumber\\
&\deg = -3\colon \quad && H^0,\ H^4, \ H^6, \ H^7,&\nonumber \\
&\deg = 0\colon \quad && H^0, \ 2H^3, \ H^4,\ H^6, \ H^7, \ H^{10}.&%\label{br(2,a)}
\end{alignat*}

For $\mathfrak{sp}(8)$, $\mathfrak{f}(4)$, $\mathfrak{o}(12)$,
$\mathfrak{o}(13)$, and $\mathfrak{e}(7)$ considered for $p=3$, we have $H^i(\fg)=0$
for $i=1,2$.

\subsection{Added in proof by A.~Lebedev}\label{AiP}
Remark on correspondence between 2-cochains obtained by means of \textsc{SuperLie} (in the form of elements of $M\otimes E^2(\fg^*)$) and $2$-cochains as anti-supersymmetric bilinear functions on $\fg$ with values in $M$.

I suggest to assign to $v\otimes \Pi\phi \wedge \Pi\psi$ the function
\[
f(x,y) = \big( (-1)^{p(x)(p(y)+1)}\phi(x)\psi(y) - (-1)^{p(y)}\phi(y)\psi(x) \big) v.
\]
That is, if the cocycle of $Z^2(\fg;\fg)$ is expressed in the form
\[
\sum_i c_i x_{j_i} \otimes \widehat{x_{k_i}} \wedge \widehat{x_{l_i}},
\]
where $c_i\in\Kee$ and $x$ with indexes is an element of a basis, then this cocycle corresponds to the infinitesimal deformation with the bracket
\begin{gather*}
[X,Y]_{\rm new}= [X,Y] \\
\hphantom{[X,Y]_{\rm new}=}{} + \bigg(\sum_i c_i \big( (-1)^{p(X)(p(Y)+1)}x_{k_i}^*(X)x_{l_i}^*(Y) - (-1)^{p(Y)}x_{k_i}^*(Y)x_{l_i}^*(X) \big) x_{j_i}\bigg) \hbar ,
\end{gather*}
(where $x^*$ with indexes is an element of the dual basis) and new squaring
\[
S_{\rm new}(X) = S(X) + \bigg(\sum_i c_i x_{k_i}^*(X)x_{l_i}^*(X) x_{j_i}\bigg) \hbar.
\]
(For $p\neq 2$, there is no need to describe squaring separately, but this formula will not hurt.)

Similarly, if the cocycle of $Z^2(\fg;\Kee)$ is expressed as
\[
\sum_i c_i \widehat{x_{k_i}} \wedge \widehat{x_{l_i}},
\]
then it corresponds to the central extension with the bracket
\[
[X,Y]_{\rm new} = [X,Y] + \bigg(\sum_i c_i \big( (-1)^{p(X)(p(Y)+1)}x_{k_i}^*(X)x_{l_i}^*(Y) - (-1)^{p(Y)}x_{k_i}^*(Y)x_{l_i}^*(X) \big)\bigg) z
\]
(where $z$ is the central element) and new squaring
\[
S_{\rm new} (X) = S(X) + \bigg(\sum_i c_i x_{k_i}^*(X)x_{l_i}^*(X)\bigg) z.
\]

Observe (although this is not needed in this text) that to 1-cochain
\[
v\otimes \Pi\phi \in M\otimes E^1(\fg^*)
\]
should correspond the function
\[
f(x) = (-1)^{p(x)+1} \phi(x) v.
\]

\begin{Remark}[added in proof]
David Cushing, George W.~Stagg and David I.~Stewart recently updated their article \cite{Cushing2022} with interesting results for $p=2$.
\end{Remark}

\subsection*{Acknowledgements}

We are thankful to A.~Krutov and A.~Lebedev for huge help. We are thankful to N.~Chebochko and M.~Kuznetsov for helpful
discussions of their unpublished results pertaining to this paper. We thank A.~Dzhumadildaev for pointing out \cite{IT} thus correcting an error. We are very thankful to the referees, carefully selected by SIGMA, for constructive criticism and extremely careful job.
D.L.\ is thankful to MPIMiS, Leipzig, where he was Sophus-Lie-Professor
(2004-07, when certain results of this paper were obtained), for financial
support and most creative environment. We are thankful to M.~Al Barwani, Director of the High Performance Computing resources at New York University Abu Dhabi for the possibility to perform the difficult computations of this research. S.B.~and D.L.~were supported by the grant AD~065~NYUAD.

\LastPageEnding

\end{document}